\documentclass[hidelinks,onefignum,onetabnum]{siamart251216}

\usepackage{lipsum}
\usepackage{amsfonts}
\usepackage{graphicx}
\usepackage{epstopdf}
\usepackage{algpseudocode}
\ifpdf
  \DeclareGraphicsExtensions{.eps,.pdf,.png,.jpg}
\else
  \DeclareGraphicsExtensions{.eps}
\fi

\usepackage{enumitem}
\setlist[enumerate]{leftmargin=.5in}
\setlist[itemize]{leftmargin=.5in}

\newsiamremark{remark}{Remark}
\newsiamremark{hypothesis}{Hypothesis}
\crefname{hypothesis}{Hypothesis}{Hypotheses}
\newsiamthm{claim}{Claim}
\newsiamremark{fact}{Fact}
\crefname{fact}{Fact}{Facts}

\headers{Causal Input--Output Modeling with Prior Information}{A. Taghvaei, D. Xu}

\title{Causal Input--Output Modeling with Prior Information\thanks{\funding{This work was funded by National Science Foundation (NSF) award 2318977 and Air Force Office of Scientific Research (AFOSR) award FA9550261B219.}}}

\author{Amirhossein Taghvaei\thanks{Department of Aeronautics and Astronautics, University of Washington, Seattle, 
		WA
  (\email{amirtag@uw.edu}).}
\and Daran Xu\thanks{Department of Applied Mathematics, University of Washington, Seattle,
		WA 
  (\email{daranxu@uw.edu}).}
}

\usepackage{amsopn}

\usepackage{adjustbox}
\usepackage{booktabs}
\newcommand{\Expect}{\mathbb{E}}
\newcommand{\bfU}{\mathbf{U}}
\newcommand{\bfY}{\mathbf{Y}}
\newcommand{\trace}{\operatorname{Tr}}
\newcommand{\Thetafeasible}{\Theta_{\mathrm{feasible}}}
\usepackage{comment}

\newcommand{\ud}{\mathrm{d}}

\newcommand{\indep}{\perp\!\!\!\perp}

\newsiamthm{assumption}{Assumption}

\DeclareMathOperator{\Cov}{Cov}

\newcommand{\numberthis}{\addtocounter{equation}{1}\tag{\theequation}}

\DeclareMathOperator*{\argmin}{arg\,min}

\ifpdf
\hypersetup{
  pdftitle={Causal Input--Output Modeling with Prior Information},
  pdfauthor={A. Taghvaei, D. Xu}
}
\fi

\begin{document}

\maketitle

\begin{abstract}
We study the problem of inferring a causal input--output model from unpaired marginal observations and a given prior coupling. We formulate this problem using entropic causal optimal transport, which seeks a coupling that preserves the observed marginals, satisfies a one--sided causality constraint, and remains close to the prior. For Gaussian marginals and general quadratic costs, we prove that the optimal coupling is Gaussian and reduce the problem to a finite--dimensional optimization over covariance matrices. Using a Cholesky transformation to simplify the causality constraint, we derive a Lagrangian dual formulation and develop a block coordinate descent algorithm with explicit update formulas. We establish linear convergence of the algorithm with an explicit convergence rate, prove stability of the optimal coupling under marginal covariance perturbations, and derive a statistical error bound. We further show its equivalence to Sinkhorn iterations. Numerical experiments demonstrate the computational efficiency of the proposed algorithm and illustrate the predictive advantages of causal couplings.
\end{abstract}

\begin{keywords}
Causal optimal transport, system identification, block coordinate descent
\end{keywords}

\begin{MSCcodes}
49Q22, 90C25, 65K05
\end{MSCcodes}

\section{Introduction} \label{sec:intro}

Suppose we are given unpaired input and output data from a dynamic system, together with a nominal (prior) model. The prior model provides information about the dependence between
 the input and output processes, but it may be imperfect or misspecified, so that its implied input-output distributions may differ
 from those observed in practice. Our goal is to infer an input--output model that agrees with the data while staying as close as possible to the prior. If the data were paired, this would be a standard system identification problem~\cite{ljung1999system}. Without pairing, however, the data reveal only the distribution of the inputs and the distribution of the outputs, not how the two are related. Many models are consistent with these distributions, and the prior is what allows us to choose among them.

We formalize this problem using optimal transport  theory~\cite{villani2021topics}.
Let ${U}:=\{U_t\}_{t=1}^T $ and ${Y}:=\{Y_t\}_{t=1}^T $ denote the input and output processes, with laws $\mu$ and $\nu$, respectively. An input--output model is then a joint law of $(U,Y)$, that is, a coupling $\pi\in\Pi(\mu,\nu)$, where $\Pi(\mu,\nu)$ is the set of probability measures with marginals $\mu$ and $\nu$. 
The prior model is likewise a coupling $\gamma$, which need not have marginals $\mu$ and $\nu$. Then we aim to find a coupling $\pi \in \Pi(\mu,\nu)$ that is closest to the prior $\gamma$ in the Kullback--Leibler (KL) divergence: 
\begin{equation}
\inf_{\pi\in\Pi(\mu,\nu)} D(\pi,\gamma),\qquad \text{where}\qquad D_{}(\pi, \gamma):=\begin{cases} 
\int \ln(\frac{\ud \pi}{\ud \gamma}) \ud \pi & \text{if } \pi \ll \gamma \\
\infty & \text{otherwise}
\end{cases}
\label{eq:SB-classical}
\end{equation}
This is the Schr\"odinger bridge problem~\cite{DaiPra1991,
léonard2013surveyschrodingerproblemconnections,
chen2014relationoptimaltransportschrodinger,
chen2014optimalsteeringlinearstochasticI,
caluya2021wassersteinproximalalgorithmsschrodinger} which is equivalent to the following entropy-regularized optimal transport problem when 
$\mu$, $\nu$, and $\gamma$ have strictly positive densities
with respect to the Lebesgue measure: 
\begin{equation}
\inf_{\pi\in\Pi(\mu,\nu)} \int c\,\ud\pi + D(\pi,\mu\otimes\nu),\qquad \text{where}\quad  c(\mathbf{u},\mathbf{y}) := - \ln \gamma(\mathbf{u},\mathbf{y}).
\label{eq:EOT}
\end{equation}

Problem~\eqref{eq:SB-classical}, however, ignores a basic property of input--output systems: causality. Its solution typically lets the output at time $t$ depend on inputs after time $t$. Such a model is not physically meaningful and cannot be used for sequential prediction. To rule this out, we follow~\cite{veraguas2017causaltransportdiscretetime} and call a coupling $\pi$ \emph{causal} if
\begin{equation}
(Y_1,\ldots,Y_t)\indep_{(U_1,\ldots,U_t)}(U_{t+1},\ldots,U_T),\qquad t=1,\ldots,T-1.
\label{eq:causal-constraint}
\end{equation}
In words, conditioned on the inputs up to time $t$,
the outputs up to time $t$ are independent of future inputs.
Let $\Pi_c(\mu,\nu)\subset\Pi(\mu,\nu)$ denote the set of causal couplings. Restricting~\eqref{eq:EOT} to this set yields
\begin{equation}
\inf_{\pi\in\Pi_c(\mu,\nu)} \int c\,\ud\pi + D(\pi,\mu\otimes\nu).
\label{eq:problem-formulation-2}
\end{equation}
We refer to~\eqref{eq:problem-formulation-2} as the entropic causal optimal transport (ECOT) problem.

Causal optimal transport has been studied in \cite{veraguas2017causaltransportdiscretetime}, where the causality structure is exploited to derive a dynamic programming formulation of the optimal coupling problem. However, this approach is restricted to settings in which the input process $\{U_t\}$ is Markovian and the cost function admits the separable form $\sum_{t=1}^T c_t(u_t,y_{1:t})$, which significantly limits the choice of prior models.  A closely related line of work studies \emph{bicausal} (adapted) transport, where causality is imposed in both directions~\cite{veraguas2017causaltransportdiscretetime,eckstein2023computationalmethodsadaptedoptimal,BartlBeiglbockPammer2023WassersteinSpace}. For Gaussian marginals with quadratic costs $\sum_{t=1}^T \|y_t - u_t\|^2$, closed-form solutions under the \emph{bicausal} constraint are available in both the unregularized setting ~\cite{gunasingam2025adaptedoptimaltransportgaussian,gunasingam2026adaptedoptimaltransportfiltered} and the entropically regularized setting~\cite{acciaio2025entropicadaptedwassersteindistance}. 
In continuous time, \cite{JiangLim2025Transfer} also derives explicit formulas for the adapted Wasserstein distance between
Gaussian processes.
However, a bi-causal constraint is not natural for an input-output model. On the computational side,  \cite{eckstein2023computationalmethodsadaptedoptimal} proposes abstract Sinkhorn-type iterations for ECOT.
However, each iteration involves a backward recursion over distributions or functions, and no explicit finite-dimensional implementation is available, even in the Gaussian setting.

In this paper, we address these gaps in the Gaussian setting in which the input distribution $\mu$, output distribution $\nu$, and prior coupling $\gamma$ are all Gaussian. We develop an iterative algorithm that solves the ECOT problem in this setting, and we establish explicit convergence rates and stability guarantees for it.
Our contributions are as follows:
\begin{itemize}
\item We  reduce the infinite-dimensional ECOT problem to a finite dimensional optimization problem over covariance matrices (Section~\ref{sec: optimization in cov space}). We introduce a change of coordinates, based on Cholesky factorization of marginal covariances, that normalizes the marginals and greatly simplifies the causality constraint (Section~\ref{section: cholesky Coordinates}).  

\item We derive the Lagrangian dual and establish strong duality (Theorem~\ref{theorem: dual problem of ecot gaussian}). Based on the dual, we propose a block coordinate descent (BCD) algorithm with closed-form updates (Algorithm~\ref{alg:bcd-dual}, Theorems~\ref{theorem: stationary condition update-upsi} and~\ref{theorem: stationary condition of update-lambday}). 
\item We prove that the BCD algorithm converges linearly, with an explicit rate that depends on the marginal and prior covariance (Theorem~\ref{theorem: linear convergence}). We further show that the algorithm coincides with the Gaussian specialization of the Sinkhorn iterations of~\cite{eckstein2023computationalmethodsadaptedoptimal} (Proposition~\ref{proposition: sinkhorn is bcd}). To our knowledge, this is the first linear convergence rate for causal Sinkhorn iterations with an unbounded cost.
\item We prove that the optimal coupling depends Lipschitz-continuously on the marginal covariances (Theorem~\ref{theorem: stability under marginal covariance perturbations}). As a consequence, when the marginals are estimated from $N$ samples, the estimated coupling is within $O(N^{-1/2})$ of the true one with high probability (Theorem~\ref{theorem: error bound}).
\item We validate the theory numerically. The experiments confirm
that the computed coupling satisfies the causality constraint
(Figure~\ref{fig:cross_covariance_plot}), that the BCD algorithm
converges linearly, and that the estimation error decays as the
sample size grows (Figure~\ref{fig:dual_convergence_and_objective}). In input-only prediction, the causal coupling outperforms
non-causal alternatives, especially when the prior itself is
non-causal (Figure~\ref{fig:simulation_input_only}). Finally, the BCD algorithm is also substantially faster than generic
convex solvers (Figure~\ref{fig:benchmark_comparison}).

\end{itemize}

Our work is also connected to a line of research on identifying dynamical models from distributional or aggregate data~\cite{bernstein2016consistently,haasler2021multimarginal,haasler2021control,mascherpa2023estimating,terpin2024dynamic,emerick2024causal,mascherpa2025convex}. Most of these approaches formulate the problem as a multi-marginal optimal transport problem over a graph in which each node corresponds to a marginal distribution, and the couplings between all nodes must be inferred. Our input--output setting can also be viewed as a graph, with nodes $U_s$ and $Y_\tau$ for $1 \leq s,\tau \leq T$. The key difference is the information available: we assume that the input and output sequences are given, so the dependence structure along each sequence is known. In the multi-marginal setting, by contrast, only the individual marginals are observed and the couplings between all nodes are unknown.

The remainder of this paper is organized as follows. After establishing our notation in Section~\ref{sec:notation}, Section~\ref{sec: optimization in cov space} proves the Gaussianity of the optimal coupling and introduces the Cholesky transformation. Section \ref{section: gaussian dual problem} derives the dual formulation and block coordinate descent (BCD) scheme. Section~\ref{sec: linear convergence} establishes the linear convergence rate of the BCD algorithm. Section~\ref{sec: stability} analyzes the stability of the optimal coupling under marginal covariance perturbations and derives a statistical error bound. We then establish the equivalence between the BCD scheme and Sinkhorn. Finally, Section~\ref{sec: numerical experiment} provides the numerical experiments.

{
\subsection{Notation}\label{sec:notation}
Suppose time is discretized into $T\in\mathbb{N}$ steps. 
For a Gaussian process $x_{1:T}$ with $x_t\in\mathbb R^d$,
we identify $x_{1:T}$ with the stacked Gaussian vector
$(x_1^\top,\ldots,x_T^\top)^\top\in\mathbb R^{Td}$,
whose covariance matrix has $(s,t)$-th block
$\operatorname{Cov}(x_s,x_t)\in\mathbb R^{d\times d}$.
Conditioning on $x_{s:t}$ corresponds to conditioning
on the associated subvector.


Let $\mathbb{S}^{n}$ and  $\mathbb{S}^{n}_{++}$ denote the space of $n\times n$ real symmetric matrices and the space of $n\times n$ positive definite matrices, respectively.
For a matrix \(M\), we denote its spectral norm and Frobenius norm by
\(\|M\|_2\) and \(\|M\|_F\), respectively. For any symmetric positive semidefinite matrix $Q\succeq0$,
$Q^{1/2}$ denotes its unique symmetric positive semidefinite
square root, satisfying $Q^{1/2}Q^{1/2}=Q$.
}

\section{Gaussian ECOT: Finite--dimensional Formulation}\label{sec: optimization in cov space}
We study the problem \eqref{eq:problem-formulation-2} in the setting where the cost is quadratic and the marginal distributions are Gaussian.  Under these assumptions, we establish that the optimal coupling is Gaussian, thus, reducing the original infinite-dimensional optimization over couplings to a finite-dimensional problem over covariance matrices. 
We then introduce a change of coordinates based on Cholesky factorization that normalizes the marginal covariance matrices and simplifies the causality constraints. The final simplified formulation is~\eqref{eqn:problem-over-cov} which forms the foundation for the subsequent sections. 


\subsection{Gaussian reduction}\label{section: Gaussian reduction}
Throughout the paper, we have the following assumption about the cost function and the marginal distributions: 
\begin{assumption}\label{assumption:cost} The cost function $c$ in problem \eqref{eq:problem-formulation-2} is of the quadratic form
\begin{equation}
c(\mathbf u,\mathbf y)
=
\frac12
\begin{bmatrix}
\mathbf u\\
\mathbf y
\end{bmatrix}^{\top}
(\Sigma^\gamma)^{-1}
\begin{bmatrix}
\mathbf u\\
\mathbf y
\end{bmatrix},\label{eqn: quad c with pd hessian}
\end{equation}
with a finite $2Td\times 2Td$ positive-definite matrix $\Sigma^\gamma\succ 0$.
\end{assumption}

\begin{assumption}\label{assumption:maringal} 
Both marginals $\mu$ and $\nu$ are centered non-degenerate Gaussian distributions, $\mu=\mathcal N(0,\Sigma^\mu),\quad \nu=\mathcal N(0,\Sigma^\nu), $
with finite positive definite $Td\times Td$ covariance matrices
$\Sigma^\mu,\Sigma^\nu\succ0$. 
\end{assumption}
Our results extend to
non-centered Gaussian marginals by simply shifting the coordinate origin (details are provided in~\Cref{sec: supplementary gaussian ecot}).
We now establish the main structural result for the optimizer of \eqref{eq:problem-formulation-2} under  Assumptions \ref{assumption:cost}-\ref{assumption:maringal}.
\begin{lemma}\label{theorem: optomality of gaussian}
Under Assumptions \ref{assumption:cost}-\ref{assumption:maringal}, the following hold:
\begin{enumerate}[label=\alph*.]
    \item The value of problem \eqref{eq:problem-formulation-2} is finite.
    \item Problem \eqref{eq:problem-formulation-2} admits a unique optimizer.
    \item The unique optimizer
    of \eqref{eq:problem-formulation-2} is a non--degenerate Gaussian coupling.
\end{enumerate}
\end{lemma}
The key step is to show that replacing a feasible
coupling by the Gaussian coupling with the same
mean and covariance preserves causality and does
not increase the objective.
The detailed proof is provided in
\Cref{sec: supplementary Gaussian reduction}.


By Lemma \ref{theorem: optomality of gaussian}, we know
the unique optimizer of problem (\ref{eq:problem-formulation-2}) is Gaussian. Therefore, the optimization problem (\ref{eq:problem-formulation-2}) can be reformulated as a finite-dimensional optimization over the mean vector and the covariance matrix of the Gaussian coupling.
Moreover, for every Gaussian coupling
$\pi\in\Pi_c(\mu,\nu)$ with $\pi\sim\mathcal N(m^\pi,\Sigma^\pi)$, the marginal constraints determine its mean:
$m^\pi=
\mathbf 0.$
Thus, problem \eqref{eq:problem-formulation-2} can be reduced to an optimization over the covariance matrix $\Sigma^\pi.$

\subsection{Cholesky Coordinate Transformation} \label{section: cholesky Coordinates}
We simplify the problem by a change of coordinates, through the Cholesky factorization of the marginal covariance matrices. Under Assumption \ref{assumption:maringal}, both matrices $\Sigma^\mu,\Sigma^\nu$ are positive--definite, hence  admit a Cholesky factorization 
\begin{equation}
\Sigma^\mu=L^\mu(L^\mu)^\top,
\qquad
\Sigma^\nu=L^\nu(L^\nu)^\top,
\label{eqn: marginal cholesky factors}
\end{equation}
where
\(L^\mu,L^\nu\in\mathbb R^{Td\times Td}\) are invertible block lower-triangular
matrices. Introduce the transformed variables
$\widetilde{\mathbf{U}}=(L^\mu)^{-1}{\mathbf{U}},
\quad
\widetilde {\mathbf{Y}}=(L^\nu)^{-1}{\mathbf{Y}}.$
By the Cholesky factorization, the transformed variables have
standard Gaussian distributions,
\begin{equation*}
\widetilde\bfU \sim \widetilde\mu=\mathcal N(0,I_{Td}),
\qquad
\widetilde\bfY \sim \widetilde\nu=\mathcal N(0,I_{Td}).
\end{equation*}
Define the Cholesky transformation
$T_{\mathrm{chol}}:\mathbb R^{2Td}\to\mathbb R^{2Td}$ by
\begin{equation}
T_{\mathrm{chol}}\left(
\begin{bmatrix}
\mathbf u\\
\mathbf y
\end{bmatrix}
\right)
:=
L^{-1}
\begin{bmatrix}
\mathbf u\\
\mathbf y
\end{bmatrix},
\qquad
L:=
\begin{bmatrix}
L^\mu & 0\\
0 & L^\nu
\end{bmatrix}.
\label{eqn: def of L for two marginal cholesky}
\end{equation}
The Cholesky transformation is useful because it normalizes the marginal covariance matrices and preserves and simplifies the causality.


\begin{lemma}
\label{corollary: causality under cholesky transformation}
Let $\pi$ be a non--degenerate Gaussian coupling of
$(U_{1:T},Y_{1:T})$, and define its Cholesky transformation by
\begin{equation}
\widetilde\pi:=(T_{\mathrm{chol}})_\#\pi.
\label{eqn: pushforward defining pi tilde from pi and map Tchol}
\end{equation}
Then the following hold.
\begin{enumerate}[label=\alph*.]
    \item The marginal constraints are preserved:
   $ \pi\in\Pi(\mu,\nu)
   \Longleftrightarrow
    \widetilde\pi\in\Pi(\widetilde\mu,\widetilde\nu).$

    \item The transformed coupling $\widetilde\pi$ is causal if and only if
    \begin{equation}
    \Cov_{\widetilde\pi}
    \left(
    \widetilde Y_s,\widetilde U_t
    \right)
    =0,
    \qquad
    \forall(s,t)\quad\text{s.t.}\quad 1\leq s<t\leq T.
    \label{eqn: causal constraint in cholesky coordinates}
    \end{equation}

     \item The causal constraints are preserved:
$ \pi\in\Pi_c(\mu,\nu)
    \Longleftrightarrow     \widetilde\pi\in\Pi_c(\widetilde\mu,\widetilde\nu).$

\end{enumerate}
\end{lemma}

\begin{proof}[Proof sketch]
Part (a) follows from the invertibility of
$T_{\mathrm{chol}}$ and the fact that it acts
separately on the two marginal variables.
For part (b), applying the Gaussian conditional covariance formula to $\widetilde\pi$, and using $\Sigma^{\widetilde\mu}=I_{Td}$, we obtain that $\widetilde\pi$ is causal if and only if
\eqref{eqn: causal constraint in cholesky coordinates} holds.
For part (c), let $\sigma(\cdot)$ denote the generated sigma-algebra.
Since $L^\mu$ and $L^\nu$ are invertible block lower-triangular matrices,
we have
$\sigma(U_{1:t})=\sigma(\widetilde U_{1:t})$,
$\sigma(Y_{1:t})=\sigma(\widetilde Y_{1:t})$, and
$\sigma(U_{1:T})=\sigma(\widetilde U_{1:T})$.
These imply that $\pi$ is causal if and only if
$\widetilde\pi$ is causal, since causality is equivalent to
$Y_{1:t}\perp\!\!\!\perp U_{1:T}\mid U_{1:t}$.
Together with part (a),
this proves part (c). The complete proof is provided in~\Cref{sec: supplementary cholesky Coordinates}.
\end{proof}


Under the Cholesky coordinates, the cost function is transformed into
\begin{equation}
\widetilde c(\widetilde{\mathbf u},\widetilde{\mathbf y})
:=c(L^\mu\widetilde{\mathbf u},L^\nu\widetilde{\mathbf y})
=
\frac12
\begin{bmatrix}
\widetilde{\mathbf u}\\
\widetilde{\mathbf y}
\end{bmatrix}^{\top}
(\widetilde\Sigma^\gamma)^{-1}
\begin{bmatrix}
\widetilde{\mathbf u}\\
\widetilde{\mathbf y}
\end{bmatrix},
\qquad
\widetilde\Sigma^\gamma
:=L^{-1}\Sigma^\gamma L^{-\top}\succ0.
\label{eqn: quadratic c tilde in cholesky}
\end{equation}
The transformed covariance matrix $\widetilde\Sigma^\gamma$
is positive definite because $\Sigma^\gamma\succ0$ and $L$ is
invertible. Thus, the cost function $\widetilde c$ remains
quadratic and satisfies Assumption~\ref{assumption:cost}.
Moreover, the entropy component of the objective of problem \eqref{eq:problem-formulation-2} is unchanged under the Cholesky transformation. 
Therefore, it is
sufficient to study the ECOT problem with
$\Sigma^\mu=\Sigma^\nu=I_{Td}.$
Once the optimizer $\widetilde\pi^*$ in the Cholesky
coordinates is obtained, the optimizer in the original
coordinates is recovered through
$\pi^*=(T_{\mathrm{chol}}^{-1})_\#\widetilde\pi^*$.

\subsection{Simplified problem formulation} \label{section: gaussian primal and dual}
We now derive the finite-dimensional covariance formulation of problem
\eqref{eq:problem-formulation-2} in the Cholesky coordinates.
Throughout the remainder of this section, we work in these coordinates and
drop the tildes for notational simplicity.
Thus,
\[
\Sigma^\mu=\Sigma^\nu=I_{Td},
\]
and $\Sigma^\gamma$ denotes the transformed covariance
$\widetilde{\Sigma}^\gamma$ defined in \eqref{eqn: quadratic c tilde in cholesky}. 

\noindent
{\bf Objective function:}
Let $\pi=\mathcal N(0,\Sigma)$ be a non--degenerate Gaussian coupling.
Under Assumption~\ref{assumption:maringal}, together with
$\Sigma^\mu=\Sigma^\nu=I_{Td}$, the marginals are
$\mu=\nu=\mathcal N(0,I_{Td})$, and hence
$\mu\otimes\nu=\mathcal N(0,I_{2Td})$.
Therefore, up to an additive constant independent of $\Sigma$,
the objective function of problem \eqref{eq:problem-formulation-2}
takes the form
\begin{equation}
\mathcal J(\Sigma)
:=
\frac12\operatorname{Tr}(Q\Sigma)
-\frac12\log\det\Sigma,
\qquad
Q:=(\Sigma^\gamma)^{-1}+I_{2Td}\succ0,
\label{eqn:primal loss}
\end{equation}
for $\Sigma\succ0$. We set $\mathcal J(\Sigma)=+\infty$
when $\Sigma$ is not positive definite.

\medskip
\noindent
{\bf Constraints.} 
Define the projection matrices
\[
P_u=
\begin{bmatrix}
I_{Td} & 0
\end{bmatrix}
\in\mathbb{R}^{Td\times 2Td},
\qquad
P_y=
\begin{bmatrix}
0 & I_{Td}
\end{bmatrix}
\in\mathbb{R}^{Td\times 2Td}.
\]
For any matrix \(\Sigma\in\mathbb{R}^{2Td\times 2Td}\),
\(P_u\Sigma P_u^\top\) extracts the upper-left \(Td\times Td\) block of
\(\Sigma\), while \(P_y\Sigma P_y^\top\) extracts the lower-right
\(Td\times Td\) block. Using the projection matrices, 
the marginal constraints are expressed as $P_u\Sigma P_u^\top=I_{Td},$ and $P_y\Sigma P_y^\top=I_{Td}.$
In order to express the causality constraints, we first note that Lemma \ref{corollary: causality under cholesky transformation}, specifically
\eqref{eqn: causal constraint in cholesky coordinates}, is equivalent to $\Sigma_{Y_{1:t},U_{t+1}}=0,$
 for $t=1,\ldots,T-1.$
In order to specify the constraint more compactly, introduce the projection matrices
\[
E_{t+1}
=
\begin{bmatrix}
0_{td\times d}\\
I_d\\
0_{(T-t-1)d\times d}
\end{bmatrix}
\in\mathbb R^{Td\times d},
\qquad
S_t
=
\begin{bmatrix}
I_{td} & 0_{td\times (T-t)d}
\end{bmatrix}
\in\mathbb R^{td\times Td}.
\]
for
\(t=1,\ldots,T-1\). The matrix \(S_t\) selects the first \(t\) \(d\)-dimensional blocks,
whereas \(E_{t+1}\) selects the \((t+1)\)-st \(d\)-dimensional block. Using these projection matrices, the causality constraint is expressed as  $ S_t P_y \Sigma P_u^\top E_{t+1} =0,$ for  $t=1,\ldots,T-1.$

\medskip
\noindent
{\bf New formulation.} 
Combining the objective and constraint representations above, problem
\eqref{eq:problem-formulation-2} reduces to
\begin{subequations}
\label{eqn:problem-over-cov}
\begin{align}
\min_{\Sigma\in\mathbb{S}^{2Td}}
\quad&
\mathcal{J}(\Sigma)
\notag\\
\text{subject to}\qquad
&P_u\Sigma P_u^\top=I_{Td},
\label{eqn:maringal-constraint-mu}\\
&P_y\Sigma P_y^\top=I_{Td},
\label{eqn:maringal-constraint-nu}\\
&S_t P_y \Sigma P_u^\top E_{t+1} =0,\qquad
t=1,\ldots,T-1.
\label{eqn:causal-constraint}
\end{align}
\end{subequations}

\section{Dual formulation and the proposed algorithm}\label{section: gaussian dual problem}
We now derive the Lagrangian dual of the covariance optimization
problem \eqref{eqn:problem-over-cov}, establish strong duality, and present the iterative block-coordinate-descent algorithm to solve it.

\subsection{Dual formulation}\label{section: dual formulation}
Define the dual space
\[
\Theta
=
\mathbb S^{Td}\times\mathbb S^{Td}
\times
\prod_{t=1}^{T-1}\mathbb R^{td\times d},
\]
and dual variables $\theta=(\Lambda_u,\Lambda_y,\Psi)\in\Theta,$ $\Psi=\{\Psi_t\}_{t=1}^{T-1},$ $\Psi_t\in\mathbb R^{td\times d}.$
The variables $\Lambda_u$, $\Lambda_y$, and $\Psi$ are dual to the constraints~\eqref{eqn:maringal-constraint-mu},~\eqref{eqn:maringal-constraint-nu}, and~\eqref{eqn:causal-constraint}, respectively. 
Define
\begin{equation}
\label{eqn: Z from theta}
Z(\theta)
:=
Q+
\begin{bmatrix}
2\Lambda_u & \Gamma_\Psi\\
\Gamma_\Psi^\top & 2\Lambda_y
\end{bmatrix},
\quad \text{with } 
 \Gamma_\Psi
:=
\sum_{t=1}^{T-1}
E_{t+1}\Psi_t^\top S_t
\in\mathbb R^{Td\times Td}.
\end{equation}
By construction, \(\Gamma_\Psi\) is strictly block
lower triangular, with \(\Psi_t^\top\) occupying
the first \(t\) block columns of the
\((t+1)\)-st block row.
The map $Z:\Theta \to \mathbb S^{2Td}$  connects the dual variables with the primal
covariance matrix in our duality result. Finally, define the dual objective
\begin{equation}
\label{eqn: dual objective}
\mathcal D(\theta)
:=
\begin{cases}
\displaystyle
\frac12\log\det Z(\theta)
-
\operatorname{Tr}(\Lambda_u)
-
\operatorname{Tr}(\Lambda_y)
+
Td,
& \theta \in \Thetafeasible,\\[2mm]
-\infty,
& \text{otherwise}.
\end{cases}
\end{equation}
 where 
$
    \Thetafeasible
:=
\{\theta\in\Theta:Z(\theta)\succ0\}.
$
$\Thetafeasible$ is convex, since $Z$ is affine and the positive-definite cone $\mathbb S_{++}^{2Td}$
is convex.
Consider the following maximization problem:
\begin{equation}
\label{eqn: dual problem}
\max_{\theta\in\Theta}
\mathcal D(\theta).
\end{equation}


\begin{theorem}
\label{theorem: dual problem of ecot gaussian}
The
following hold:
\begin{enumerate}[label=\alph*.]
\item Problem \eqref{eqn: dual problem} is the Lagrangian dual of
the covariance optimization problem \eqref{eqn:problem-over-cov}.
Strong duality holds, i.e.,
\[
\min_{\substack{\Sigma\in\mathbb S^{2Td}\\
\text{satisfying }\eqref{eqn:maringal-constraint-mu},
\eqref{eqn:maringal-constraint-nu},
\eqref{eqn:causal-constraint}}}
\mathcal J(\Sigma)
=
\max_{\theta\in\Theta}\mathcal D(\theta),
\]
and the dual optimum is attained.
    \item The dual objective \(\mathcal D\) is strictly concave on
    \(\Thetafeasible:=
 \{\theta\in\Theta:Z(\theta)\succ0\}\).

    \item The dual optimizer \(\theta^*\) is unique.

    \item The primal and dual optimizers satisfy
    \(
    \Sigma^*=Z(\theta^*)^{-1}.
    \)
\end{enumerate}
\end{theorem}
The proof follows from Lagrangian duality and
Slater's condition. The strict concavity of the
dual objective follows from the injectivity of
$Z$ and the strict concavity of $\log\det$.
The detailed proof is provided in
\Cref{sec: supp for dual formulation}.
 
\subsection{Block Coordinate Descent (BCD) Algorithm} \label{section:bcd}
We present an iterative algorithm to solve the dual problem~\eqref{eqn: dual problem}. 
Throughout the rest of the paper, we adopt the convention $f(\theta):=-\mathcal{D}(\theta)$
in order to express the dual problem \eqref{eqn: dual problem}   as the minimization
\[
\min_{\theta\in\Theta} f(\theta).
\]
\begin{algorithm}[t]
\caption{block coordinate descent for the Dual Problem}
\label{alg:bcd-dual}
\begin{algorithmic}[1]
\State \textbf{Input:} Stopping tolerance $\varepsilon_{\mathrm{tol}}>0$

\State Initialize $\theta^{\{0\}} = (\Lambda_u^{\{0\}}, \Lambda_y^{\{0\}}, \Psi^{\{0\}}) = 0$ \Comment{component-wise zero}

\For{$k = 0,1,2,\dots$} 
     \State Set $ \theta^{\{k\}} =
    (\Lambda_u^{\{k\}}, \Lambda_y^{\{k\}}, \Psi^{\{k\}})$
   \State \textbf{Update $(\Lambda_u,\Psi)$ using
\eqref{eqn: schrodinger eqn 3} and
\eqref{eqn: schrodinger eqn 1}:}
\begin{equation}
\label{eq:update-upsi}
(\Lambda_u^{\{k+1\}},\Psi^{\{k+1\}})
=
\arg\min_{(\Lambda_u,\Psi)}
f(\Lambda_u,\Lambda_y^{\{k\}},\Psi)
\end{equation}


    \State \textbf{Update $\Lambda_y$ using
\eqref{eqn: schrodinger eqn 2}:}
\begin{equation}
\label{eq:update-lambday}
\Lambda_y^{\{k+1\}}
=
\arg\min_{\Lambda_y}
f(\Lambda_u^{\{k+1\}},\Lambda_y,\Psi^{\{k+1\}})
\end{equation}

    \State Set $ \theta^{\{k+1\}} =
    (\Lambda_u^{\{k+1\}}, \Lambda_y^{\{k+1\}}, \Psi^{\{k+1\}})$
\If{$\|\theta^{\{k+1\}}-\theta^{\{k\}}\|_\Theta
<\varepsilon_{\mathrm{tol}}$}
    \State \textbf{break}
\EndIf
\EndFor
\State \textbf{Return:} $\theta^{\{k+1\}}$
\end{algorithmic}
\end{algorithm}

To solve this problem computationally, we apply a block coordinate descent scheme, decomposing the full dual variable $\theta$ into distinct coordinate blocks and minimizing the objective $f$ sequentially over each subset,  as summarized in Algorithm \ref{alg:bcd-dual}.
This iterative method is initialized at
\[
\theta^{\{0\}}=\mathbf 0.
\]
which is feasible, i.e. 
\(\theta^{\{0\}}\in \Thetafeasible\), because
\(
Z(\theta^{\{0\}})
=
Q
\succ0.
\) Then, the algorithm proceeds by solving the two subproblems~\eqref{eq:update-upsi} and~\eqref{eq:update-lambday}. Each of the subproblems is strictly convex, because, by Theorem~\ref{theorem: dual problem of ecot gaussian}, $f$ is strictly convex on $\Thetafeasible$, so its restriction to each feasible block-coordinate slice is strictly convex. 
The following propositions provide explicit solutions to
the subproblems by deriving the stationary conditions
and showing that they uniquely determine feasible block updates.
 To state the result, we introduce the following notation:
For any matrix $M\in\mathbb{R}^{2Td\times2Td}$, we use the block partition
\[
M=\begin{bmatrix}
M_{uu} & M_{uy}\\
M_{yu} & M_{yy}
\end{bmatrix},
\qquad
M_{uu},M_{uy},M_{yu},M_{yy}\in\mathbb{R}^{Td\times Td}.
\]

\begin{theorem}
\label{theorem: stationary condition update-upsi}
Suppose that
$\theta^{\{k\}}
\in\Thetafeasible.$
Then the block subproblem \eqref{eq:update-upsi} admits a unique minimizer.
Moreover, the following hold.
\begin{enumerate}[label=\alph*.]
    \item For \(t=1,\ldots,T-1\),
    \begin{equation}
    \Psi_t^{\{k+1\}}
    =
    -
    \left(
    S_tZ_{yy}^{-1}S_t^\top
    \right)^{-1}
    S_tZ_{yy}^{-1}Q_{yu}E_{t+1},
    \label{eqn: schrodinger eqn 3}
    \end{equation}
    where
    \[
    Z_{yy}
    =
    Q_{yy}+2\Lambda_y^{\{k\}}.
    \]

    \item After \(\Psi^{\{k+1\}}\) is obtained,
    \begin{equation}
    2\Lambda_u^{\{k+1\}}
    =
    -Q_{uu}
    +
    I_{Td}
    +
    Z_{uy}Z_{yy}^{-1}Z_{uy}^\top,
    \label{eqn: schrodinger eqn 1}
    \end{equation}
    where
    \[
    Z_{uy}
    =
    Q_{uy}+\Gamma_{\Psi^{\{k+1\}}}.
    \]

    \item The updated parameters satisfy
    \[
    \theta^{\{k+\frac12\}}
    =
    (\Lambda_u^{\{k+1\}},
    \Lambda_y^{\{k\}},
    \Psi^{\{k+1\}})
    \in\Thetafeasible.
    \]
\end{enumerate}
\end{theorem}
The update formulas follow from the stationarity conditions and the Schur complement identities. 
The detailed proof is provided in
\Cref{section: supp bcd}.


In \eqref{eqn: schrodinger eqn 3}, all \(\{\Psi_t^{\{k+1\}}\}_{t=1}^{T-1}\) can be
computed independently. The resulting \(\Psi^{\{k+1\}}\) is then used in
\eqref{eqn: schrodinger eqn 1} to determine
\(\Lambda_u^{\{k+1\}}\).  Let $\theta^{\{k+\frac12\}}
:=
(\Lambda_u^{\{k+1\}},
\Lambda_y^{\{k\}},
\Psi^{\{k+1\}})$
\begin{theorem}
\label{theorem: stationary condition of update-lambday}
Suppose that
$\theta^{\{k+\frac12\}}
\in\Thetafeasible.$
Then the block subproblem \eqref{eq:update-lambday}
admits a unique minimizer. Moreover, the following hold.
\begin{enumerate}[label=\alph*.]
    \item The minimizer is given by
    \begin{equation}
    2\Lambda_y^{\{k+1\}}
    =
    -Q_{yy}
    +
    I_{Td}
    +
    Z_{yu}Z_{uu}^{-1}Z_{yu}^\top,
    \label{eqn: schrodinger eqn 2}
    \end{equation}
    where
    \[
    Z_{uu}
    =
    Q_{uu}+2\Lambda_u^{\{k+1\}},
    \qquad
    Z_{yu}
    =
    Q_{yu}+\Gamma_{\Psi^{\{k+1\}}}^{\top}.
    \]

    \item The updated parameters satisfy
    \[
    \theta^{\{k+1\}}
    =
    (\Lambda_u^{\{k+1\}},
    \Lambda_y^{\{k+1\}},
    \Psi^{\{k+1\}})
    \in\Thetafeasible.
    \]
\end{enumerate}
\end{theorem}
The proof is provided in~\Cref{section: supp bcd}.
Algorithm~\ref{alg:bcd-dual} can be implemented by iteratively applying
Theorems~\ref{theorem: stationary condition update-upsi} and
\ref{theorem: stationary condition of update-lambday}.
\section{Convergence Analysis}
\label{sec: linear convergence}
We now establish the convergence properties of the BCD scheme in
Section~\ref{section:bcd}. 
The function \(f\) is \(C^\infty\) on
\(\Thetafeasible\), since \(Z(\theta)\) is affine in \(\theta\) and
\(\log\det(\cdot)\) is \(C^\infty\) on
\(\mathbb S_{++}^{2Td}\). Hence, all derivatives used below are
well--defined and continuous.

\subsection{Smoothness and strong convexity over compact domains}\label{section: smooth and str convx over compact} We establish the smoothness and strong convexity of the function $f$ over compact subsets of $\Thetafeasible$. In order to define these properties, 
we  equip \(\Theta\) with an inner product and its induced norm, which will be used throughout the convergence
analysis. For
$
\theta=(\Lambda_u,\Lambda_y,\Psi) \in \Theta$  and $
\bar\theta=(\bar\Lambda_u,\bar\Lambda_y,\bar\Psi)
\in\Theta
$, define the inner-product  using the Frobenius inner product on each matrix component:
\begin{equation}
    \label{eq: def of inner produce of Theta space}
\langle\theta,\bar\theta\rangle_{\Theta}
:=
\langle\Lambda_u,\bar\Lambda_u\rangle_F
+
\langle\Lambda_y,\bar\Lambda_y\rangle_F
+
\sum_{t=1}^{T-1}
\langle\Psi_t,\bar\Psi_t\rangle_F.
\end{equation}
and denote the 
 induced norm by \(\|\theta\|_\Theta:=\sqrt{ \langle\theta,\theta\rangle_\Theta}\). 

For the convergence analysis, it is convenient to introduce a notation for the linear part
of the affine map \(Z(\theta)\):
\begin{equation}
\label{eqn: A from theta}
A(\theta)
:=
Z(\theta)-Z(0)
=
\begin{bmatrix}
2\Lambda_u & \Gamma_\Psi\\
\Gamma_\Psi^\top & 2\Lambda_y
\end{bmatrix}.
\end{equation}
Since $Z$ is affine, for any $\theta,\delta\theta\in\Theta$,
we have $Z(\theta+\delta\theta)-Z(\theta)=A(\delta\theta)$.
Moreover, $A$ is injective.
We use $A$ to express the Hessian of the dual objective.

\begin{lemma}
\label{lemma: hessian of f}
For every \(\theta\in\Thetafeasible\) and
\(\delta\theta_1,\delta\theta_2\in\Theta\),
\begin{equation}
\nabla^2 f(\theta)[\delta\theta_1,\delta\theta_2]
=
\frac12
\left\langle
Z(\theta)^{-1/2}A(\delta\theta_1)Z(\theta)^{-1/2},
Z(\theta)^{-1/2}A(\delta\theta_2)Z(\theta)^{-1/2}
\right\rangle_F.
\label{eqn: hessian of f}
\end{equation}
In particular,
\begin{equation}
\nabla^2 f(\theta)[\delta\theta,\delta\theta]=
\frac12
\left\|
Z(\theta)^{-1/2}A(\delta\theta)Z(\theta)^{-1/2}
\right\|_F^2>0
\qquad
\text{for every }\delta\theta\neq0.
\label{eqn: strictly positive hessian of f}
\end{equation}
\end{lemma}
To describe the two blocks used in the BCD scheme, we decompose the domain $\Theta = \Theta_1 \times \Theta_2$ where 
\begin{equation*}
    \Theta_1
:=
\mathbb S^{Td}
\times
\prod_{t=1}^{T-1}\mathbb R^{td\times d},
\qquad
\Theta_2
:=
\mathbb S^{Td}.
\end{equation*}
and for
\(\theta=(\Lambda_u,\Lambda_y,\Psi)=(\theta_1,\theta_2) \in \Theta_1 \times\Theta_2\), we denote the coordinates in the subdomain by $\theta_1:=(\Lambda_u,\Psi)\in\Theta_1,
$ $
\theta_2:=\Lambda_y\in\Theta_2.$
The spaces \(\Theta_1\) and \(\Theta_2\) inherit their
inner products and norms from \(\Theta\), so that $\|\theta\|_\Theta^2
=
\|\theta_1\|_{\Theta_1}^2
+
\|\theta_2\|_{\Theta_2}^2.$
We write \(f(\theta)=f(\theta_1,\theta_2)\) and denote
the corresponding block gradients by
\(\nabla_{\theta_1}f\) and \(\nabla_{\theta_2}f\):
\(\nabla f(\theta_1,\theta_2)
=
\bigl(
\nabla_{\theta_1}f(\theta_1,\theta_2),
\nabla_{\theta_2}f(\theta_1,\theta_2)
\bigr).\)



We next introduce two regularity properties used in the convergence
analysis. Let \(S\subset\Thetafeasible\), and we say that the gradient of \(f\)
is block--coordinate Lipschitz continuous on \(S\) with constants
\(L_1,L_2>0\) if 
\begin{align}
\left\|
\nabla_{\theta_1}f(\widetilde\theta_1,\theta_2)
-
\nabla_{\theta_1}f(\theta_1,\theta_2)
\right\|_{\Theta_1}
&\leq
L_1
\|\widetilde\theta_1-\theta_1\|_{\Theta_1},
\quad
\forall
(\theta_1,\theta_2),
(\widetilde\theta_1,\theta_2)\in S,
\label{eqn: block coordinate Lipschitz theta1}
\\
\left\|
\nabla_{\theta_2}f(\theta_1,\widetilde\theta_2)
-
\nabla_{\theta_2}f(\theta_1,\theta_2)
\right\|_{\Theta_2}
&\leq
L_2
\|\widetilde\theta_2-\theta_2\|_{\Theta_2},
\quad
\forall
(\theta_1,\theta_2),
(\theta_1,\widetilde\theta_2)\in S.
\label{eqn: block coordinate Lipschitz theta2}
\end{align}
Let \(S\subset\Thetafeasible\) be convex. We say that \(f\) is strongly convex with constant \(\sigma>0\) on \(S\) if, for any
\(\theta,\widetilde\theta\in S\),
\begin{equation}
f(\widetilde\theta)
\geq
f(\theta)
+
\left\langle
\nabla f(\theta),
\widetilde\theta-\theta
\right\rangle_{\Theta}
+
\frac{\sigma}{2}
\|\widetilde\theta-\theta\|_{\Theta}^2.
\label{eqn: definition of strong convexity}
\end{equation}
For a nonempty set  \(K\subset\Thetafeasible\), define
\begin{equation}
\label{eqn: def of min and max eigven of Z over a compact set}
m_K
:=
\inf_{\theta\in K}
\lambda_{\min}(Z(\theta)),
\qquad
M_K
:=
\sup_{\theta\in K}
\lambda_{\max}(Z(\theta)).
\end{equation}
Here, $\lambda_{\min}$ and $\lambda_{\max}$ denote the
smallest and largest eigenvalues of a positive-definite matrix, respectively.
If $K$ is compact, continuity of $Z$ and of the
extreme eigenvalues ensures that the infimum and
supremum are attained. Since $Z(\theta)\succ0$
for every $\theta\in K$, it follows that
$0<m_K\leq M_K<\infty$.
The following lemma establishes the  block–coordinate Lipschitz and strong convexity 
on compact convex subsets $K$ with explicit constants that depend on the smallest and largest eigenvalues~\eqref{eqn: def of min and max eigven of Z over a compact set}.

\begin{lemma}
\label{lemma: local lip and str convex}
Let \(K\subset\Thetafeasible\) be nonempty, compact, and convex.
Then:
\begin{enumerate}[label=\alph*.]
    \item The gradient of \(f\) is block--coordinate Lipschitz continuous
    on \(K\) with constants
    \[
    L_1=L_2
=
2\,m_K^{-2}.
    \]

    \item The function \(f\) is strongly convex on \(K\) with constant
    \[\sigma_K
=M_K^{-2}.
    \]
\end{enumerate}
\end{lemma}
The proof uses the Hessian in
Lemma~\ref{lemma: hessian of f}, together with the spectral bounds on $Z(\theta)$.
These yield the stated block Lipschitz and
strong convexity constants.
The detailed proof is provided in
\Cref{section: supp smooth and str convx over compact}.


\subsection{Sublevel Set restriction} \label{section: sublevel set restriction}
We restrict the convergence analysis to the initial sublevel set
\begin{equation*}
L(\theta^{\{0\}})
:=
\left\{
\theta\in\Thetafeasible:
f(\theta)\leq f(\theta^{\{0\}})
\right\},
\end{equation*}
which is nonempty since \(\theta^{\{0\}}=\mathbf 0\in\Thetafeasible\).
We first establish uniform spectral bounds for $Z(\theta)$
on this set, which will be used to obtain explicit regularity
constants and prove compactness.
Define
\begin{equation*}
\begin{aligned}
C_Q
&:=
\operatorname{Tr}(Q)-\log\det Q-(2Td-1),\\
a_Q
&:=-W_0(-e^{-C_Q}),\qquad
b_Q:=-W_{-1}(-e^{-C_Q}),
\end{aligned}
\end{equation*}
where $W_0$ and $W_{-1}$ are the two real branches
of the Lambert $W$ function
\cite{corless1996lambert}\footnote{The Lambert $W$ function is defined by the solutions $w=W(z)$ to the equation $we^w=z$.}. The constants $a_Q$ and $b_Q$ provide lower and upper bounds for the minimum and maximum eigenvalues in~\eqref{eqn: def of min and max eigven of Z over a compact set} when the domain is the initial sublevel set. 
\begin{lemma}
\label{lemma: uniform spectral bound on level set}
The constants $a_Q$ and $b_Q$ are well-defined and satisfy $    0<a_Q<1<b_Q<\infty.$
 Moreover, for
$\theta^{\{0\}}=0$, $a_Q
\leq
m_{L(\theta^{\{0\}})}
\leq
M_{L(\theta^{\{0\}})}
\leq
b_Q.$
\end{lemma}
\begin{proof}[Proof sketch]
By definition of $Q$, we have $-e^{-C_Q} \in (-e^{-1},0)$,
so $a_Q$ and $b_Q$ are well-defined.
For $\theta\in L(\theta^{\{0\}})$, the sublevel
condition and the expression for $f$ yield
\begin{equation*}
\sum_{i=1}^{2Td}
\left(
\lambda_i(Z(\theta))
-
\log\lambda_i(Z(\theta))
\right)
\leq
\operatorname{Tr}(Q)-\log\det Q.
\end{equation*}
Since $x-\log x\geq1$ for $x>0$, each eigenvalue
$\lambda_i(Z(\theta))$ satisfies
\begin{equation*}
\lambda_i(Z(\theta))
-
\log\lambda_i(Z(\theta))
\leq C_Q.
\end{equation*}
The equation $x-\log x=C_Q$ has two positive
solutions, $a_Q$ and $b_Q$, expressed through
the two real branches of the Lambert $W$ function.
Hence, every eigenvalue of $Z(\theta)$
lies in $[a_Q,b_Q]$.
The complete proof is in
 \Cref{section: supp sublevel set restriction}.
\end{proof}
\begin{lemma} \label{lemma: compact level set}
$L(\theta^{\{0\}})$ is compact and convex.
\end{lemma}


Moreover, the exact BCD updates
satisfy
$f(\theta^{\{k+1\}})
\leq f(\theta^{\{k+\frac12\}})
\leq f(\theta^{\{k\}})$
for every $k\geq0$.
Consequently, all BCD iterates remain in
$L(\theta^{\{0\}})$.
These properties are established in~\Cref{section: supp sublevel set restriction}.
The following corollary applies Lemma~\ref{lemma: local lip and str convex}
to \(L(\theta^{\{0\}})\) and gives explicit constants in terms of \(a_Q\) and \(b_Q\).
\begin{corollary}
\label{corollary: lip and strong convexity on sublevel set}
\begin{enumerate}[label=\alph*.]
    \item The gradient of \(f\) is block--coordinate Lipschitz continuous on
\(L(\theta^{\{0\}})\) with constants $L^Q_1,L^Q_2$, where $L_1^Q=L_2^Q
=
2\,a_Q^{-2}=:L^Q.$
\item \(f\) is strongly convex on \(L(\theta^{\{0\}})\) with constant
$\sigma^Q
=
b_Q^{-2}.$
\end{enumerate}
\end{corollary}

\subsection{Linear Convergence}
\label{subsection: linear convergence final rate}
We establish linear convergence of the BCD algorithm by adapting the argument of \cite[Theorem~5.2]{BeckAmirConvergenceofBlockCoordinateDescent}.
Since strong convexity and block--coordinate Lipschitz continuity hold only on the initial sublevel set, we must additionally verify that the block gradient steps used in the descent estimates remain in this set. The technical results supporting this argument are established in~\Cref{subsec: linear convergence final rate}.
\begin{theorem} \label{theorem: linear convergence} 
Let
$\left\{
\theta^{\{k\}},
\theta^{\{k+\frac12\}}
\right\}_{k\geq0}$
be the sequence of dual iterates generated by
Algorithm~\ref{alg:bcd-dual}. Then,
\begin{align}
f(\theta^{\{k+1\}})-f(\theta^*)
&\leq
(\rho^Q)^{2k+1}
\left[
f(\theta^{\{0\}})-f(\theta^*)
\right],
\qquad k\geq0,
\label{eqn: linear conv bound integer iterate}\\
f(\theta^{\{k+\frac12\}})-f(\theta^*)
&\leq
(\rho^Q)^{2k}
\left[
f(\theta^{\{0\}})-f(\theta^*)
\right],
\qquad k\geq0.
\label{eqn: linear conv bound half iterate}
\end{align}
where the linear convergence rate is given by $\rho^Q
:=
1-\frac{\sigma^Q}{L^Q}.$
Moreover, $\rho^Q\in(\frac12,1).$
\end{theorem}
The main step is to show that the block gradient
steps remain in the initial sublevel set
$L(\theta^{\{0\}})$, where the block Lipschitz and strong convexity bounds hold.
The linear convergence estimates then follow from the block coordinate descent analysis of \cite{BeckAmirConvergenceofBlockCoordinateDescent}.
The detailed proof is provided in
\Cref{subsec: linear convergence final rate}.

The following corollary translates the linear convergence of the dual
objective into linear convergence of the dual iterates.
\begin{corollary}
\label{corollary: linear convergence of the dual iterations}
Let \(\left\{
\theta^{\{k\}},
\theta^{\{k+\frac12\}}
\right\}_{k\geq0}\)
be the sequence of dual iterates generated by
Algorithm~\ref{alg:bcd-dual}. Then there exists a constant
\(C_0>0\) such that, for every \(k\geq0\),
\begin{equation*}
\left\|
\theta^{\{k\}}-\theta^*
\right\|_\Theta^2
\leq
C_0\rho^k, \qquad \left\|
\theta^{\{k+\frac12\}}-\theta^*
\right\|_\Theta^2
\leq
C_0\rho^k,
\end{equation*}
where $\rho:=(\rho^Q)^2\in(\frac14,1).$
\end{corollary}
This follows directly from the strong convexity of $f$
on $L(\theta^{\{0\}})$ and
Theorem~\ref{theorem: linear convergence}.

\begin{remark}
\label{remark: dependence of convergence rate on reference covariance}
We briefly discuss how the convergence-rate bound depends on the reference
 covariance \(\Sigma^\gamma\). 
Let \(s_1,\ldots,s_{2Td}>0\)
be the eigenvalues of \(\Sigma^\gamma\). Then,  by definition,
\(C_Q=
1+
\sum_{i=1}^{2Td}
\left(
\frac{1}{s_i}
-
\log\left(1+\frac{1}{s_i}\right)
\right).\)
Decreasing any \(s_i\) while keeping the others fixed increases 
\(C_Q\). An increase in \(C_Q\) moves the lower-bound \(a_Q\) toward zero and
 the upper-bound \(b_Q\) toward infinity.  Thus, the convergence factor
$\rho^Q=1-\frac12(a_Q/b_Q)^2$ increases toward $1$,
yielding a weaker convergence guarantee.

In our ECOT formulation, decreasing an eigenvalue of
\(\Sigma^\gamma\) increases the corresponding eigenvalue of
\((\Sigma^\gamma)^{-1}\), and hence strengthens the quadratic cost in
that direction while leaving the coefficient of the entropy term
unchanged. Thus, the entropy regularization becomes relatively weaker
compared with the quadratic cost in that direction. Intuitively, with
weaker relative entropy regularization, the problem moves toward a more
nearly degenerate regime.
 This provides an interpretation of the weakening of the theoretical convergence guarantee
as $\rho^Q$ approaches $1$.
\end{remark}

\section{Stability}
\label{sec: stability}
We now return to the original coordinates in order to study perturbations
of the marginal covariance matrices. Recall that in
Section~\ref{section: gaussian primal and dual}, after applying the
Cholesky transformation, we suppressed the tildes for notational
simplicity. In this section, we restore this notation to distinguish between the original and Cholesky coordinates.
We first establish a stability result under perturbations of the marginal covariance matrices. The stability analysis relies on perturbation bounds for Cholesky factors and local Lipschitz continuity
of the dual optimizer with respect to the Cholesky factors.
These technical results are established in \Cref{sec: supp of stabiility}.
Let
\begin{equation*}
    \Sigma_0
:=
\begin{bmatrix}
\Sigma^\mu&0\\
0&\Sigma^\nu
\end{bmatrix},
\quad
\Delta\Sigma_0
:=
\begin{bmatrix}
\Delta\Sigma^\mu&0\\
0&\Delta\Sigma^\nu
\end{bmatrix},
\end{equation*}
where
\(\Sigma^\mu,\Sigma^\nu\) satisfy
Assumption~\ref{assumption:maringal}, and
\(\Delta\Sigma^\mu,\Delta\Sigma^\nu\in\mathbb S^{Td}\).
Let \(\pi^*\) denote the Gaussian ECOT optimizer corresponding to
\((\Sigma^\mu,\Sigma^\nu)\), as characterized in Problem \eqref{eqn:problem-over-cov} in the Cholesky coordinates defined by \eqref{eqn: marginal cholesky factors} and \eqref{eqn: quadratic c tilde in cholesky}. 
\begin{theorem}
\label{theorem: stability under marginal covariance perturbations}
There exist constants $r_\Sigma,C_\Sigma>0$, depending only on
$\Sigma_0$, $\Sigma^\gamma$, $T$, and $d$, such that, 
\begin{equation}
\left\|
\Sigma^{\pi_\Delta^*}
-
\Sigma^{\pi^*}
\right\|_F
\leq
C_\Sigma
\|\Delta\Sigma_0\|_F,\quad \text{if}\quad \|\Delta\Sigma_0\|_F\leq r_\Sigma,
\label{eqn: stability of optimal covariance under marginal perturbation}
\end{equation}
where $\pi_\Delta^*$ denotes the optimizer corresponding to the perturbed marginals. 

\end{theorem}
The proof first controls perturbations of the marginal Cholesky factors and then applies the implicit function theorem to establish local Lipschitz continuity of the optimal
dual variables. Perturbation bounds for the inverse matrix $Z(\theta^*(L);L)^{-1}$, followed by transformation to the original coordinates, yield the stated covariance stability estimate.
The detailed proof is provided in
\Cref{sec: supp of stabiility}.

Under Assumption~\ref{assumption:maringal}, let
 $\mu=\mathcal N(0,\Sigma^\mu)$, $\nu=\mathcal N(0,\Sigma^\nu)$,
and suppose that the covariance matrices \(\Sigma^\mu\) and
\(\Sigma^\nu\) are estimated from data.
Given \(N\geq Td\) i.i.d. samples
\(\{\mathbf u^i\}_{i=1}^N\) from \(\mu\) and
\(\{\mathbf y^i\}_{i=1}^N\) from \(\nu\), define
\begin{equation*}
\hat\Sigma^\mu
=
\frac1N\sum_{i=1}^N
\mathbf u^i(\mathbf u^i)^\top,
\qquad
\hat\Sigma^\nu
=
\frac1N\sum_{i=1}^N
\mathbf y^i(\mathbf y^i)^\top,
\end{equation*}
 \(\hat\Sigma^\mu\) and \(\hat\Sigma^\nu\) are positive
definite almost surely since $N\geq Td$.
And the corresponding empirical Gaussian marginals are
 $\hat\mu=\mathcal N(0,\hat\Sigma^\mu),
$ $
\hat\nu=\mathcal N(0,\hat\Sigma^\nu).$
In the remainder of this section, we take the marginal
covariance perturbations to be
$\Delta\Sigma^\mu:=\hat\Sigma^\mu-\Sigma^\mu$ and
$\Delta\Sigma^\nu:=\hat\Sigma^\nu-\Sigma^\nu$.
Accordingly, define
\begin{equation*}
\hat\Sigma_0
:=
\begin{bmatrix}
\hat\Sigma^\mu&0\\
0&\hat\Sigma^\nu
\end{bmatrix},
\qquad
\Delta\Sigma_0=\hat\Sigma_0-\Sigma_0.
\end{equation*}
Let $\pi^*$ and $\hat\pi^*$ denote the optimal
couplings of problem \eqref{eq:problem-formulation-2}
corresponding to the true marginals $(\mu,\nu)$
and empirical Gaussian marginals $(\hat\mu,\hat\nu)$,
respectively.
Our goal is to bound
$\|\Sigma^{\hat\pi^*}-\Sigma^{\pi^*}\|_F$.
Combining the deterministic stability result in Theorem~\ref{theorem: stability under marginal covariance perturbations} with
Gaussian covariance concentration bounds
yields the following statistical guarantee.
\begin{theorem}
\label{theorem: error bound}
Let
\(
r_\Sigma
=
r_\Sigma(\Sigma_0,\Sigma^\gamma,T,d)
\)
and
\(
C_\Sigma
=
C_\Sigma(\Sigma_0,\Sigma^\gamma,T,d)
\)
be the quantities defined in
Theorem~\ref{theorem: stability under marginal covariance perturbations}.

There exists a constant \(C>0\) such that,
for every \(\epsilon\geq0\) and every \(N\geq Td\)
satisfying
\begin{equation}
\psi(N,\epsilon,C)
\leq
r_\Sigma,
\label{eqn: sample size condition for statistical stability}
\end{equation}
where for \(N\geq1\), \(\epsilon\geq0\), and \(C>0\), we define
\begin{equation}
\psi(N,\epsilon,C)
:=
\sqrt{Td}\,C
\sqrt{
\|\Sigma^\mu\|_2^2
+
\|\Sigma^\nu\|_2^2
}
\left(
\sqrt{\frac{Td+\epsilon}{N}}
+
\frac{Td+\epsilon}{N}
\right),
\label{eqn: statistical marginal covariance perturbation bound}
\end{equation}
the following bound holds 
\begin{equation}
\|\Sigma^{\hat\pi^*}-\Sigma^{\pi^*}\|_F
\leq
C_\Sigma\psi(N,\epsilon,C),\quad \text{with probability at least \(1-4e^{-\epsilon}\).}
\label{eqn: statistical stability of optimal covariance}
\end{equation}
\end{theorem}
\begin{proof}[Proof sketch]
Applying the Gaussian sample covariance concentration bounds (see \cite[Chapter~2]{Vershynin_2026})
to both marginals gives
\begin{equation*}
\|\hat\Sigma_0-\Sigma_0\|_F
\leq
\psi(N,\epsilon,C)
\end{equation*}
with probability at least $1-4e^{-\epsilon}$.
Under the sample-size condition in the theorem,
Theorem~\ref{theorem: stability under marginal covariance perturbations}
then yields \eqref{eqn: statistical stability of optimal covariance}.
The complete proof is provided in~\Cref{sec: supp of stabiility}.
\end{proof}
For fixed $T$, $d$, and $\epsilon$, Theorem~\ref{theorem: error bound}
establishes an $O(N^{-1/2})$ statistical convergence rate
as $N\to\infty$.

\section{Equivalence with Sinkhorn Iterations}
\label{sec: sinkhorn equivalence}
We establish the equivalence between the BCD algorithm
developed in Section~\ref{section:bcd}
and the Gaussian specialization of the Sinkhorn iterations
studied in \cite{eckstein2023computationalmethodsadaptedoptimal}.
The Sinkhorn algorithm is formulated as an iterative two-step procedure consisting of alternating projections onto the marginal constraint sets. For $k=0,1,2,\ldots$, the coupling is updated via the following two optimization steps:
\begin{align}
	\text{(odd iteration)}\quad &\pi^{ \{k+ \frac{1}{2} \} } :=  \argmin_{\pi \in\Pi_c(\mu, *)} D_{}(\pi, \pi^{ \{k\} }) \label{eqn: odd sinkhorn}, \\
		\text{(even iteration)}\quad &\pi^{ \{ k+1 \} } := \argmin_{\pi \in \Pi(*, \nu)} D_{}(\pi, \pi^{ \{k+ \frac{1}{2}\} }) \label{eqn: even sinkhorn},
\end{align}
where $\Pi_c(\mu, *)$ denotes the set of all causal couplings of $(U,Y)$ with the fixed $U$-marginal equal to $\mu$, and $\Pi(*,\nu)$ denotes the set of all couplings with the fixed $Y$-marginal equal to $\nu$. The odd and even steps can therefore be interpreted as projections onto $\Pi_c(\mu, *)$ and $\Pi(*,\nu)$, respectively. Let \(P=\mu\otimes\nu\). The iteration is initialized by the
coupling \(\pi^{\{0\}}\)
defined by
\begin{equation}
\frac{d\pi^{\{0\}}}{dP}(\mathbf u,\mathbf y)
=
\frac{
\exp\!\left(-c(\mathbf u,\mathbf y)\right)
}{
\displaystyle
\int_{\mathcal U\times\mathcal Y}
\exp\!\left(-c(\mathbf u,\mathbf y)\right)
\,P(d\mathbf u,d\mathbf y)
}.
\label{eqn: init of sinkhorn}
\end{equation}
For comparison with Algorithm~\ref{alg:bcd-dual}, throughout this section
we formulate the Sinkhorn iterations in the Cholesky coordinates introduced
in Section~\ref{section: cholesky Coordinates}. As in
Section~\ref{section: gaussian primal and dual}, we drop the tildes for
notational simplicity. Thus, $\mu=\nu=\mathcal N(0,I_{Td}),$ and hence  $P=\mu\otimes\nu=\mathcal N(0,I_{2Td}),$
and \(c\) denotes the transformed quadratic cost in
\eqref{eqn: quadratic c tilde in cholesky}.  Before establishing the equivalence with the BCD iterations, we first show that the Gaussian family is preserved under both projection steps of the Sinkhorn algorithm.
\begin{lemma}
\label{lemma: gaussianity preserved under sinkhorn projections}
Let
\(\rho\in\mathbb P(\mathcal U\times\mathcal Y)\)
be a non--degenerate Gaussian distribution, and let
\(\mu\) and \(\nu\) satisfy Assumption \ref{assumption:maringal}.
Then the solutions to
\begin{equation*}
\min_{\pi\in\Pi_c(\mu,*)}D(\pi,\rho)
\qquad\text{and}\qquad
\min_{\pi\in\Pi(*,\nu)}D(\pi,\rho)
\end{equation*}
are both non--degenerate Gaussian distributions.
\end{lemma}
\begin{proof}[Proof sketch]
For the causal projection, the backward recursive
construction in
\cite[Lemma~6.2]{eckstein2023computationalmethodsadaptedoptimal}
preserves Gaussianity by induction: starting from the Gaussian terminal conditional distribution, the relative entropy terms remain quadratic in
the conditioning variables, so each resulting conditional distribution remains Gaussian.
For the second projection, the optimizer takes the form
$\pi(d\mathbf u,d\mathbf y)
=\rho(d\mathbf u\mid\mathbf y)\nu(d\mathbf y)$,
which is Gaussian because both $\rho$, $\nu$ are Gaussian.
Moreover, each projection problem has finite optimal
value, since $\mu\otimes\nu$ is feasible and has finite
relative entropy with respect to $\rho$.
Thus, the optimizers are absolutely continuous with
respect to $\rho$ and cannot be degenerate.
The complete proof is provided in~\Cref{sec: supp sinkhorn equivalence}.
\end{proof}
We now establish the equivalence between the
Sinkhorn iterations and Algorithm~\ref{alg:bcd-dual}.
\begin{proposition}
\label{proposition: sinkhorn is bcd}
Under Assumptions~\ref{assumption:cost} and
\ref{assumption:maringal}, let
\(\{\pi^{\{k\}},\pi^{\{k+\frac12\}}\}_{k\geq0}\)
be the Sinkhorn iterates defined by
\eqref{eqn: odd sinkhorn}--\eqref{eqn: even sinkhorn},
initialized at \eqref{eqn: init of sinkhorn}.
Let \(\{\theta^{\{k\}} ,\theta^{\{k+\frac12\}}\}_{k\geq0}\)
be the BCD iterates from Algorithm~\ref{alg:bcd-dual},
initialized at \(\theta^{\{0\}}=\mathbf 0\).
Then 
for every \(k\geq0\),
\begin{equation*}
\left(\Sigma^{\pi^{\{k\}}}\right)^{-1}
=
Z(\theta^{\{k\}}),
\qquad
\left(\Sigma^{\pi^{\{k+\frac12\}}}\right)^{-1}
=
Z(\theta^{\{k+\frac12\}}).
\end{equation*}
\end{proposition}
\begin{proof}[Proof sketch]
We prove the correspondence by induction.
At initialization, \eqref{eqn: init of sinkhorn}
gives a centered Gaussian distribution with
precision matrix
\begin{equation*}
\left(\Sigma^{\pi^{\{0\}}}\right)^{-1}
=
(\Sigma^\gamma)^{-1}+I_{2Td}
=
Z(\theta^{\{0\}}).
\end{equation*}
Suppose the correspondence holds at iteration $k$.
Since the Sinkhorn iterates remain non--degenerate Gaussian
by Lemma~\ref{lemma: gaussianity preserved under sinkhorn projections},
each projection can be expressed as an optimization over
Gaussian covariance matrices.

For the odd projection, rewriting
$D(\pi,\pi^{\{k\}})$ relative to
$P=\mu\otimes\nu$ yields an entropy-regularized
problem with a modified quadratic cost.
After a change of dual variables, its dual formulation
agrees with the first BCD block update.
Similarly, the even projection corresponds to the
second BCD block update.

Consequently, the inverse covariance matrices of
the two Sinkhorn iterates coincide with
$Z(\theta^{\{k+\frac12\}})$ and
$Z(\theta^{\{k+1\}})$, respectively.
The complete proof is provided in~\Cref{sec: supp sinkhorn equivalence}.
\end{proof}

\section{Numerical Experiments}
\label{sec: numerical experiment}


In this section, we evaluate the proposed algorithms on synthetic experiments. The synthetic dataset is generated by a  Gaussian Process (GP) for the input marginal $\mu$, which then propagates through a linear Gaussian State-Space Model (SSM) to define the output marginal $\nu$. We consider a time interval $[0,\text{Total Time}]$ discretized into
$T$ steps, with time step size $\Delta t=\frac{\text{Total Time}}{T}$.
The input marginal \(\mu\) is modeled as a
\(d\)-dimensional Gaussian process with a time-varying
mean and an exponential temporal covariance kernel,
with independent  components.
The detailed specification is provided in~\Cref{section: detail of the input gaussina process}.
The output distribution $\nu$ is generated through a linear Gaussian State-Space Model:
\begin{align}
X_t
&=
F_*X_{t-1}
+
B_*U_t
+
(\Delta tQ_*)^{\frac{1}{2}}V_t, \quad X_0=\mathbf 0
\label{eqn: spring system}\\
Y_t
&=
H_*X_t
+
(\Delta tR_*)^{\frac{1}{2}}W_t,
\nonumber
\end{align}
where
$X_t\in\mathbb R^{d_x}$ and $U_t,Y_t\in\mathbb R^d$,
with
$F_*\in\mathbb R^{d_x\times d_x}$,
$B_*\in\mathbb R^{d_x\times d}$,
$H_*\in\mathbb R^{d\times d_x}$,
$Q_*\in\mathbb R^{d_x\times d_x}$, and
$R_*\in\mathbb R^{d\times d}$, assuming that $Q_*\succeq0$ and $R_*\succ0$.
The random variables
$V_t\sim\mathcal N(0,I_{d_x})$ and
$W_t\sim\mathcal N(0,I_d)$
are independent across time, independent of each other, and independent of $U_{1:T}$. The noise terms are scaled by $\sqrt{\Delta t}$,
consistent with the variance scaling of Brownian increments.
Since the model \eqref{eqn: spring system} is linear Gaussian and
$U_{1:T}$ is Gaussian, the resulting output path $Y_{1:T}$ is also Gaussian.
Unless otherwise specified, we draw $N=10000$
independent trajectories from $\mu$.
Independently, another $N$ input trajectories are sampled from $\mu$ and propagated through \eqref{eqn: spring system}, with independent realizations of $V_{1:T}$ and $W_{1:T}$, to generate samples from the output marginal $\nu$.
For the experiments based on empirical marginal estimates, only the resulting unpaired marginal samples are used by the proposed method, and the input-output pairing generated by the true SSM is not used. 
The marginal means and covariances are estimated from
unpaired samples, and the non-centered marginals are
handled by shifting the coordinate origin using the
estimated means, as in Section~\ref{section: Gaussian reduction}.
Since the true SSM \eqref{eqn: spring system} is linear Gaussian, the true moments of $\nu$ can be obtained directly from the true moments of $\mu$ and the parameters in \eqref{eqn: spring system} above. See \Cref{section: for the simpler true ssm experiment-the ref coupling for stability numeical plot} for its derivation, and how the assumptions here ensure that $\Sigma^\nu\succ0$, as required by Assumption \ref{assumption:maringal}.

We first consider two prior models: 
The first is an SSM prior similar to~\eqref{eqn: spring system} with model parameters different from those of the true
SSM. 
The second is a boundary value problem (BVP) prior with second--order latent dynamics that depend on both past and future states, thereby introducing noncausal dependencies,
motivated by \cite{58529}.
The definitions of both priors and the construction of the corresponding Gaussian couplings appear in~\Cref{sec: From a prior dynamic to a prior coupling}.

In the experiments, we consider the time interval $[0,8]$ discretized into $T=32$ steps. The parameters of the true SSM and the two prior models
are provided in \Cref{section: numercial experiment parameters}. 
Problem~\eqref{eq:problem-formulation-2} is solved using
Algorithm~\ref{alg:bcd-dual} with stopping tolerance
$\varepsilon_{\mathrm{tol}}=10^{-10}$.
We denote the resulting coupling by $\hat{\pi}_c$ and regard it as the optimal causal coupling
up to numerical tolerance.
For comparison, we also solve the entropic optimal transport problem
over the enlarged feasible set $\Pi(\mu,\nu)$, without the causality
constraint. This can be solved by the same BCD scheme with $\Psi$
fixed at zero, alternating only the updates of $\Lambda_u$ and
$\Lambda_y$. We denote the resulting coupling by $\hat{\pi}_{nc}$ and regard it as the optimal non--causal coupling
up to numerical tolerance.

\begin{figure}[!tbp]
   \centering
    \vspace{-0.2cm}
    \begin{adjustbox}{width=0.7\linewidth, center}
\includegraphics{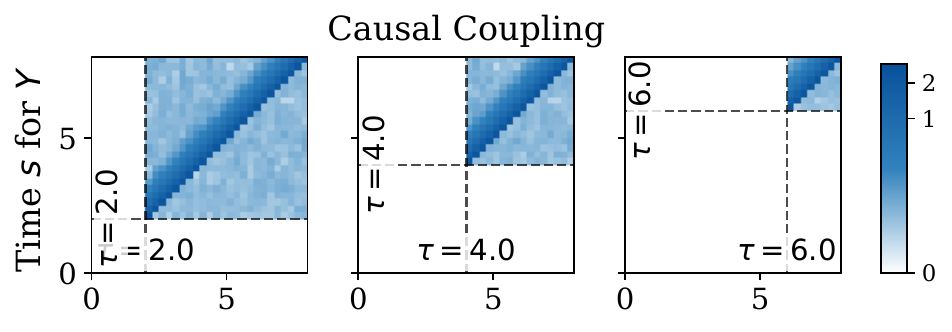}
    \end{adjustbox} \\
   \vspace{-0.2cm}
    \begin{adjustbox}{width=0.7\linewidth, center}
\includegraphics{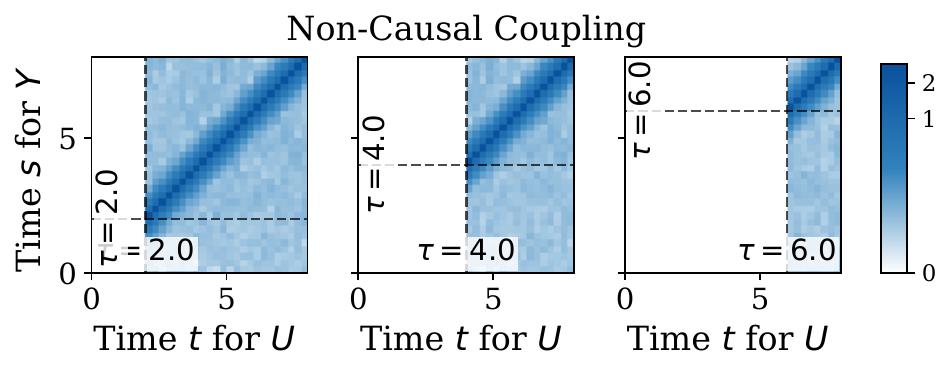}
    \end{adjustbox}
    \vspace{-0.2cm}
   \caption{Conditional cross-covariance $||\text{Cov}(Y_s, U_t \mid U_{1:n_\tau})||_F$ of the optimal couplings associated with the BVP prior (on a symmetrical logarithmic scale), for $\tau=2,4,6$ (columns), with the
conditioning time indicated by dashed lines. Top row: Causal optimal transport coupling $\hat{\pi}_c$. Bottom row: Non-causal optimal transport coupling $\hat{\pi}_{nc}$. The vanishing conditional covariance in the bottom-right quadrant of the causal subplots illustrates the causality empirically.}\label{fig:cross_covariance_plot}
\end{figure}

Figure~\ref{fig:cross_covariance_plot} compares the
conditional cross-covariance
\(\mathrm{Cov}^{\hat\pi}(Y_s,U_t\mid U_{1:n_\tau})\)
where $\hat \pi\in\{\hat{\pi}_c,\hat{\pi}_{nc}\}$, under the BVP prior,
with \(\tau=2,4,6\).
Here $\tau$ denotes physical time and
$n_\tau=\tau/\Delta t$ is the corresponding
discrete-time index.
For \(t>n_\tau\) and \(s\leq n_\tau\), the conditional covariance
vanishes under the causal coupling. Since the couplings
are Gaussian, this implies conditional independence,
confirming the causality constraint.
In contrast, the non-causal coupling has nonzero conditional covariance there.

\begin{figure}[!tbp]
    \centering
    \includegraphics[width=\textwidth]{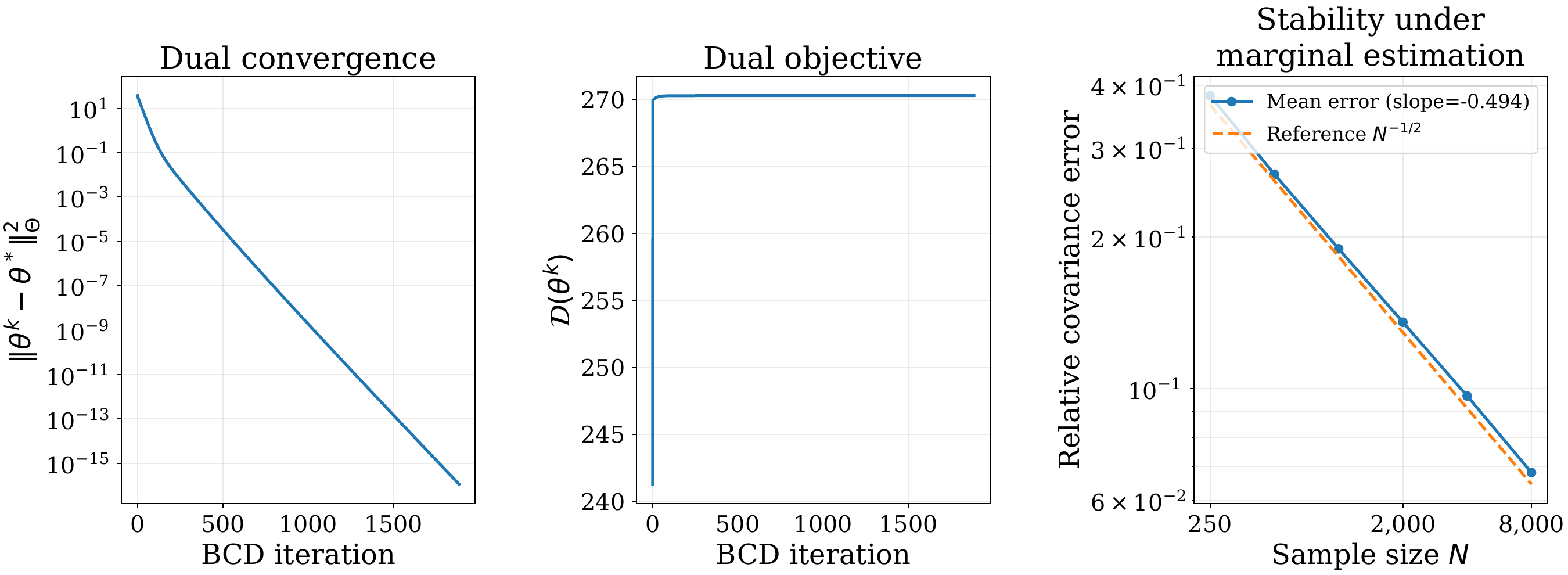}
    \caption{Numerical convergence and stability of the BCD scheme
with $T=32$ and $d=2$.
The left and middle panels use the Gaussian prior with
a known optimizer, while the right panel uses the BVP prior.
    Left: Squared dual-variable error along the BCD iterations,
    exhibiting approximately geometric decay.
    Middle: Dual objective value along the BCD iterations,
    exhibiting monotone convergence.
   Right: Mean relative covariance error under marginal estimation
for different sample sizes $N$, averaged over three
independent repetitions. The dots represent the averaged errors,
and the dashed line shows the convergence
rate $\mathcal{O}(N^{-1/2})$.}
    \label{fig:dual_convergence_and_objective}
\end{figure}

To study numerical convergence, we consider a Gaussian prior
for which the associated ECOT problem has a known exact optimizer.
The construction is detailed in
\Cref{section: a problem with known optimizer}.
To evaluate the numerical convergence of the proposed method, we examine the convergence of the dual variables of the BCD scheme. The dual-variable error
\(e_k=\left\|\theta^{\{k\}}-\theta^*\right\|_\Theta^2,\)
where $\theta^*$ is the exact dual optimizer,
is recorded during the iterations. We also evaluate the dual objective $\mathcal D(\theta^{\{k\}})$ at each BCD iterate. The left panel of Figure \ref{fig:dual_convergence_and_objective} displays the dual-variable error against iteration $k$ on a semi-log scale and exhibits approximately geometric decay.  The middle panel of Figure \ref{fig:dual_convergence_and_objective} shows that the dual objective is monotone along the BCD iterations. These numerical results are consistent with the linear convergence established in Section \ref{sec: linear convergence}.

We also evaluate the empirical convergence of the optimal coupling covariance with respect to the Monte Carlo sample size $N$. We study the relative error 
\begin{equation}
\label{eqn: relatvie error def for the stablity plot}
\frac{
\|\Sigma^{\hat{\pi}_{c,N}}-
\Sigma^{\pi_{c,\mathrm{ref}}}\|_F
}{
\|\Sigma^{\pi_{c,\mathrm{ref}}}\|_F
}.
\end{equation}
where $\hat{\pi}_{c,N}$ denotes the causal coupling
computed using marginals estimated from $N$
independent simulated trajectories.
The reference coupling $\pi_{c,\mathrm{ref}}$ is
computed using the true marginal moments
$\Sigma^\mu$ and $\Sigma^\nu$
(see~\Cref{section: for the simpler true ssm experiment-the ref coupling for stability numeical plot} for their derivation),
rather than moments estimated from finite
Monte Carlo samples.
The same prior covariance \(\Sigma^\gamma\), constructed from the BVP prior and used in the corresponding empirical experiment,
is kept fixed.
For each sample size $N$, we repeat the experiment three times
with independently generated empirical marginals. In each repetition,
we compute the relative error defined in \eqref{eqn: relatvie error def for the stablity plot} and then average these three errors.
The right panel of Figure~\ref{fig:dual_convergence_and_objective} reports the resulting log–log plot of the averaged error, together with a line with the convergence rate $\mathcal{O}(N^{-1/2})$. A linear regression in log–log scale yields a slope close to $-0.5$, consistent with the theoretical convergence rate established in Theorem \ref{theorem: error bound}. 

To evaluate the predictive performance of the learned transport plans, we consider an input-only prediction task, where the model observes \(U_{1:t}\) to predict \(Y_t\) at each time $t$. We simulate a test dataset comprising independent sample paths of coupled inputs and outputs generated via the true SSM \eqref{eqn: spring system}. Besides  $\hat{\pi}_{nc}$, we consider a prefix non-causal coupling as another benchmark. Instead of solving one non--causal optimal transport problem over the full time horizon, we independently solve a non--causal problem over $(U_{1:t},Y_{1:t})$ for each $t=1,\ldots,T$. The reference for each problem is obtained by marginalizing the original full-horizon prior $\gamma$ onto the corresponding prefix, and the marginal constraints are the corresponding prefixes of $\mu$ and $\nu$. The coupling fitted at time $t$ is used only to predict $Y_t$. The detailed formulation in Cholesky coordinates is given in \Cref{section: definition of the prefix non--causal}.
We evaluate the causal, non-causal, and prefix non-causal couplings,
together with the prior $\gamma$. All three fitted methods use the same
estimated marginal distributions and the same prior reference. The prior
itself does not incorporate the observed marginal samples.
For each coupling, the corresponding linear Gaussian
predictor is obtained from its conditional distribution.
The explicit formulas are given in \Cref{eqn: section: prediction task setup}.

The prior $\gamma$ is also evaluated using the corresponding
conditional prediction formula. Since it does not incorporate the
empirical marginal distributions and has substantially worse
prediction performance, we omit its curve from the figures to make
the differences among the three fitted couplings more visible.

\begin{figure}[!tbp]
    \centering
    \includegraphics[width=0.9\textwidth]{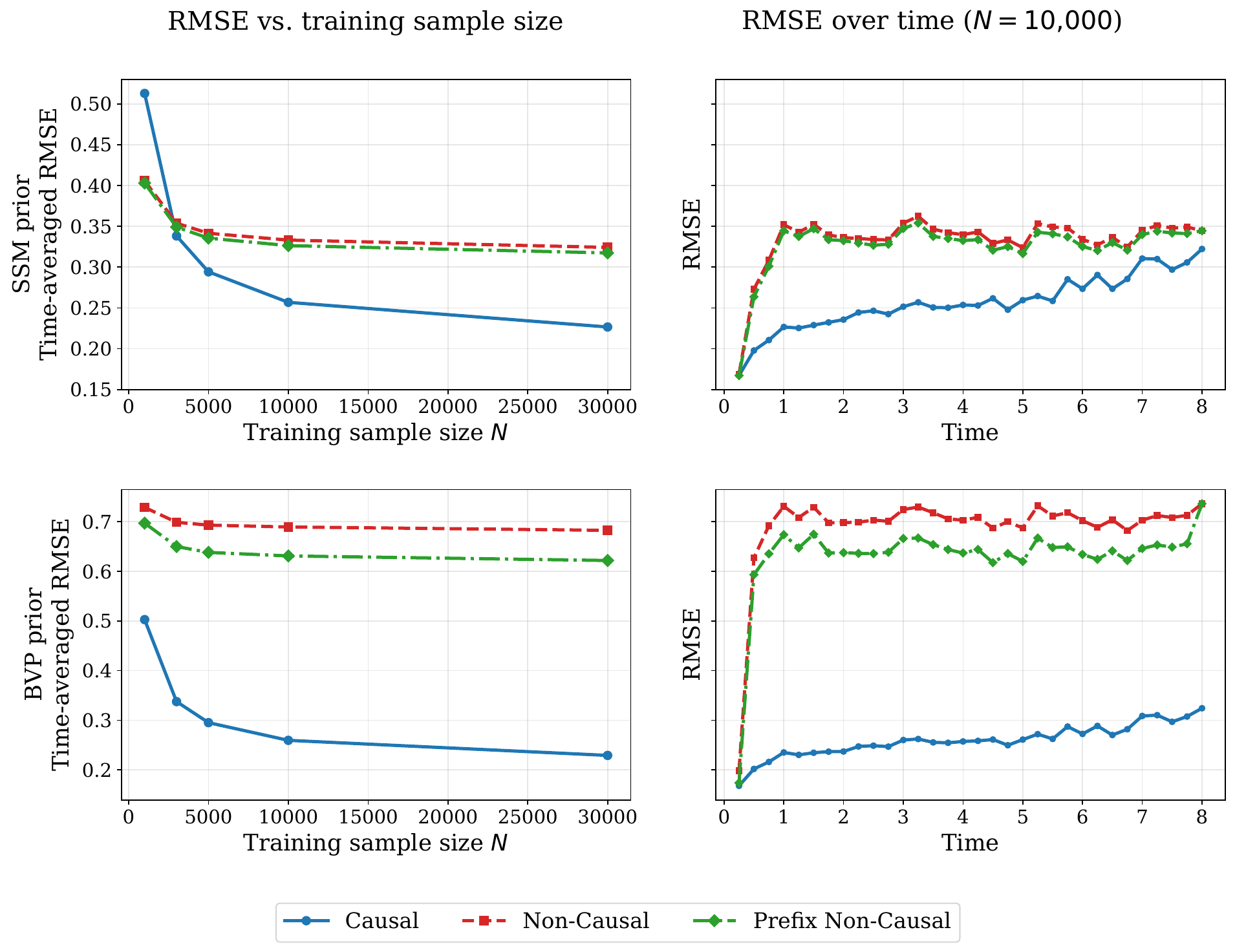}
    \caption{Prediction RMSE for the causal, non-causal, and prefix
    non-causal couplings under input-only prediction.
    The top row uses the SSM prior and the bottom row uses the BVP prior.
    Left: Time-averaged RMSE against the marginal estimation sample size $N$.
    Right: RMSE over time with $N=10000$.
    The total time horizon is 8, with $T=32$ time steps and
    data dimension $d=2$.}
    \label{fig:simulation_input_only}
\end{figure}

To capture prediction errors across the temporal horizon, we analyze
the Root Mean Squared Error (RMSE) at each time $t$.
For a given model, the residual for sample path $i$, at time step $t$,
and dimension component $j\in\{1,\dots,d\}$ is defined as
$e_{t,j}^{(i)}=y_{t,j}^{(i)}-\hat y_{t,j}^{(i)}$,
where $\hat y_{t,j}^{(i)}$ is the prediction generated using
the input history $U_{1:t}$.
The RMSE is computed by averaging the squared residuals across the
$N^{\text{test}}=10000$ test paths and across all $d$ 
dimensions:
\begin{equation*}
\text{RMSE}(t)
=
\sqrt{
\frac{1}{d}\sum_{j=1}^d
\left(
\frac{1}{N^{\text{test}}}
\sum_{i=1}^{N^{\text{test}}}
\left|y_{t,j}^{(i)}-\hat y_{t,j}^{(i)}\right|^2
\right)
}.
\end{equation*}
We further examine the effect of the marginal estimation sample size
$N$ on the time-averaged RMSE, defined as
\(\overline{\text{RMSE}}
=\frac{1}{T}\sum_{t=1}^T\text{RMSE}(t).\)
Figure~\ref{fig:simulation_input_only} presents the time-averaged
RMSE as a function of $N$ in the left column and the RMSE over time
for $N=10000$ in the right column.
The full numerical values, including those for the prior \(\gamma\),
are reported in~\Cref{extra section: numerical results}.

For the SSM prior, which is structurally same as the true model, the difference among the three fitted couplings depends
more strongly on the sample size used to estimate the marginals.
For smaller sample sizes, the two non-causal methods can outperform
the causal coupling under input-only prediction.
However, as the sample size increases, the causal coupling becomes
increasingly competitive and eventually outperforms both non-causal
alternatives, as illustrated in the left column of
Figure~\ref{fig:simulation_input_only}.
This trend suggests that the competitive performance of the non-causal
methods at smaller sample sizes may be related to differences in
their sensitivity to marginal estimation errors, while the causal
coupling provides better prediction when the marginal estimates
are sufficiently accurate.

For the BVP prior, which is structurally different from the true model, the causal coupling maintains considerably lower
RMSE than both non-causal couplings,
as shown in the bottom row of Figure~\ref{fig:simulation_input_only}.
In particular, the full-horizon non-causal coupling exhibits
substantially higher prediction errors, while the prefix non-causal
coupling reduces these errors but remains less accurate than the
causal coupling. This behavior persists across the sample sizes
considered in the left column of Figure~\ref{fig:simulation_input_only}.
This illustrates that the causal coupling can mitigate the effects
of prior model misspecification and achieve better prediction performance.

We also compare the computational performance of the proposed BCD algorithm
with CVXPY for different time horizons $T$, with $d=2$. 
We use CVXPY with its default solver settings to solve the primal
problem \eqref{eqn:problem-over-cov} and the dual problem
\eqref{eqn: dual problem}, referred to as CVXPY primal and
CVXPY dual, respectively. 
Both formulations are solved using SCS (\cite{scs}), as selected by CVXPY.
For BCD, we consider stopping tolerances $\varepsilon_{\mathrm{tol}}\in\{10^{-12},10^{-13},10^{-14}\}$.
 For each $T$, we use the same construction of a Gaussian prior
with a known exact optimizer as in the convergence experiment above.
Both BCD and CVXPY solve the corresponding covariance problem
in Cholesky coordinates, using the true marginal covariances.
The experiments were conducted on a server equipped with
Intel Xeon E5-2698 v4 processors and 503 GiB of memory.

\begin{figure}[!tbp]
\centering\includegraphics[width=\textwidth]{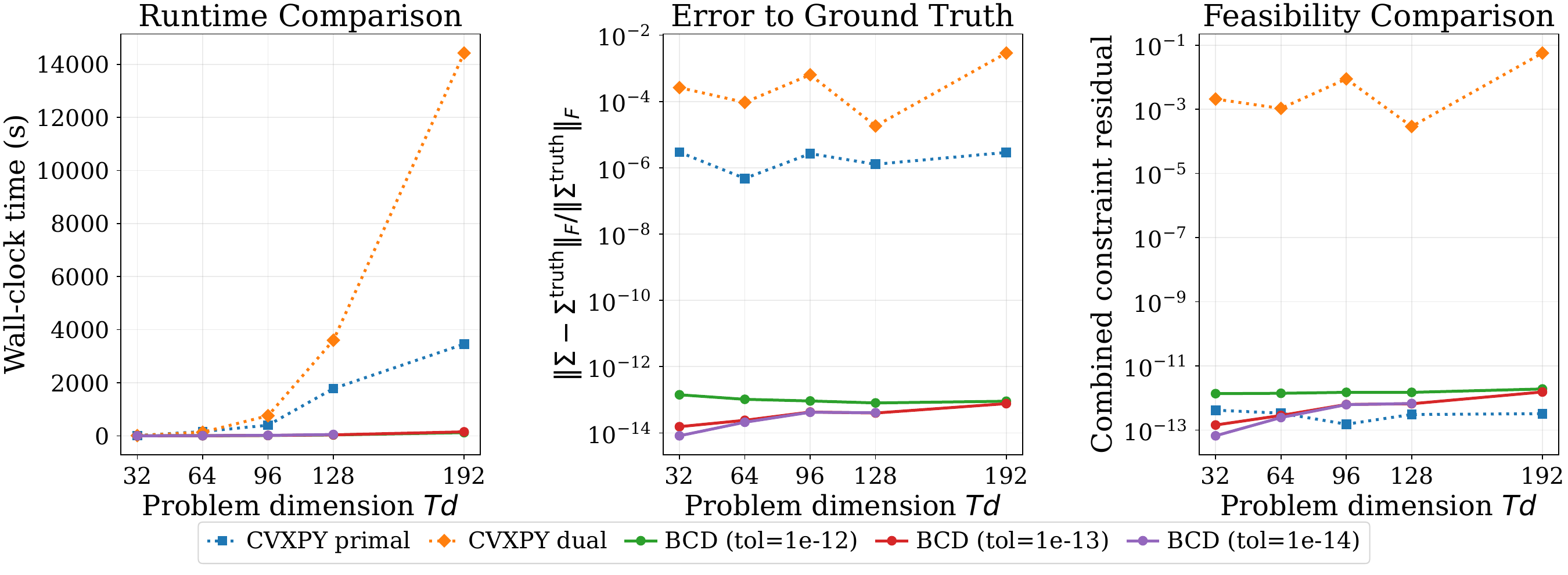}
    \caption{Comparison of BCD and CVXPY for different problem
    dimensions $Td$, with $d=2$.
    Left: Wall-clock runtimes, reported as medians over three runs.
  Middle: Relative Frobenius error of the covariance
with respect to the exact optimizer $\Sigma^{\text{truth}}$.
Right: Combined Frobenius norm of the marginal and causality
constraint violations.
}
    \label{fig:benchmark_comparison}
\end{figure}

The left panel of Figure~\ref{fig:benchmark_comparison}
compares the wall-clock runtimes as a function of the problem
dimension $Td$. As the problem dimension increases, BCD becomes
substantially faster than both the primal and dual CVXPY
formulations, even with a tighter stopping tolerance of $10^{-14}$.
Moreover, the runtimes of BCD remain similar across the three
stopping tolerances.
Unlike the proposed BCD algorithm, which exploits the structure
of problem \eqref{eqn:problem-over-cov}, CVXPY, when using SCS,
reformulates the log-determinant objectives using auxiliary
variables and convex cone constraints.
This reformulation may contribute to the higher computational
cost observed for the CVXPY formulations.
The middle panel of Figure~\ref{fig:benchmark_comparison}
compares the relative covariance errors
of the five solver configurations with respect to the exact optimizer.
The right panel of Figure~\ref{fig:benchmark_comparison}
reports the Frobenius norm of the constraint violations, computed
as the square root of the sum of the squared Frobenius norms of the
two marginal covariance constraint violations and the individual
causality constraint violations.
The constraint residuals of BCD generally decrease as the
stopping tolerance is tightened, reaching approximately
$10^{-13}$ or lower for smaller problem dimensions.
For larger dimensions, the residuals increase slightly,
while remaining below $10^{-12}$ for the tightest tolerance.
The CVXPY primal formulation achieves residuals of a similar
order of magnitude, whereas the CVXPY dual formulation
exhibits substantially larger residuals that fluctuate
with the problem dimension $Td$.
The reported runtimes are medians over three repeated
runs, while  the relative covariance errors  and constraint residuals
are evaluated from the first run.
These results demonstrate the computational efficiency and accuracy
of BCD.

\section{Conclusion}
In this work, we studied how to infer an input--output
system from unpaired marginal observations and a given prior model while enforcing a one-sided causality
constraint through entropic causal optimal transport.
In the Gaussian setting with general quadratic costs, we showed that the optimal coupling is Gaussian and developed an efficient block coordinate descent algorithm for computing it.
We established linear convergence of the algorithm and derived
stability and statistical error bounds under marginal estimation.
Numerical experiments demonstrated the computational efficiency
and accuracy of the proposed algorithm and investigated the
predictive performance of the inferred causal coupling.
Future research may explore applications of the proposed framework to time-series forecasting, sequential generative modeling \cite{XuEtAl2020COTGAN,AcciaioEcksteinHou2024TCVAE}, and the adaptation of non-causal dynamical models to causal settings.


\bibliographystyle{siamplain}
\bibliography{references}
\appendix

Detailed proofs and additional technical details are provided in the following appendices.
\section{Differential Notation and Conventions}
{
We first recall the differential notation and Hessian convention that will be used in the calculations below. \(D\) denotes the Fr\'echet derivative. Let
\(\mathcal X\) be a finite-dimensional real inner-product space with inner
product \(\langle\cdot,\cdot\rangle_{\mathcal X}\), and let
\(\Omega\subset\mathcal X\) be open. For a differentiable function
\[
\varphi:\Omega\to\mathbb R,
\]
the Fr\'echet derivative at \(x\in\Omega\) is denoted by
\[
D\varphi(x)[\cdot]:\mathcal X\to\mathbb R.
\]
The gradient \(\nabla\varphi(x)\in\mathcal X\) is the unique element
satisfying
\[
D\varphi(x)[\delta x]
=
\left\langle
\nabla\varphi(x),
\delta x
\right\rangle_{\mathcal X},
\qquad
\delta x\in\mathcal X.
\]

If \(\varphi\) is twice Fr\'echet differentiable, its second Fr\'echet
derivative is the bilinear form defined by
\[
D^2\varphi(x)
[\delta x_1,\delta x_2]
:=
D(D\varphi)(x)[\delta x_2][\delta x_1],
\qquad
\delta x_1,\delta x_2\in\mathcal X.
\]
We use the Hessian in its bilinear form and write
\[
\nabla^2\varphi(x)
[\delta x_1,\delta x_2]
:=
D^2\varphi(x)
[\delta x_1,\delta x_2].
\]
Let  \(\mathcal Y\) be a finite--dimensional real normed space. 
Let \(I\subset\mathbb R\) be open. Let
\(r:I\to\mathcal X\), \(\theta:I\to\Omega\), and \(q:\Omega\to\mathcal Y\) be Fr\'echet differentiable.  We will use the following standard identities whenever the corresponding
derivatives are defined.

\begin{equation}
\frac{d}{dt}r(t)
=
Dr(t)[1],
\qquad t\in I.
\label{eqn: ordinary derivative equals Frechet derivative at 1}
\end{equation}

\begin{equation}
\frac{d}{dt}q(\theta(t))
=
Dq(\theta(t))
\left[
\frac{d}{dt}\theta(t)
\right],
\qquad t\in I.
\label{eqn: time derivative of composition is D  in directionn of inner function time derivative}
\end{equation}

\begin{equation}
\nabla^2\varphi(x)
[\delta x_1,\delta x_2]
=
\left\langle
\delta x_1,
D\nabla\varphi(x)[\delta x_2]
\right\rangle_{\mathcal X}, \quad \forall x\in\Omega, \delta x_1,\delta x_2\in\mathcal X
\label{eqn: Hessian bilinear form via gradient derivative}
\end{equation}

\section{Supplementary Analysis for Section~\ref{sec: optimization in cov space}}
\label{sec: supplementary gaussian ecot}
\subsection{Supplementary Analysis for Section~\ref{section: Gaussian reduction}} \label{sec: supplementary Gaussian reduction}
Following Remark 2.2 in \cite{eckstein2023computationalmethodsadaptedoptimal}, the causality constraint \eqref{eq:causal-constraint} can be simplified into a sequential single-step constraint,
under $\pi\in\mathbb{P}(\mathcal{U}\times\mathcal{Y})$,
\[\forall t\in \{1,\cdots,T-1\},\quad (Y_1, \dots, Y_t)  \indep_{(U_1,\dots,U_t)}(U_{t+1}, \dots, U_T )\]
if and only if  
\begin{equation}
    \forall t\in \{1,\cdots,T-1\}, \quad  (Y_1, \dots, Y_t)\indep_{(U_1,\dots,U_t)} U_{t+1}. \label{eqn: light causal constraint}
\end{equation}

{The following lemma shows that, for a quadratic cost, the non--centered ECOT problem \eqref{eq:problem-formulation-2} can be reduced to the centered case.
\begin{lemma}\label{lemma: centering gaussian marginals}
    Consider the case when the marginals $\mathbf U \sim \mu$ and $\mathbf Y \sim \nu$ have non-zero means, $m^\mu$ and $m^\nu$, respectively. Define the centered variables 
\[ \mathbf U^0= \mathbf U-m^\mu,
\qquad
 \mathbf Y^0= \mathbf Y-m^\nu,
\]
admitting the centered marginals $\mu^0$ and $\nu^0$, respectively. Also, for any coupling 
$\pi$ on $(\mathbf U,\mathbf Y)$, define its centered version $\pi^0$ as the coupling between  $(\mathbf U^0,\mathbf Y^0) = (
    \mathbf U - m^\mu, \mathbf Y - m^\nu )$.  
Then, under the quadratic cost Assumption~\ref{assumption:cost}, 
\begin{equation*}
    \pi^*\in
\argmin_{\pi\in\Pi_c(\mu,\nu)}
\{
\int c\,d\pi+D(\pi,\mu\otimes\nu)
\}  \Leftrightarrow (\pi^*)^0\in
\argmin_{\pi^0\in\Pi_c(\mu^0,\nu^0)}
\{
\int c\,d\pi^0+D(\pi^0,\mu^0\otimes\nu^0)
\}.
\end{equation*}
\end{lemma}

}
\begin{proof}[Proof of Lemma~\ref{lemma: centering gaussian marginals}]
First, we show that 
\begin{equation}\label{eq:pi0-constraint-equivalency}
\pi\in\Pi_c(\mu,\nu)
\quad\Longleftrightarrow\quad
\pi^0\in\Pi_c(\mu^0,\nu^0).
\end{equation}
$\pi^0 \in \Pi(\mu^0,\nu^0)$  by definition. 
The causality constraint is preserved because conditional independence is not changed by adding or subtracting deterministic constants, i.e., for every
\(t\in\{1,\ldots,T-1\}\), we have
\[
U_{t+1}\indep_{U_{1:t}}Y_{1:t}
\quad\Longleftrightarrow\quad
U^0_{t+1}\indep_{U^0_{1:t}}Y^0_{1:t},
\]
Next, we show that the optimization objective function changes by a constant when the coupling becomes centered. The relative entropy term $D(\pi,\mu \otimes \nu)=D(\pi^0,\mu^0 \otimes \nu^0)$ because both distribution arguments change by the same translation. The cost term 
\begin{align*}
    \int c \,\ud \pi &= \Expect_{(\mathbf U,\mathbf Y)\sim \pi}\left[ \frac12
\begin{bmatrix}
\bfU\\
\bfY
\end{bmatrix}^{\top}
(\Sigma^\gamma)^{-1}
\begin{bmatrix}
\bfU\\
\bfY
\end{bmatrix}\right]
\\&=
\Expect_{(\mathbf U^0,\mathbf Y^0)\sim \pi^0}\left[ \frac12
\begin{bmatrix}
\bfU^0 + m^\mu\\
\bfY^0 + m^\nu
\end{bmatrix}^{\top}
(\Sigma^\gamma)^{-1}
\begin{bmatrix}
\bfU^0+ m^\mu\\
\bfY^0+ m^\nu
\end{bmatrix}\right]\\
&=
\Expect_{(\mathbf U^0,\mathbf Y^0)\sim \pi^0}\left[ \frac12
\begin{bmatrix}
\bfU^0\\
\bfY^0 
\end{bmatrix}^{\top}
(\Sigma^\gamma)^{-1}
\begin{bmatrix}
\bfU^0\\
\bfY^0
\end{bmatrix}\right] +
\frac12
\begin{bmatrix}
m^\mu\\\
m^\nu
\end{bmatrix}^{\top}
(\Sigma^\gamma)^{-1}
\begin{bmatrix}
m^\mu\\\
m^\nu
\end{bmatrix}\\
&= \int  c \, \ud \pi^0 + \text{const.}
\end{align*} 
where the linear term vanishes because $\bfU^0$ and $\bfY^0$ are centered. This observation, along with~\eqref{eq:pi0-constraint-equivalency} concludes the proof.
\end{proof}
Since Gaussian couplings will play a central role in this work, we first characterize the causality constraint in \eqref{eqn: light causal constraint} for a Gaussian coupling in terms of its covariance matrix. In order to do so, we introduce the following notations. We use 
$\Sigma^\pi_{A,B}$ to denote the cross-covariance matrix between random vectors $A$ and $B$ under $\pi$, i.e. $\Sigma^\pi_{A,B}= \Cov_\pi(A,B):= \Expect_\pi[(A-\Expect_\pi[A])(B-\Expect_\pi[B])^\top]$ . We omit the superscript or subscript \(\pi\) when the underlying coupling is clear from the context. Moreover, we use the  notation
$\Sigma^\mu_{U_{1:t},U_{1:t}}=\text{Cov}(U_{1:t},U_{1:t})$
and
$\Sigma^\mu_{U_{1:t},U_{t+1}}=\text{Cov}(U_{1:t},U_{t+1})$ where $U_{1:T} \sim \mu$.
\begin{lemma}
\label{lemma: gaussian causal constraint}
Let $\pi$ be a non--degenerate Gaussian coupling of $\begin{bmatrix}
    U_{1:T}\\Y_{1:T}
\end{bmatrix}$ with covariance matrix
$\Sigma$. $\pi$ is causal if and only if
\begin{equation}
0 =\Sigma_{Y_{1:t},U_{t+1}}
-
\Sigma_{Y_{1:t},U_{1:t}}
\left(\Sigma^\mu_{U_{1:t},U_{1:t}}\right)^{-1}
\Sigma^\mu_{U_{1:t},U_{t+1}},
\qquad
\forall t\in\{1,\ldots,T-1\}.
\label{eqn: causal constraint in the gaussian joint cov}
\end{equation}
\end{lemma}
\begin{proof}[Proof of Lemma~\ref{lemma: gaussian causal constraint}]
Under a Gaussian coupling, the conditional distribution of
$\begin{bmatrix}
Y_{1:t}\\
U_{t+1}
\end{bmatrix}$ given $U_{1:t}$ is Gaussian. For a jointly Gaussian vector, independence is equivalent to
zero covariance. In particular,
\eqref{eqn: light causal constraint} is equivalent to
\[
\Cov(Y_{1:t},U_{t+1}\mid U_{1:t})=0,
\qquad
\forall t\in\{1,\ldots,T-1\}.
\]
By properties of conditional Gaussian distributions,
\begin{align*}
&\Cov\left(
\begin{bmatrix}
Y_{1:t}\\
U_{t+1}
\end{bmatrix},\begin{bmatrix}
Y_{1:t}\\
U_{t+1}
\end{bmatrix}
\Bigg|\,U_{1:t}
\right)\\
&=
\begin{bmatrix}
\Cov(Y_{1:t},Y_{1:t}\mid U_{1:t})
&
\Cov(Y_{1:t},U_{t+1}\mid U_{1:t})
\\
\Cov(U_{t+1},Y_{1:t}\mid U_{1:t})
&
\Cov(U_{t+1},U_{t+1}\mid U_{1:t})
\end{bmatrix}\\
&=
\begin{bmatrix}
\Sigma_{Y_{1:t},Y_{1:t}}
&
\Sigma_{Y_{1:t},U_{t+1}}
\\
\Sigma_{U_{t+1},Y_{1:t}}
&
\Sigma^\mu_{U_{t+1},U_{t+1}}
\end{bmatrix}-
\begin{bmatrix}
\Sigma_{Y_{1:t},U_{1:t}}
\\
\Sigma^\mu_{U_{t+1},U_{1:t}}
\end{bmatrix}
\left(\Sigma^\mu_{U_{1:t},U_{1:t}}\right)^{-1}
\begin{bmatrix}
\Sigma_{U_{1:t},Y_{1:t}}
&
\Sigma^\mu_{U_{1:t},U_{t+1}}
\end{bmatrix}.
\end{align*}
Reading off the upper-right block gives
\[
\Cov(Y_{1:t},U_{t+1}\mid U_{1:t})
=
\Sigma_{Y_{1:t},U_{t+1}}
-
\Sigma_{Y_{1:t},U_{1:t}}
\left(\Sigma^\mu_{U_{1:t},U_{1:t}}\right)^{-1}
\Sigma^\mu_{U_{1:t},U_{t+1}},
\]
which proves the result.
\end{proof}

\begin{proof}[Proof of Lemma~\ref{theorem: optomality of gaussian}] \leavevmode\\
\begin{enumerate}[label=\alph*.]
    \item 
    First, to see why the value in
\eqref{eq:problem-formulation-2} is bounded above, consider the product coupling $\mu\otimes\nu$, which is feasible
\(
\mu\otimes\nu\in\Pi_c(\mu,\nu).
\)
The entropy part in \eqref{eq:problem-formulation-2} is simply zero,
\[
D(\mu\otimes\nu,\mu\otimes\nu)=0.
\]
By \eqref{eqn: quad c with pd hessian},
\begin{align*}
\int c(\mathbf u,\mathbf y)\,\ud(\mu\otimes\nu)
&=
\frac12
\int
\begin{bmatrix}
\mathbf u\\
\mathbf y
\end{bmatrix}^{\top}
(\Sigma^\gamma)^{-1}
\begin{bmatrix}
\mathbf u\\
\mathbf y
\end{bmatrix}
\,\ud(\mu\otimes\nu)\\
&\leq
\frac12\|(\Sigma^\gamma)^{-1}\|_2
\left(
\int\|\mathbf u\|^2\,\ud\mu
+
\int\|\mathbf y\|^2\,\ud\nu
\right)\\
&<\infty,
\end{align*}
where the last inequality follows from the finiteness of
\(\Sigma^\mu\), \(\Sigma^\nu\), and \(\Sigma^\gamma\).
Next, the objective in
\eqref{eq:problem-formulation-2} is bounded from below, because \(c\) is bounded from below and
\(
D(\pi,\mu\otimes\nu)\geq0
\)
for every feasible \(\pi\). Consequently,
the value of problem \eqref{eq:problem-formulation-2} is finite.

    \item Since the cost $c$ is continuous and bounded from below, existence and uniqueness of the optimizer of \eqref{eq:problem-formulation-2} follow from \cite[Remark 2.3]{eckstein2023computationalmethodsadaptedoptimal}.
    \item Let \(\pi\in\Pi_c(\mu,\nu)\), with covariance $\Sigma^\pi$, have finite objective value in
\eqref{eq:problem-formulation-2}, and let \(\pi_G\) be the Gaussian distribution having the same mean and covariance as \(\pi\).
\begin{itemize}
    \item 
{$\Sigma^\pi$ is finite.} By definition
\begin{align*}
    \Sigma^\pi = \begin{bmatrix}
        \Expect^\pi[\mathbf U \mathbf U^\top] & \Expect^\pi[\mathbf U \mathbf Y^\top] \\
        \Expect^\pi[\mathbf Y \mathbf U^\top] & 
        \Expect^\pi[\mathbf Y \mathbf Y^\top]
    \end{bmatrix} = \begin{bmatrix}
        \Sigma^\mu & \Expect^\pi[\mathbf U \mathbf Y^\top] \\
        \Expect^\pi[\mathbf Y \mathbf U^\top] & 
        \Sigma^\nu
    \end{bmatrix}
\end{align*}
The diagonal blocks $\Sigma^\mu$ and $\Sigma^\nu$ are finite according to Assumption~\ref{assumption:maringal}. The off-diagonal blocks are finite by Jensen's inequality and the Cauchy--Schwarz inequality:  $\forall\, i,j=1,\ldots,T,$
\begin{align*}
\left\|\Expect^\pi[U_iY_j^\top]\right\|_F
&\leq
\Expect^\pi\!\left[\|U_iY_j^\top\|_F\right]\\
&=
\Expect^\pi\!\left[\|U_i\|\,\|Y_j\|\right]\\
&\leq
\bigl(\Expect^\pi[\|U_i\|^2]\bigr)^{1/2}
\bigl(\Expect^\pi[\|Y_j\|^2]\bigr)^{1/2}\\
&=
\bigl(\Expect^\mu[\|U_i\|^2]\bigr)^{1/2}
\bigl(\Expect^\nu[\|Y_j\|^2]\bigr)^{1/2}
<\infty.
\end{align*}
\item {$\pi_G$ is non-degenerate Gaussian.}  Since $\pi$ has finite objective value and $c\geq0$, $D(\pi,\mu\otimes\nu)<\infty.$
Hence,
\[
\pi\ll\mu\otimes\nu.
\]
By Assumption~\ref{assumption:maringal}, $\mu\otimes\nu$ is a
non--degenerate Gaussian distribution and therefore has a strictly positive
density with respect to the Lebesgue measure. Thus,
\[
\pi\ll\mathrm{Leb}.
\]
If there exists a nonzero vector $a\in \Re^{2Td}$ such that
\[
0=a^\top\Sigma^\pi a
=\text{Var}_{X\sim \pi}(a^\top X),
\]
then $a^\top X$ is constant $\pi$-a.s.. Thus, $\pi$ would assign probability one to a hyperplane, contradicting absolute continuity of $\pi$. Hence $\Sigma^\pi\succ0$, concluding  $\pi_G$ is non-degenerate Gaussian.
\item $\pi_G\in\Pi(\mu,\nu)$. Since $\pi_G$ has the same mean vector and covariance matrix as $\pi$, its
$U$ and $Y$ marginals have means and covariance matrices
$(0,\Sigma^\mu)$ and $(0,\Sigma^\nu)$, respectively. Since these
marginals are Gaussian, they coincide with $\mu$ and $\nu$. Hence,
$\pi_G\in\Pi(\mu,\nu)$.

\item $\pi_G$ is causal. Since $\pi_G$ is Gaussian, by Lemma \ref{lemma: gaussian causal constraint}, it is sufficient to show 
\begin{align*}
    &\Cov_{\pi_G}(Y_{1:t},U_{t+1}\mid U_{1:t}) \\ &=\Sigma^{\pi}_{Y_{1:t},U_{t+1}} - \Sigma^{\pi}_{Y_{1:t},U_{1:t}} \bigl(\Sigma^\mu_{U_{1:t},U_{1:t}}\bigr)^{-1} \Sigma^\mu_{U_{1:t},U_{t+1}}, 
\end{align*}
since \(\pi_G\) and \(\pi\) have the same covariance matrix.
We simplify the first term to show that the conditional covariance is zero. 
\begin{align*}
\Sigma^\pi_{Y_{1:t},U_{t+1}}
&=
\Expect_\pi[Y_{1:t}U_{t+1}^\top]\\
&=
\Expect_\pi\left[
Y_{1:t}
\Expect_\pi[U_{t+1}\mid U_{1:t},Y_{1:t}]^\top
\right].
\end{align*}
Since $\pi$ is causal and $\mu$ is Gaussian,
\begin{align*}
\Expect_\pi[U_{t+1}\mid U_{1:t},Y_{1:t}]
&=
\Expect_\pi[U_{t+1}\mid U_{1:t}]\\
&=
\Expect_\mu[U_{t+1}\mid U_{1:t}]\\
&=
\Sigma^\mu_{U_{t+1},U_{1:t}}
\left(\Sigma^\mu_{U_{1:t},U_{1:t}}\right)^{-1}
U_{1:t}.
\end{align*}
Therefore,
\begin{align*}
\Sigma^\pi_{Y_{1:t},U_{t+1}}
&=
\Expect_\pi\left[
Y_{1:t}U_{1:t}^\top
\right]
\left(\Sigma^\mu_{U_{1:t},U_{1:t}}\right)^{-1}
\Sigma^\mu_{U_{1:t},U_{t+1}}\\
&=
\Sigma^\pi_{Y_{1:t},U_{1:t}}
\left(\Sigma^\mu_{U_{1:t},U_{1:t}}\right)^{-1}
\Sigma^\mu_{U_{1:t},U_{t+1}}.
\end{align*}
Hence,
\[
\Cov_{\pi_G}(Y_{1:t},U_{t+1}\mid U_{1:t})=0,
\]
proving that $\pi_G$ is causal.
\item Finally, we show that
\[
\int c\,d\pi_G+D(\pi_G,\mu\otimes\nu)
\leq
\int c\,d\pi+D(\pi,\mu\otimes\nu),
\]
proving that at optimality, we must have $\pi=\pi_G$. 
In order to prove the claim, we use the identities 
\begin{align}
   & \int c\,d\pi=\int c\,d\pi_G, \label{eqn: c integral is same for pi and piG}\\
    &\int \log\frac{d\pi_G}{d(\mu\otimes\nu)}\,d\pi
=
\int \log\frac{d\pi_G}{d(\mu\otimes\nu)}\,d\pi_G, \label{eqn: log density integral is same for pi and piG}
\end{align}
which are true because $c$ and $\log(\frac{d\pi_G}{d(\mu\otimes\nu)})$ are quadratic functions and $\pi$ and $\pi_G$ share the same first two moments.
Since $\pi_G$ is non--degenerate Gaussian, it has a strictly positive density
with respect to the Lebesgue measure. Together with $\pi\ll\mathrm{Leb}$, this gives
\[
\pi\ll\pi_G.
\]
We use \eqref{eqn: log density integral is same for pi and piG} to decompose the KL-divergence according to
\begin{align*}
D_{}(\pi,\mu\otimes\nu)
&=
\int \log\frac{d\pi}{d(\mu\otimes\nu)}\,d\pi\\
&=
\int \log\frac{d\pi}{d\pi_G}\,d\pi
+
\int \log\frac{d\pi_G}{d(\mu\otimes\nu)}\,d\pi\\
&=
D_{}(\pi,\pi_G)+D_{}(\pi_G,\mu\otimes\nu)\\
&\geq
D_{}(\pi_G,\mu\otimes\nu),
\end{align*}
with equality if and only if $\pi=\pi_G$. Combining this with \eqref{eqn: c integral is same for pi and piG}, we obtain
\begin{align*}
\int c\,d\pi+D_{}(\pi,\mu\otimes\nu)
&=
\int c\,d\pi_G
+
D_{}(\pi,\pi_G)
+
D_{}(\pi_G,\mu\otimes\nu)\\
&\geq
\int c\,d\pi_G+D_{}(\pi_G,\mu\otimes\nu),
\end{align*}
with equality if and only if $\pi=\pi_G$. Thus, any non-Gaussian feasible
coupling with finite objective value can be strictly improved by its Gaussian
counterpart. Therefore, the unique optimizer of
\eqref{eq:problem-formulation-2} is a non--degenerate Gaussian coupling.
\end{itemize}
\end{enumerate}

\end{proof}
\subsection{Supplementary Analysis for Section \ref{section: cholesky Coordinates}}\label{sec: supplementary cholesky Coordinates}
\begin{proof}[Proof of Lemma~\ref{corollary: causality under cholesky transformation}]
Part (a) follows from the invertibility of $T_{\mathrm{chol}}$
and the fact that it acts separately on the two marginal variables.
For part (b), since the input covariance in the Cholesky coordinates is
$\Sigma^{\widetilde\mu}=I_{Td}$, we have $\Sigma^{\widetilde\mu}_{\widetilde U_{1:t},\widetilde U_{t+1}}=0.$
By Lemma~\ref{lemma: gaussian causal constraint},
$\widetilde\pi$ is causal if and only if
\begin{equation*}\Sigma^{\widetilde\pi}_{\widetilde Y_{1:t},\widetilde U_{t+1}}=0,
\qquad t=1,\ldots,T-1,
\end{equation*}
which is equivalent to
\eqref{eqn: causal constraint in cholesky coordinates}.

For part (c), let $\sigma(\cdot)$ denote the generated sigma-algebra.
Since $L^\mu$ and $L^\nu$ are invertible block lower-triangular matrices,
\[
\sigma(U_{1:t})=\sigma(\widetilde U_{1:t}),\quad
\sigma(Y_{1:t})=\sigma(\widetilde Y_{1:t}),\quad
\sigma(U_{1:T})=\sigma(\widetilde U_{1:T}).
\]
Since causality is equivalent to
\[
Y_{1:t}\perp\!\!\!\perp U_{1:T}\mid U_{1:t},
\qquad t=1,\ldots,T-1,
\]
these equalities imply that $\pi$ is causal if and only if
$\widetilde\pi$ is causal. This proves part (c).
\end{proof}

The following two lemmas establish the invariance of relative entropy under the Cholesky transformation and the equivalence of the original and transformed optimization problems.
\begin{lemma}
\label{lemma: relative entropy under cholesky transformation}
Let $\pi\in\Pi(\mu,\nu)$ be a non--degenerate Gaussian coupling. Then, under Assumption \ref{assumption:maringal}, 
\[
D(\pi,\mu\otimes\nu)
=
D(\widetilde\pi,\widetilde\mu\otimes\widetilde\nu).
\]
\end{lemma}
\begin{proof}[Proof of Lemma~\ref{lemma: relative entropy under cholesky transformation}]
    Since $\pi$, $\mu$, and $\nu$ are non--degenerate Gaussian distributions,
we have $\pi\ll\mu\otimes\nu$.
\begin{align*}
D(\pi,\mu\otimes\nu)
&=
\int
\log\frac{d\pi}{d(\mu\otimes\nu)}(\mathbf u,\mathbf y)\,
\pi(d\mathbf u,d\mathbf y)\\
&=
\int
\log\frac{d\widetilde\pi}
{d(\widetilde\mu\otimes\widetilde\nu)}
\bigl(T_{\mathrm{chol}}(\mathbf u,\mathbf y)\bigr)\,
\pi(d\mathbf u,d\mathbf y)\\
&=
\int
\log\frac{d\widetilde\pi}
{d(\widetilde\mu\otimes\widetilde\nu)}
(\widetilde{\mathbf u},\widetilde{\mathbf y})\,
\widetilde\pi(d\widetilde{\mathbf u},d\widetilde{\mathbf y})\\
&=
D(\widetilde\pi,\widetilde\mu\otimes\widetilde\nu),
\end{align*}
where the second equality uses
\[
\frac{d\widetilde\pi}
{d(\widetilde\mu\otimes\widetilde\nu)}
(\widetilde{\mathbf u},\widetilde{\mathbf y})
=
\frac{d\pi}{d(\mu\otimes\nu)}
\left(
T_{\mathrm{chol}}^{-1}
(\widetilde{\mathbf u},\widetilde{\mathbf y})
\right),
\]
which follows since, for any measurable set $A$, using the definition of the pushforward,
\begin{align*}
\widetilde\pi(A)
&=
\pi(T_{\mathrm{chol}}^{-1}(A))\\
&=
\int_{T_{\mathrm{chol}}^{-1}(A)}
\frac{d\pi}{d(\mu\otimes\nu)}(\mathbf u,\mathbf y)\,
(\mu\otimes\nu)(d\mathbf u,d\mathbf y)\\
&=
\int_A
\frac{d\pi}{d(\mu\otimes\nu)}
\left(T_{\mathrm{chol}}^{-1}
(\widetilde{\mathbf u},\widetilde{\mathbf y})\right)
(\widetilde\mu\otimes\widetilde\nu)
(d\widetilde{\mathbf u},d\widetilde{\mathbf y}).
\end{align*}
\end{proof}
\begin{lemma} \label{lemma: pi start is minimizer of original problem if and ony if its pushfowrad by T chol is optimizer of the problem in chol coordinate}
Under Assumptions~\ref{assumption:cost}-\ref{assumption:maringal}, we have 
\[
\pi^*\in
\arg\min_{\pi\in\Pi_c(\mu,\nu)}
\left\{
\int c\,d\pi+D(\pi,\mu\otimes\nu)
\right\}
\]
if and only if
\[
(T_{\mathrm{chol}})_\#\pi^*
\in
\arg\min_{\widetilde\pi\in\Pi_c(\widetilde\mu,\widetilde\nu)}
\left\{
\int \widetilde c\,d\widetilde\pi
+
D(\widetilde\pi,\widetilde\mu\otimes\widetilde\nu)
\right\}.
\]
\end{lemma}
\begin{proof}[Proof of Lemma~\ref{lemma: pi start is minimizer of original problem if and ony if its pushfowrad by T chol is optimizer of the problem in chol coordinate}]

By Lemma~\ref{theorem: optomality of gaussian}, it is sufficient to consider
non--degenerate Gaussian couplings. Since
$\widetilde c$, $\widetilde\mu$, and $\widetilde\nu$ satisfy
Assumptions~\ref{assumption:cost}-\ref{assumption:maringal}, the same result
also applies to the transformed problem. For such couplings,
Lemma~\ref{corollary: causality under cholesky transformation} give
\[
\pi\in\Pi_c(\mu,\nu)
\quad\Longleftrightarrow\quad
\widetilde\pi\in\Pi_c(\widetilde\mu,\widetilde\nu).
\]
Since \(T_{\mathrm{chol}}\) is invertible, this transformation gives a
one-to-one correspondence between the feasible Gaussian couplings in the
two coordinate systems.
Moreover, by the definition of the pushforward and Lemma \ref{lemma: relative entropy under cholesky transformation}, for every corresponding pair $\pi,\widetilde \pi$, they have the same objective value
\begin{align*}
    \int c(\mathbf u,\mathbf y)\,\pi(d\mathbf u,d\mathbf y)+D(\pi,\mu\otimes\nu)
&=
\int c(L^\mu\widetilde{\mathbf u},L^\nu\widetilde{\mathbf y})\,
\widetilde\pi(d\widetilde{\mathbf u},d\widetilde{\mathbf y})+D(\widetilde\pi,\widetilde\mu\otimes\widetilde\nu)\\&= \int \widetilde c(\widetilde{\mathbf u},\widetilde{\mathbf y})\,
\widetilde\pi(d\widetilde{\mathbf u},d\widetilde{\mathbf y})+D(\widetilde\pi,\widetilde\mu\otimes\widetilde\nu).
\end{align*}
where we used  the definition of \(\widetilde c\). This concludes the result.
\end{proof}

\section{Supplementary Analysis for Section \ref{section: gaussian dual problem}}
\subsection{Supplementary Analysis for Section  \ref{section: dual formulation}}\label{sec: supp for dual formulation}
The following lemma establishes that \(Z\) is affine and injective, which will be used in the proof of Theorem~\ref{theorem: dual problem of ecot gaussian} to establish the strict concavity of the dual objective.
\begin{lemma}
\label{lemma: Z injective}
The map \(Z:\Theta\to\mathbb S^{2Td}\) defined in
\eqref{eqn: Z from theta} is affine and injective.
\end{lemma}
\begin{proof}[Proof of Lemma~\ref{lemma: Z injective}]
For any \(\theta_1,\theta_2\in\Theta\) and \(\lambda\in\mathbb R\), the linear dependence of
\(\Gamma_\Psi\) on \(\Psi\) gives
\[
Z\bigl(\lambda\theta_1+(1-\lambda)\theta_2\bigr)
=
\lambda Z(\theta_1)+(1-\lambda)Z(\theta_2),
\]
because the constant term \(Q\) is unchanged. Hence, \(Z\) is affine.

To see why it is injective, suppose
\[
Z(\theta)=Z(\theta'),
\]
where
\[
\theta=(\Lambda_u,\Lambda_y,\Psi),
\qquad
\theta'=(\Lambda_u',\Lambda_y',\Psi').
\]
Since \(Q\) is fixed, \eqref{eqn: Z from theta} gives
\[
\begin{bmatrix}
2(\Lambda_u-\Lambda_u')
&
\Gamma_\Psi-\Gamma_{\Psi'}
\\
(\Gamma_\Psi-\Gamma_{\Psi'})^\top
&
2(\Lambda_y-\Lambda_y')
\end{bmatrix}
=0.
\]
The diagonal blocks therefore give
\[
\Lambda_u=\Lambda_u',
\qquad
\Lambda_y=\Lambda_y',
\]
while the off-diagonal block gives
\[
\Gamma_\Psi=\Gamma_{\Psi'}.
\]
Upon using the definition of $E_{t+1}$ and $S_t$, 
the structural layout of \(\Gamma_\Psi\) is given by
\begin{equation}
    \label{eqn: def of Gamma Psi from Psi t s}
\Gamma_\Psi
=
\left[
\begin{array}{ccccc}
0 & 0 & 0 & \cdots & 0\\[2mm]
\multicolumn{1}{c}{
\underbrace{\makebox[1.5cm][c]{$\Psi_1^\top$}}_{d\times d}}
& 0 & 0 & \cdots & 0\\[3mm]
\multicolumn{2}{c}{
\underbrace{\makebox[3cm][c]{$\Psi_2^\top$}}_{d\times 2d}}
& 0 & \cdots & 0\\[3mm]
\ldots & \ldots & \ldots & \ddots & \vdots\\[2mm]
\multicolumn{4}{c}{
\underbrace{\makebox[6cm][c]{$\Psi_{T-1}^\top$}}_{d\times (T-1)d}}
& 0
\end{array}
\right].
\end{equation}
Here, \(\Psi_t^\top\in\mathbb R^{d\times td}\) occupies the first \(t\)
block columns of the \((t+1)\)-st block row. 
By the block-row representation of \(\Gamma_\Psi\), for every
\(t=1,\ldots,T-1\),
\[
\begin{bmatrix}
\Psi_t^\top & 0_{d\times(T-t)d}
\end{bmatrix}
=
\begin{bmatrix}
(\Psi_t')^\top & 0_{d\times(T-t)d}
\end{bmatrix}.
\]
Hence,
\[
\Psi_t=\Psi_t',
\qquad
t=1,\ldots,T-1.
\]
Thus, \(\theta=\theta'\), proving that \(Z\) is injective.
\end{proof}

\begin{lemma}
\label{lemma: theta star open convex}
\(\Thetafeasible\)  is open and convex.
\end{lemma}

\begin{proof}[Proof of Lemma~\ref{lemma: theta star open convex}]
Since \(Z:\Theta\to\mathbb S^{2Td}\) is continuous and
\(\mathbb S_{++}^{2Td}\) is open,
\[
\Thetafeasible
=
Z^{-1}\!\left(\mathbb S_{++}^{2Td}\right)
\]
is open.

Moreover, let \(\theta_1,\theta_2\in\Thetafeasible\) and \(t\in[0,1]\).
Since \(Z\) is affine,
\[
Z\bigl((1-t)\theta_1+t\theta_2\bigr)
=
(1-t)Z(\theta_1)+tZ(\theta_2).
\]
We have
\[
Z(\theta_1)\succ0,
\qquad
Z(\theta_2)\succ0.
\]
Since the cone \(\mathbb S_{++}^{2Td}\) is convex,
\[
(1-t)Z(\theta_1)+tZ(\theta_2)\succ0.
\]
Therefore,
\[
(1-t)\theta_1+t\theta_2\in\Thetafeasible,
\]
and hence \(\Thetafeasible\) is convex.
\end{proof}

\begin{proof}[Proof of Theorem~\ref{theorem: dual problem of ecot gaussian}]
We use the dual variables to introduce the Lagrangian
\begin{align*}
&\mathcal L(\Sigma;\theta)\numberthis\label{eqn: lagrangian of cov optimizatio problem} =\mathcal J(\Sigma) \\&\quad + \trace(\Lambda_u(P_u\Sigma P_u^\top-I_{Td})) + \trace(\Lambda_y(P_y\Sigma P_y^\top-I_{Td})) \\&\quad+ \sum_{t=1}^{T-1} \trace(\Psi_t^\top S_t P_y \Sigma P_u^\top E_{t+1})\\
&\quad = \frac12 \trace(\Sigma Z(\theta)) -  \frac12 \log\det{\Sigma} - \trace(\Lambda_u) - \trace(\Lambda_y) .
\end{align*}
\begin{enumerate}[label=\alph*.]
    \item
    We first establish strong duality and derive the dual objective.
\eqref{eqn:problem-over-cov} is a convex optimization problem, since all constraints in \eqref{eqn:problem-over-cov} are affine and \(\mathcal J\) is convex on its effective domain, as \(-\log\det\Sigma\) is convex on \(\mathbb S_{++}^{2Td}\) and \(\operatorname{Tr}(Q\Sigma)\) is linear in \(\Sigma\).
Therefore, Slater's
condition can be used to establish strong duality and dual attainment.
Consider
    \[
    \Sigma=I_{2Td},
    \]
    which satisfies the marginal constraints in \eqref{eqn:maringal-constraint-mu},
\eqref{eqn:maringal-constraint-nu}. It also satisfies the causal constraint in \eqref{eqn:causal-constraint}, since all of its \(Y\)-\(U\) cross-covariance blocks vanish. Since \(I_{2Td}\succ0\), this
    feasible point lies in the interior of the effective domain of
    \(\mathcal J\). Therefore, Slater's condition holds. Consequently,
    strong duality holds, the dual optimum is attained, and the KKT
    conditions are necessary and sufficient for optimality.

    It remains to identify the corresponding dual function. By the definition
of the Lagrangian above, for fixed \(\theta\in\Theta\), consider the inner problem
\[
\inf_{\Sigma\in\mathbb S^{2Td}}
\mathcal L(\Sigma;\theta)
=
-\operatorname{Tr}(\Lambda_u)
-\operatorname{Tr}(\Lambda_y)
+
\inf_{\Sigma\succ0}
\left\{
\frac12\operatorname{Tr}\!\left(
Z(\theta)\Sigma
\right)
-
\frac12\log\det\Sigma
\right\}.
\]
    If \(Z(\theta)\succ0\), this objective is strictly convex in
    \(\Sigma\), and its stationary condition is
    \begin{equation}
    \label{eqn: stationary condition for Sigma in dual}
    Z(\theta)-\Sigma^{-1}=0.
    \end{equation}
    Hence, the unique minimizer is
    \[
    \Sigma=Z(\theta)^{-1}.
    \]
    Substituting this minimizer gives
    \begin{align*}
    &\inf_{\Sigma\succ0}
    \left\{
    \frac12\operatorname{Tr}\!\left(
    Z(\theta)\Sigma
    \right)
    -
    \frac12\log\det\Sigma
    \right\}\\
    &=
    Td+\frac12\log\det Z(\theta).
    \end{align*}
    If \(Z(\theta)\not\succ0\), the inner infimum is \(-\infty\).\\
Therefore, the dual function
\[
\theta\longmapsto
\inf_{\Sigma\in\mathbb S^{2Td}}
\mathcal L(\Sigma;\theta)
\]
is exactly the function \(\mathcal D\) defined in
\eqref{eqn: dual objective}. Hence, maximizing the dual function over
\(\theta\in\Theta\) gives problem \eqref{eqn: dual problem}, proving that
\eqref{eqn: dual problem} is the Lagrangian dual of
\eqref{eqn:problem-over-cov}.
    \item
For any distinct \(\theta_1,\theta_2\in\Thetafeasible\) and
\(\lambda\in(0,1)\), Lemma \ref{lemma: Z injective} gives
\(Z(\theta_1)\neq Z(\theta_2)\). Therefore,
\begin{align*}
&\log\det Z\!\left(\lambda\theta_1+(1-\lambda)\theta_2\right)\\
&=
\log\det\!\left(
\lambda Z(\theta_1)+(1-\lambda)Z(\theta_2)
\right)\\
&>
\lambda\log\det Z(\theta_1)
+
(1-\lambda)\log\det Z(\theta_2),
\end{align*}
where the equality follows from the affinity of \(Z\), and the strict
inequality follows from the strict concavity of \(\log\det\) together with
\(Z(\theta_1)\neq Z(\theta_2)\). Hence, $    \theta\longmapsto\log\det Z(\theta)$ is strictly concave on \(\Thetafeasible\). The remaining terms $ -\operatorname{Tr}(\Lambda_u)
    -
    \operatorname{Tr}(\Lambda_y)$
    are linear in \(\theta\). Hence, \(\mathcal D\) is strictly concave
    on \(\Thetafeasible\).

    \item
    By part (a), the dual optimum is attained and, by strong duality, its
optimal value is finite. Hence, every dual optimizer belongs to \(\Thetafeasible\). The strict
    concavity established in part (b) therefore implies that the dual
    optimizer \(\theta^*\) is unique.

    \item
    By the KKT conditions established in part (a), the primal-dual
    optimizers satisfy the stationary condition
    \eqref{eqn: stationary condition for Sigma in dual}. Therefore,
    \[
    Z(\theta^*)-(\Sigma^*)^{-1}=0,
    \]
    and hence
    \[
    \Sigma^*=Z(\theta^*)^{-1}.
    \]
\end{enumerate}
\end{proof}
\subsection{Supplementary Analysis for Section  \ref{section:bcd}} \label{section: supp bcd}
\begin{lemma}
\label{lemma: block inverse formula}
Let
\[
M=
\begin{bmatrix}
M_{uu} & M_{uy}\\
M_{yu} & M_{yy}
\end{bmatrix}.
\]
If \(M\), \(M_{uu}\), and \(M_{yy}\) are invertible, then
\begin{align}
[M^{-1}]_{uu}
&=
\left(
M_{uu}-M_{uy}M_{yy}^{-1}M_{yu}
\right)^{-1},
\label{eqn: block inverse uu}\\
[M^{-1}]_{yy}
&=
\left(
M_{yy}-M_{yu}M_{uu}^{-1}M_{uy}
\right)^{-1},
\label{eqn: block inverse yy}\\
[M^{-1}]_{yu}
&=
-M_{yy}^{-1}M_{yu}[M^{-1}]_{uu}.
\label{eqn: block inverse yu}
\end{align}
\end{lemma}

\begin{proof}[Proof of Lemma~\ref{lemma: block inverse formula}]
By the inversion formula in
\cite[Eq.~(6.0.8)]{zhang2005schur},
{\footnotesize
\begin{align*}  
&M^{-1}
=\\
&\begin{bmatrix}
(M_{uu}-M_{uy}M_{yy}^{-1}M_{yu})^{-1}
&
-(M_{uu}-M_{uy}M_{yy}^{-1}M_{yu})^{-1}
M_{uy}M_{yy}^{-1}
\\[2mm]
-M_{yy}^{-1}M_{yu}
(M_{uu}-M_{uy}M_{yy}^{-1}M_{yu})^{-1}
&
M_{yy}^{-1}
+
M_{yy}^{-1}M_{yu}
(M_{uu}-M_{uy}M_{yy}^{-1}M_{yu})^{-1}
M_{uy}M_{yy}^{-1}
\end{bmatrix}.
\end{align*}
}
Reading the top--left block gives
\eqref{eqn: block inverse uu}. Reading the bottom--left block and using
\eqref{eqn: block inverse uu} gives
\eqref{eqn: block inverse yu}.
Similarly, by the inversion formula in
\cite[Eq.~(6.0.7)]{zhang2005schur},
the bottom--right block gives
\eqref{eqn: block inverse yy}.
\end{proof}

\begin{lemma}
\label{lemma: gradient of log det of linear map}
Let
\[
X\in\mathbb R^{p\times q},
\]
and suppose
\[
M(X)=C+AXB+GX^\top E\in\mathbb R^{m\times m},
\]
where
\[
A\in\mathbb R^{m\times p},\qquad
B\in\mathbb R^{q\times m},\qquad
G\in\mathbb R^{m\times q},\qquad
E\in\mathbb R^{p\times m},
\]
and \(C\in\mathbb R^{m\times m}\). Assume that \(M(X)\succ0\), and let
\[
h(X)=g(M(X)),
\qquad
g(M)=-\frac12\log\det M.
\]
Then
\begin{equation}
\nabla h(X)
=
-\frac12\left(
A^\top M(X)^{-1}B^\top
+
EM(X)^{-1}G
\right).
\label{eqn: graident of log det of liner map}
\end{equation}
\end{lemma}
\begin{proof}[Proof of Lemma~\ref{lemma: gradient of log det of linear map}]
Since \(M\) is affine, for any \(U\in\mathbb R^{p\times q}\),
\[
DM(X)[U]
=
AUB+GU^\top E.
\]
By Theorem 8.2 of \cite{magnus2019matrix},
\[
d\log\det M
=
\operatorname{Tr}\!\left(M^{-1}dM\right),
\]
so, equivalently,
\[
Dg(M)[V]
=
-\frac12\operatorname{Tr}\!\left(M^{-1}V\right).
\]
Since \(h=g\circ M\), Theorem 5.12 of \cite{magnus2019matrix} gives
\[
Dh(X)
=
Dg(M(X))\,DM(X).
\]
Thus, for any \(U\in\mathbb R^{p\times q}\),
\begin{align*}
Dh(X)[U]
&=
Dg(M(X))\!\left[DM(X)[U]\right]\\
&=
-\frac12
\operatorname{Tr}\!\left(
M(X)^{-1}\left(AUB+GU^\top E\right)
\right)\\
&=
-\frac12
\operatorname{Tr}\!\left(
M(X)^{-1}AUB
\right)
-\frac12
\operatorname{Tr}\!\left(
M(X)^{-1}GU^\top E
\right).
\end{align*}
Using cyclicity of the trace,
\begin{align*}
\operatorname{Tr}\!\left(
M(X)^{-1}AUB
\right)
&=
\operatorname{Tr}\!\left(
BM(X)^{-1}AU
\right)\\
&=
\operatorname{Tr}\!\left(
\left(A^\top M(X)^{-\top}B^\top\right)^\top U
\right),
\end{align*}
and
\begin{align*}
\operatorname{Tr}\!\left(
M(X)^{-1}GU^\top E
\right)
&=
\operatorname{Tr}\!\left(
EM(X)^{-1}GU^\top
\right)\\
&=
\operatorname{Tr}\!\left(
\left(EM(X)^{-1}G\right)^\top U
\right).
\end{align*}
Therefore,
\[
Dh(X)[U]
=
\operatorname{Tr}\!\left(
\left[
-\frac12 A^\top M(X)^{-\top}B^\top
-\frac12 EM(X)^{-1}G
\right]^\top U
\right).
\]
Since
\[
Dh(X)[U]
=
\operatorname{Tr}\!\left(
\nabla h(X)^\top U
\right),
\]
we obtain
\[
\nabla h(X)
=
-\frac12 A^\top M(X)^{-\top}B^\top
-\frac12 EM(X)^{-1}G.
\]
Since \(M(X)\succ0\), we have \(M(X)^{-\top}=M(X)^{-1}\). Therefore,
\[
\nabla h(X)
=
-\frac12\left(
A^\top M(X)^{-1}B^\top
+
EM(X)^{-1}G
\right),
\]
which proves \eqref{eqn: graident of log det of liner map}.
\end{proof}

\begin{proof}[Proof of Theorem~\ref{theorem: stationary condition update-upsi}]
We construct a feasible stationary point of the block subproblem by first
solving the stationary conditions for \(\Psi\) and \(\Lambda_u\), and then
verifying that the resulting parameters remain in \(\Thetafeasible\).
The conditions are
\[
\nabla_{\Lambda_u}f=0,
\qquad
\nabla_{\Psi_t}f=0,
\quad
t=1,\ldots,T-1.
\]
By \eqref{eqn: dual objective}, the objective on $ \theta \in \Thetafeasible$ is
\[
f(\theta)
=
-\frac12\log\det Z(\theta)
+
\operatorname{Tr}(\Lambda_u)
+
\operatorname{Tr}(\Lambda_y)
-
Td.
\]
We first apply Lemma \ref{lemma: gradient of log det of linear map} to compute
the derivative of the log-determinant term in \(f(\theta)\) with respect to
\(\Lambda_u\): 
\[
\nabla_{\Lambda_u}
\left(
-\frac12\log\det Z(\theta)
\right)
=
-P_uZ(\theta)^{-1}P_u^\top.
\]
where we used $Z(\theta)
=
C+2P_u^\top \Lambda_u P_u$
with \(C\) collecting the terms independent of \(\Lambda_u\).
Moreover, 
\[
\nabla_{\Lambda_u}\operatorname{Tr}(\Lambda_u)=I_{Td},
\]
concluding,
\[
\nabla_{\Lambda_u}f(\theta)
=
-P_uZ(\theta)^{-1}P_u^\top+I_{Td}.
\]
As a result, the stationary condition
\(\nabla_{\Lambda_u}f=0\) yields
\begin{equation}
    [Z(\theta)^{-1}]_{uu}=I_{Td}. \label{eqn: in bcd for u psi block the Z inv uu block is maringal I}
\end{equation}
where we used $
P_uZ(\theta)^{-1}P_u^\top=[Z(\theta)^{-1}]_{uu}$. Finally, we use this identity and~\eqref{eqn: block inverse yu} in Lemma \ref{lemma: block inverse formula} applied to $Z(\theta)$ to conclude \begin{equation}
[Z(\theta)^{-1}]_{yu}
=
-[Z_{yy}]^{-1}[Z(\theta)]_{yu}.
\label{eqn: Z inv yu from the bcd stationary condition on Lambda_u}
\end{equation}
This identity is used next to simplify the stationary conditions for
\(\{\Psi_t\}_{t=1}^{T-1}\).


\begin{enumerate}[label=\alph*.]

\item 
Applying Lemma \ref{lemma: gradient of log det of linear map} to compute
the derivative of the log-determinant term in \(f(\theta)\) with respect to
\(\Psi_t\) yields
\begin{align*}
\nabla_{\Psi_t}
\left(
-\frac12\log\det Z(\theta)
\right)
&=
-\frac12\left(
S_tP_yZ(\theta)^{-1}P_u^\top E_{t+1}
+
S_tP_yZ(\theta)^{-1}P_u^\top E_{t+1}
\right)\\
&=
-
S_t[Z(\theta)^{-1}]_{yu}E_{t+1}.
\end{align*}
where we used
\[
Z(\theta) = P_y^\top S_t^\top \Psi_t E_{t+1}^\top P_u
+
P_u^\top E_{t+1}\Psi_t^\top S_tP_y + \text{(const. in $\Psi_t$)}.
\]
Since the remaining terms in \(f\) do not depend on \(\Psi_t\),
\[
\nabla_{\Psi_t}f
=
-
S_t[Z(\theta)^{-1}]_{yu}E_{t+1},
\]
concluding the stationary condition 
\[
S_t[Z(\theta)^{-1}]_{yu}E_{t+1}=0.
\]
Using
\eqref{eqn: Z inv yu from the bcd stationary condition on Lambda_u},
this becomes
\[
S_t[Z_{yy}]^{-1}[Z(\theta)]_{yu}E_{t+1}=0.
\]

By the definition of \(Z(\theta)\),
\[
[Z(\theta)]_{yu}
=
Q_{yu}+\Gamma_\Psi^\top,
\]
where
\[
\Gamma_\Psi^\top
=
\sum_{s=1}^{T-1}
S_s^\top\Psi_sE_{s+1}^\top.
\]
Since the matrices \(E_{s+1}\) select distinct block columns,
\[
E_{s+1}^\top E_{t+1}
=
\begin{cases}
I_d,&s=t,\\
0,&s\neq t.
\end{cases}
\]
Therefore,
\[
\Gamma_\Psi^\top E_{t+1}
=
S_t^\top\Psi_t.
\]
Thus, the stationary condition corresponding to time \(t\) depends only on
\(\Psi_t\)and reduces to
\[
S_t[Z_{yy}]^{-1}
\left(
Q_{yu}E_{t+1}
+
S_t^\top\Psi_t
\right)
=0,
\]
or equivalently,
\[
\left(
S_t[Z_{yy}]^{-1}S_t^\top
\right)\Psi_t
=
-
S_t[Z_{yy}]^{-1}Q_{yu}E_{t+1}.
\]

We have \(S_t^\top x\neq0\) whenever \(x\neq0\), because, for any
\(x\in\mathbb R^{td}\),
\(
S_t^\top x
=
\begin{bmatrix}
x\\
0
\end{bmatrix}.
\) Moreover, since \(Z_{yy}^{-1}\succ0\), for any \(x\in\mathbb R^{td}\setminus\{0\}\),
\[
x^\top S_t[Z_{yy}]^{-1}S_t^\top x
=
(S_t^\top x)^\top [Z_{yy}]^{-1}(S_t^\top x)
>0.
\]
Hence,
\[
S_t[Z_{yy}]^{-1}S_t^\top\succ0,
\]
and therefore it is invertible.
Thus,
\[
\Psi_t^{\{k+1\}}
=
-
\left(
S_t[Z_{yy}]^{-1}S_t^\top
\right)^{-1}
S_t[Z_{yy}]^{-1}Q_{yu}E_{t+1},
\]
which is exactly \eqref{eqn: schrodinger eqn 3}. Hence, each
\(\Psi_t^{\{k+1\}}\) is uniquely determined. 
\item
Once \(\Psi^{\{k+1\}}\) has been determined, the stationary condition
\[
[Z(\theta)^{-1}]_{uu}=I_{Td}
\]
determines \(\Lambda_u\). 
By Lemma \ref{lemma: block inverse formula},
\[
[Z(\theta)^{-1}]_{uu}
=
\left(
[Z(\theta)]_{uu}
-
[Z(\theta)]_{uy}Z_{yy}^{-1}[Z(\theta)]_{yu}
\right)^{-1}.
\]
Therefore,
\[
[Z(\theta)]_{uu}
-
[Z(\theta)]_{uy}Z_{yy}^{-1}[Z(\theta)]_{yu}
=
I_{Td}.
\]
By part (a), any feasible stationary point must satisfy
\(\Psi=\Psi^{\{k+1\}}\), and therefore
\[
[Z(\theta)]_{uy}
=
Z_{uy},
\qquad
[Z(\theta)]_{yu}
=
Z_{uy}^\top,
\]
where
\[
Z_{uy}
=
Q_{uy}+\Gamma_{\Psi^{\{k+1\}}}.
\]
Moreover,
\[
[Z(\theta)]_{uu}
=
Q_{uu}+2\Lambda_u.
\]
Therefore, the stationary condition uniquely determines
\(\Lambda_u\) as
\[
2\Lambda_u^{\{k+1\}}
=
-Q_{uu}
+
I_{Td}
+
Z_{uy}Z_{yy}^{-1}Z_{uy}^\top,
\]
which gives \eqref{eqn: schrodinger eqn 1}.
\item
It remains to verify that the point constructed in parts (a) and (b)
belongs to \(\Thetafeasible\). For
\[
\theta^{\{k+\frac12\}}
=
(\Lambda_u^{\{k+1\}},
\Lambda_y^{\{k\}},
\Psi^{\{k+1\}}),
\]
we have
\[
[Z(\theta^{\{k+\frac12\}})]_{yy}=Z_{yy},
\qquad
[Z(\theta^{\{k+\frac12\}})]_{uy}=Z_{uy},
\qquad
[Z(\theta^{\{k+\frac12\}})]_{yu}=Z_{uy}^\top.
\]
Moreover, by \eqref{eqn: schrodinger eqn 1},
\[
[Z(\theta^{\{k+\frac12\}})]_{uu}
=
Q_{uu}+2\Lambda_u^{\{k+1\}}
=
I_{Td}
+
Z_{uy}Z_{yy}^{-1}Z_{uy}^\top.
\]
Therefore, the Schur complement of \(Z_{yy}\) in
\(Z(\theta^{\{k+\frac12\}})\) satisfies
\[
[Z(\theta^{\{k+\frac12\}})]_{uu}
-
Z_{uy}Z_{yy}^{-1}Z_{uy}^\top
=
I_{Td}\succ0.
\]
Together with \(Z_{yy}\succ0\), the Schur complement characterization
therefore gives
\[
Z(\theta^{\{k+\frac12\}})\succ0.
\]
Thus,
\[
\theta^{\{k+\frac12\}}\in\Thetafeasible.
\]
\end{enumerate}
Finally, by Theorem~\ref{theorem: dual problem of ecot gaussian},
\(f=-\mathcal D\) is strictly convex on \(\Thetafeasible\), and hence its
restriction to the \((\Lambda_u,\Psi)\)-block is strictly convex on its
feasible domain. Since the point constructed above is feasible and satisfies
all stationary conditions of the block subproblem, convexity implies that it
is a minimizer. Strict convexity then implies that this minimizer is unique.
\end{proof}

\begin{proof}[Proof of Theorem~\ref{theorem: stationary condition of update-lambday}]
We construct a feasible stationary point of the block subproblem by solving
the stationary condition for \(\Lambda_y\) and then verifying that the
resulting parameters remain in \(\Thetafeasible\). The stationary condition is
\[
\nabla_{\Lambda_y}f=0.
\]

We apply Lemma \ref{lemma: gradient of log det of linear map} to compute
the derivative of the log-determinant term in \(f(\theta)\) with respect to
\(\Lambda_y\):
\[
\nabla_{\Lambda_y}
\left(
-\frac12\log\det Z(\theta)
\right)
=
-P_yZ(\theta)^{-1}P_y^\top.
\]
where we used 
\[
Z(\theta)
=
C+2P_y^\top\Lambda_yP_y,
\]
with \(C\) collecting the terms independent of \(\Lambda_y\).
Moreover,
\[
\nabla_{\Lambda_y}\operatorname{Tr}(\Lambda_y)=I_{Td},
\]
while the remaining terms in \(f\) do not depend on \(\Lambda_y\).
Hence,
\[
\nabla_{\Lambda_y}f
=
-P_yZ(\theta)^{-1}P_y^\top+I_{Td}.
\]
Therefore, the stationary condition
\(\nabla_{\Lambda_y}f=0\) gives
\[
P_yZ(\theta)^{-1}P_y^\top=I_{Td},
\]
or equivalently,
\begin{equation}
[Z(\theta)^{-1}]_{yy}=I_{Td}.
\label{eqn: in bcd for y block the Z inv yy block is marginal I}
\end{equation}

Since
\(\theta^{\{k+\frac12\}}\in\Thetafeasible\),
we have
\[
Z(\theta^{\{k+\frac12\}})\succ0.
\]
Hence its principal \(uu\)-block is positive definite, so
\[
Z_{uu}
=
Q_{uu}+2\Lambda_u^{\{k+1\}}
\succ0.
\]
Moreover,
\[
Z_{yu}
=
Q_{yu}+\Gamma_{\Psi^{\{k+1\}}}^{\top},
\qquad
Z_{uy}=Z_{yu}^{\top}.
\]
The matrices \(Z_{uu}\), \(Z_{yu}\), and \(Z_{uy}\) depend only on
\(\Lambda_u^{\{k+1\}}\) and \(\Psi^{\{k+1\}}\), which are fixed throughout
the \(\Lambda_y\)-block subproblem, and therefore remain unchanged during
this update.
\begin{enumerate}[label=\alph*.]
\item 
To solve the stationary condition explicitly for \(\Lambda_y\), by
Lemma \ref{lemma: block inverse formula},
\[
[Z(\theta)^{-1}]_{yy}
=
\left(
[Z(\theta)]_{yy}
-
Z_{yu}Z_{uu}^{-1}Z_{uy}
\right)^{-1}.
\]
Together with
\eqref{eqn: in bcd for y block the Z inv yy block is marginal I},
this gives
\[
\left(
[Z(\theta)]_{yy}
-
Z_{yu}Z_{uu}^{-1}Z_{uy}
\right)^{-1}
=
I_{Td}.
\]
Therefore,
\[
[Z(\theta)]_{yy}
-
Z_{yu}Z_{uu}^{-1}Z_{uy}
=
I_{Td}.
\]
Since
\[
[Z(\theta)]_{yy}
=
Q_{yy}+2\Lambda_y,
\]
the stationary condition uniquely determines \(\Lambda_y\) as
\[
2\Lambda_y^{\{k+1\}}
=
-Q_{yy}
+
I_{Td}
+
Z_{yu}Z_{uu}^{-1}Z_{uy}.
\]
We obtain \eqref{eqn: schrodinger eqn 2} by $Z_{uy}=Z_{yu}^\top.$
\item
It remains to verify that the updated parameters belong to \(\Thetafeasible\).
For
\[
\theta^{\{k+1\}}
=
(\Lambda_u^{\{k+1\}},
\Lambda_y^{\{k+1\}},
\Psi^{\{k+1\}}),
\]
we have
\[
[Z(\theta^{\{k+1\}})]_{yy}
=
Q_{yy}+2\Lambda_y^{\{k+1\}}.
\]
By \eqref{eqn: schrodinger eqn 2},
\[
[Z(\theta^{\{k+1\}})]_{yy}
-
Z_{yu}Z_{uu}^{-1}Z_{uy}
=
I_{Td}\succ0.
\]
Since \(Z_{uu}\succ0\), the Schur complement characterization gives
\[
Z(\theta^{\{k+1\}})\succ0.
\]
Thus,
\[
\theta^{\{k+1\}}\in\Thetafeasible.
\]

\end{enumerate}

Finally, by Theorem~\ref{theorem: dual problem of ecot gaussian},
\(f\) is strictly convex on \(\Thetafeasible\), and hence its restriction to the
\(\Lambda_y\)-block is strictly convex on its feasible domain. Since the
point constructed above is feasible and satisfies the stationary condition
of the block subproblem, convexity implies that it is a minimizer. Strict
convexity then implies that this minimizer is unique.
\end{proof}
\section{Supplementary Analysis for Section~\ref{sec: linear convergence}}\label{sec: supp for linear convergence}

\subsection{Supplementary Analysis for Section  \ref{section: smooth and str convx over compact}}\label{section: supp smooth and str convx over compact}
\begin{proof}[Proof of Lemma~\ref{lemma: hessian of f}]
For any $\theta \in \Thetafeasible$, express 
\begin{align*}
    f(\theta)= \frac{1}{2}g(Z(\theta)) + \trace(\Lambda_u) + \trace(\Lambda_y) - Td
\end{align*}
where 
\[
g(Z):=-\log\det Z,
\qquad
Z\in\mathbb S_{++}^{2Td}.
\]
Therefore,
\begin{align*}
    \nabla^2 f(\theta)
[\delta\theta_1,\delta\theta_2]
&=
\frac12
\nabla^2 g(Z(\theta))
\left[
A (\delta\theta_1),
A(\delta\theta_2)
\right]
\end{align*}
where we used the affinity of $Z$ and the definition of $A$
in \eqref{eqn: A from theta}. Using the formula for the second derivative of log-determinant~\cite[Appendix~A.4.3]{boyd2004convex}, 
\[
\nabla^2 g(Z)[\delta Z_1,\delta Z_2]
=
\operatorname{Tr}
\left(
Z^{-1}\delta Z_1 Z^{-1}\delta Z_2
\right),
\qquad
\delta Z_1,\delta Z_2\in\mathbb S^{2Td};
\]
we arrive at \eqref{eqn: hessian of f}. 
The strict positivity of the Hessian~\eqref{eqn: strictly positive hessian of f} follows from the fact that 
\[
Z(\theta)^{-1/2}
A(\delta\theta)
Z(\theta)^{-1/2}
\neq0,\quad \forall \delta\theta\neq0
\]
because \(A(\delta\theta)\neq0\) by Lemma \ref{lemma: Z injective}, and
\(Z(\theta)^{-1/2}\) is invertible on $\Thetafeasible$.
\end{proof}

\begin{lemma}
\label{lemma: explicit bounds for A}
The linear map \(A\) defined in \eqref{eqn: A from theta} satisfies the following bounds:
\begin{align}
\|A(h_1,0)\|_F
&\leq
2\|h_1\|_{\Theta_1},
\qquad
\forall h_1\in\Theta_1,
\label{eqn: A first block upper bound}\\
\|A(0,h_2)\|_F
&=
2\|h_2\|_{\Theta_2},
\qquad
\forall h_2\in\Theta_2,
\label{eqn: A second block upper bound}\\
\|A(\theta)\|_F
&\geq\sqrt{2}
\|\theta\|_\Theta,
\qquad
\forall \theta\in\Theta.
\label{eqn: A full space lower bound}
\end{align}
\end{lemma}

\begin{proof}[Proof of Lemma~\ref{lemma: explicit bounds for A}]
Throughout the proof, we use  
\begin{equation}
\|\Gamma_\Psi\|_F^2
=
\sum_{t=1}^{T-1}\|\Psi_t\|_F^2.
\label{eqn: Gamma Psi Frobenius norm}
\end{equation}
which is obtained from the identity~\eqref{eqn: def of Gamma Psi from Psi t s}.
\begin{enumerate}[label=\alph*.]
    \item Proof of \eqref{eqn: A first block upper bound}:
    Taking $h_1=(\Lambda_u,\Psi)\in\Theta_1$, we have
\[
A(h_1,0)
=
\begin{bmatrix}
2\Lambda_u & \Gamma_\Psi\\
\Gamma_\Psi^\top & 0
\end{bmatrix}.
\]
Therefore,
\begin{align*}
    \|A(h_1,0)\|^2_F
&= 4 \| \Lambda_u\|^2_F + 2\| \Gamma_\Psi \|^2_F\\
&\leq 4\left(
\|\Lambda_u\|_F^2
+
\sum_{t=1}^{T-1}\|\Psi_t\|_F^2
\right)\\
&= 
4\|h_1\|_{\Theta_1}^2.
\end{align*}

 Taking square roots proves
 \eqref{eqn: A first block upper bound}.
\item Proof of \eqref{eqn: A second block upper bound}: Taking \(h_2\in\Theta_2\), we have
\[
A(0,h_2)
=
\begin{bmatrix}
0&0\\
0&2h_2
\end{bmatrix}.
\]
Therefore, 
\begin{align*}
\|A(0,h_2)\|_F
&= 2 \|h_2\|_F =
2\|h_2\|_{\Theta_2}.
\end{align*}
This proves \eqref{eqn: A second block upper bound}.
\item Proof of \eqref{eqn: A full space lower bound}:
Take $\theta=(\Lambda_u,\Lambda_y,\Psi)$, by the definition of the Frobenius norm and \eqref{eqn: Gamma Psi Frobenius norm},
\begin{align*}
&\|A(\theta)\|_F^2
=
4\|\Lambda_u\|_F^2
+
4\|\Lambda_y\|_F^2
+
2\|\Gamma_\Psi\|_F^2  = 4\|\Lambda_u\|_F^2
+
4\|\Lambda_y\|_F^2
+
2\sum_{t=1}^{T-1}\|\Psi_t\|_F^2 
\\&
\geq 2 \|\theta\|_\Theta^2
\end{align*}
Taking square roots proves \eqref{eqn: A full space lower bound}.
\end{enumerate}

\end{proof}

The following lemma collects two standard finite-dimensional consequences of Hessian bounds. The strong convexity statement is the usual implication from a uniform lower Hessian bound. For the block--coordinate Lipschitz statement, we apply the standard implication from an upper Hessian bound to Lipschitz continuity of the gradient after restricting the function to one coordinate block while keeping the other block fixed. We include the proof for completeness.
\begin{lemma}
\label{lemma: Hessian bounds imply block Lipschitz and strong convexity}
Let \(\mathcal X_1\) and \(\mathcal X_2\) be finite-dimensional real
inner-product spaces, and let
\(\mathcal X:=\mathcal X_1\times\mathcal X_2\)
be equipped with the product inner product. Let
\(\Omega\subset\mathcal X\) be open, let
\(\varphi\in C^2(\Omega)\), and let
\(K\subset\Omega\) be nonempty and convex.

For \(i=1,2\), define \begin{subequations}
\label{eqn: def of embedding selection filled with zero}
\begin{align}
J_1:\mathcal X_1&\to\mathcal X,
\qquad
J_1(h_1):=(h_1,0),\\
J_2:\mathcal X_2&\to\mathcal X,
\qquad
J_2(h_2):=(0,h_2).
\end{align}
\end{subequations}

Then:
\begin{enumerate}[label=\alph*.]

\item Suppose that, for \(i=1,2\), there exists \(L_i>0\) such that
\begin{equation}
0
\leq
\nabla^2\varphi(x)
\left[
J_i(h_i),J_i(h_i)
\right]
\leq
L_i\|h_i\|_{\mathcal X_i}^2
\label{eqn: general block Hessian upper bound}
\end{equation}
for every \(x\in K\) and \(h_i\in\mathcal X_i\).
Then the gradient of \(\varphi\) is block--coordinate Lipschitz continuous
on \(K\) with constants \(L_1,L_2\).
\item Suppose that, for some \(\sigma>0\),
\begin{equation}
\nabla^2\varphi(x)
[h,h]
\geq
\sigma\|h\|_{\mathcal X}^2
\label{eqn: general Hessian lower bound}
\end{equation}
for every \(x\in K\) and \(h\in\mathcal X\).
Then \(\varphi\) is strongly convex on \(K\) with constant \(\sigma\).
\end{enumerate}
\end{lemma}
\begin{proof}[Proof of Lemma~\ref{lemma: Hessian bounds imply block Lipschitz and strong convexity}]
    Let
\(x,\widetilde{x}\in K\), and consider the line segment
\begin{equation}
x(t)
:=
x+t(\widetilde{x}-x).
\label{eqn: line segment used for lip and str convex}
\end{equation}
We have
\(x(0)=x\), \(x(1)=\widetilde{x}\), and
$x(t)
\in K
\text{ for every }t\in[0,1],$
since \(K\) is convex. 

We first compute the derivative of \(\nabla\varphi(x(t))\) along this line
segment, which will be used in both parts of the proof. We have
\begin{equation}
\frac{d}{dt}x(t)
=
\widetilde{x}-x.
\label{eqn: time derivatvie along theta t line segement in lip and str cvx proof}
\end{equation}

And $\{t\in\mathbb R:x(t)\in\Omega\}$
is open, since \(x(\cdot):\mathbb R\to\mathcal X\) is continuous and
\(\Omega\) is open. Moreover, it contains \([0,1]\), since
\(x(t)\in K\subset\Omega\) for every
\(t\in[0,1]\).
Therefore, for every \(t\in[0,1]\), there exists an open interval
\(I_t\) containing \(t\) such that
\begin{equation}
x(r)\in\Omega
\qquad
\text{for every }r\in I_t.
\label{eqn: open set at any t in 0,1 st feasible}
\end{equation}
Applying
\eqref{eqn: time derivative of composition is D  in directionn of inner function time derivative}
on the open interval \(I_t\) with \(q=\nabla\varphi\) gives at \(t\),
\begin{equation}
\frac{d}{dt}\nabla\varphi(x(t))
=
D\nabla\varphi(x(t))
\left[\widetilde{x}-x\right].
\label{eqn:time derivative of nabla f theta t}
\end{equation}
Since \(t\in[0,1]\) was arbitrary, this holds for every \(t\in[0,1]\).
\begin{enumerate} [label=\alph*.]
\item
To prove the block--coordinate Lipschitz continuity, fix \(i=1,2\),
and take \(x,\widetilde{x}\in K\) that differ only in the \(i\)-th
block. If \(i=1\), write
\begin{equation*}
x=(x_1,x_2),
\qquad
\widetilde{x}=(\widetilde{x}_1,x_2),
\end{equation*}
and if \(i=2\), write
\begin{equation*}
x=(x_1,x_2),
\qquad
\widetilde{x}=(x_1,\widetilde{x}_2).
\end{equation*}
In either case, define
\begin{equation*}
h_i:=\widetilde{x}_i-x_i\in\mathcal X_i.
\end{equation*}
Then
\begin{equation}
\widetilde{x}-x
=
J_i(h_i).
\label{eqn: in block lip proof theta perturb is embedded hi}
\end{equation}

For \(x=(x_1,x_2)\in\Omega\), we write
\begin{equation*}
\nabla\varphi(x)
=
\left(
\nabla_{x_1}\varphi(x),
\nabla_{x_2}\varphi(x)
\right),
\end{equation*}
where
\(\nabla_{x_i}\varphi(x)\in\mathcal X_i\) denotes the gradient with
respect to the \(i\)-th block.

We now bound the change of the \(i\)-th block gradient along the line
segment. Denote
\begin{equation*}
q(t):=\nabla_{x_i}\varphi(x(t))\in\mathcal X_i.
\end{equation*}
We first show that
\begin{equation}
\|q'(t)\|_{\mathcal X_i}
\leq
L_i\|h_i\|_{\mathcal X_i},
\qquad
t\in[0,1].
\label{eqn: q prime bounded by Li hi}
\end{equation}
Fix \(t\in[0,1]\) and take arbitrary \(g_i\in\mathcal X_i\).
We have
\begin{equation}
\left\langle
g_i,q(t)
\right\rangle_{\mathcal X_i}
=
\left\langle
J_i(g_i),
\nabla\varphi(x(t))
\right\rangle_{\mathcal X}.
\label{eqn: inner of gi an qt is inner of Jgi and nablda f theta t}
\end{equation}
To see why, we take \(i=1\) for instance,
\begin{align*}
\left\langle
J_1(g_1),
\nabla\varphi(x(t))
\right\rangle_{\mathcal X}
&=
\left\langle
g_1,
\nabla_{x_1}\varphi(x(t))
\right\rangle_{\mathcal X_1},
\end{align*}
since
\(J_1(g_1)=(g_1,0)\) and
\(
\nabla\varphi(x(t))
=
\left(
\nabla_{x_1}\varphi(x(t)),
\nabla_{x_2}\varphi(x(t))
\right).
\)
The same argument gives the corresponding identity for \(i=2\).

Since \(g_i\) is fixed, differentiating
\eqref{eqn: inner of gi an qt is inner of Jgi and nablda f theta t}
with respect to \(t\) gives
\begin{align*}
\left\langle
g_i,q'(t)
\right\rangle_{\mathcal X_i}
&=
\left\langle
J_i(g_i),
\frac{d}{dt}\nabla\varphi(x(t))
\right\rangle_{\mathcal X}\\
&=
\left\langle
J_i(g_i),
D\nabla\varphi(x(t))
\left[J_i(h_i)\right]
\right\rangle_{\mathcal X}\\
&=
\nabla^2\varphi(x(t))
\left[
J_i(g_i),
J_i(h_i)
\right].
\end{align*}
Here, the second equality follows from
\eqref{eqn:time derivative of nabla f theta t} and
\eqref{eqn: in block lip proof theta perturb is embedded hi}, and the
third equality follows from
\eqref{eqn: Hessian bilinear form via gradient derivative}.

For fixed \(x(t)\in K\), the map
\begin{equation*}
(g_i,\widetilde g_i)
\longmapsto
\nabla^2\varphi(x(t))
\left[
J_i(g_i),J_i(\widetilde g_i)
\right]
\end{equation*}
is a symmetric positive semidefinite bilinear form.
It is symmetric because
\(\varphi\in C^2(\Omega)\), so
\(\nabla^2\varphi(x(t))=D^2\varphi(x(t))\) is symmetric.
It is positive semidefinite because
\eqref{eqn: general block Hessian upper bound} gives
\begin{equation*}
\nabla^2\varphi(x(t))
\left[
J_i(g_i),J_i(g_i)
\right]
\geq0,
\qquad
\forall g_i\in\mathcal X_i.
\end{equation*}
Hence, by the Cauchy--Schwarz inequality for positive semidefinite
bilinear forms,
\begin{align*}
\numberthis\label{eqn: inner produce of gi and q prime is bouned by Li gi norma nd hi norm}
\left|
\left\langle
g_i,q'(t)
\right\rangle_{\mathcal X_i}
\right|
&=
\left|
\nabla^2\varphi(x(t))
\left[
J_i(g_i),J_i(h_i)
\right]
\right|\\
&\leq
\sqrt{
\nabla^2\varphi(x(t))
\left[
J_i(g_i),J_i(g_i)
\right]
}
\sqrt{
\nabla^2\varphi(x(t))
\left[
J_i(h_i),J_i(h_i)
\right]
}\\
&\leq
L_i
\|g_i\|_{\mathcal X_i}
\|h_i\|_{\mathcal X_i}, 
\end{align*}
where the last inequality follows from
\eqref{eqn: general block Hessian upper bound}.
Since \(g_i\in\mathcal X_i\) was arbitrary, the preceding inequality
holds for every \(g_i\in\mathcal X_i\) at the fixed \(t\).
If \(q'(t)=0\), then
\begin{equation*}
\|q'(t)\|_{\mathcal X_i}
\leq
L_i\|h_i\|_{\mathcal X_i}
\end{equation*}
holds trivially. Otherwise, take
\begin{equation*}
g_i
=
\frac{q'(t)}
{\|q'(t)\|_{\mathcal X_i}}.
\end{equation*}
Then
\begin{align*}
\|q'(t)\|_{\mathcal X_i}
&=
\left|
\left\langle
g_i,q'(t)
\right\rangle_{\mathcal X_i}
\right|\\
&\leq
L_i\|h_i\|_{\mathcal X_i},
\end{align*}
where the inequality uses \eqref{eqn: inner produce of gi and q prime is bouned by Li gi norma nd hi norm} and  \(\|g_i\|_{\mathcal X_i}=1\). This proves \eqref{eqn: q prime bounded by Li hi}.

Finally, integrating this bound along the line segment gives us the bound
on the block-gradient difference. The function \(q(\cdot)\) is continuously
differentiable on the open set
$\{t\in\mathbb R:x(t)\in\Omega\},$ which contains \([0,1]\), because 
\(\varphi\in C^2(\Omega)\). Therefore, by the Fundamental Theorem of
Calculus,
\begin{align*}
\left\|
\nabla_{x_i}\varphi(\widetilde{x})
-
\nabla_{x_i}\varphi(x)
\right\|_{\mathcal X_i}
&=
\|q(1)-q(0)\|_{\mathcal X_i}\\
&=
\left\|
\int_0^1q'(t)\,dt
\right\|_{\mathcal X_i}\\
&\leq
\int_0^1\|q'(t)\|_{\mathcal X_i}\,dt\\
&\leq
L_i\|h_i\|_{\mathcal X_i}.
\end{align*}
Thus, since \(i=1,2\) was arbitrary, the gradient of \(\varphi\) is
block--coordinate Lipschitz continuous on \(K\) with constants
\(L_1,L_2\).

\item 
We next use
\eqref{eqn: general Hessian lower bound} to prove the
strong convexity.
Take arbitrary
\(x,\widetilde{x}\in K\), and let \(x(t)\) be the line
segment defined in \eqref{eqn: line segment used for lip and str convex}. Define
\begin{equation*}
\phi(t):=\varphi(x(t)).
\end{equation*}
We first compute the first and second derivatives of \(\phi\) along the line
segment.
Fix an arbitrary \(t\in[0,1]\). Recall from
\eqref{eqn: open set at any t in 0,1 st feasible} that there exists an open interval
\(I_t\) containing \(t\) such that
\(x(r)\in\Omega\) for every \(r\in I_t\).
Applying \eqref{eqn: time derivative of composition is D  in directionn of inner function time derivative}
on the open interval \(I_t\) with \(q=\varphi\) gives, in particular at \(t\),
\begin{align*}
\numberthis
\label{eqn: phi time derivative is inner produce of nabla f theta t and diff of thet and thetat tilde}
\phi'(t)
&=
D\varphi(x(t))
\left[
\frac{d}{dt}x(t)
\right]\\
&=
D\varphi(x(t))
\left[
\widetilde{x}-x
\right]
\text{ (by \eqref{eqn: time derivatvie along theta t line segement in lip and str cvx proof})}\\
&=
\left\langle
\nabla \varphi(x(t)),
\widetilde{x}-x
\right\rangle_{\mathcal X},
\end{align*}
where the last equality is by the definition of the gradient.\\
Since \(\widetilde{x}-x\) is fixed, differentiating the preceding
identity with respect to \(t\) gives
\begin{align*}
\phi''(t)
&=
\left\langle
\frac{d}{dt}\nabla \varphi(x(t)),
\widetilde{x}-x
\right\rangle_{\mathcal X}\\
&=
\left\langle
D\nabla \varphi(x(t))
\left[
\widetilde{x}-x
\right],
\widetilde{x}-x
\right\rangle_{\mathcal X}\\
&=
\nabla^2\varphi(x(t))
[\widetilde{x}-x,\widetilde{x}-x],
\end{align*}
where the second equality uses
\eqref{eqn:time derivative of nabla f theta t}, and the last equality
follows from
\eqref{eqn: Hessian bilinear form via gradient derivative}.

Since
\(x(t)\in K\) for every \(t\in[0,1]\),
\eqref{eqn: general Hessian lower bound} gives
\begin{align*}
\numberthis
\label{eqn: phi prime prime is bounded by strong convex and norm sqruare of theta diff}
\phi''(t)
&=
\nabla^2\varphi(x(t))
[\widetilde{x}-x,\widetilde{x}-x]\\
&\geq
\sigma
\|\widetilde{x}-x\|_{\mathcal X}^2.
\end{align*}
Since \(\varphi\in C^2(\Omega)\) and \(x(\cdot)\) is affine,
the function \(\phi(\cdot)=\varphi(x(\cdot))\) is \(C^2\) on an open set containing \([0,1]\). Using the Fundamental Theorem of Calculus,
\begin{align*}
\phi(1)\numberthis
\label{eqn: twice ftc to give phi(1)-phi(0)}
&=\phi(0)+
\int_0^1 \phi'(s)\,ds\\
&=\phi(0)+
\int_0^1
\left(
\phi'(0)+\int_0^s\phi''(t)\,dt
\right)ds\\
&=\phi(0)+
\phi'(0)
+
\int_0^1\int_0^s\phi''(t)\,dt\,ds\\
&=\phi(0)+
\phi'(0)
+
\int_0^1\int_t^1\phi''(t)\,ds\,dt\\
&=\phi(0)+
\phi'(0)
+
\int_0^1(1-t)\phi''(t)\,dt.
\end{align*}
Here, the second equality follows again from the Fundamental Theorem of
Calculus,
and the fourth equality follows by changing the order of integration.
Moreover, by
\eqref{eqn: phi time derivative is inner produce of nabla f theta t and diff of thet and thetat tilde},
\begin{equation*}
\phi'(0)
=
\left\langle
\nabla \varphi(x(0)),
\widetilde{x}-x
\right\rangle_{\mathcal X}
=
\left\langle
\nabla \varphi(x),
\widetilde{x}-x
\right\rangle_{\mathcal X}.
\end{equation*}
Recalling
\(\varphi(\widetilde{x})=\phi(1)\),
\(\varphi(x)=\phi(0)\),
\eqref{eqn: twice ftc to give phi(1)-phi(0)} becomes
\begin{align*}
\varphi(\widetilde{x})
&=
\varphi(x)
+
\left\langle
\nabla \varphi(x),
\widetilde{x}-x
\right\rangle_{\mathcal X}
+
\int_0^1
(1-t)\phi''(t)\,dt\\
&\geq
\varphi(x)
+
\left\langle
\nabla \varphi(x),
\widetilde{x}-x
\right\rangle_{\mathcal X}
+
\sigma
\|\widetilde{x}-x\|_{\mathcal X}^2
\int_0^1(1-t)\,dt\\
&=
\varphi(x)
+
\left\langle
\nabla \varphi(x),
\widetilde{x}-x
\right\rangle_{\mathcal X}
+
\frac{\sigma}{2}
\|\widetilde{x}-x\|_{\mathcal X}^2,
\end{align*}
where the inequality uses
\eqref{eqn: phi prime prime is bounded by strong convex and norm sqruare of theta diff}.
Thus, \(\varphi\) is strongly convex on \(K\) with constant
\(\sigma\).
\end{enumerate}
\end{proof}
}

\begin{proof}[Proof of Lemma~\ref{lemma: local lip and str convex}]
   For \(i=1,2\), let \(J_i:\Theta_i \to\Theta\) be the embedding such that $J_1(\theta_1) = (\theta_1,0)$ and $J_2(\theta_2) = (0,\theta_2)$. 
\begin{enumerate}[label=\alph*.]
    \item 
We establish a uniform upper bound for the Hessian evaluated along perturbations in the \(i\)-th block, for $i=1,2$.
For every \(\vartheta\in K\) and \(g_i\in\Theta_i\),
by Lemma \ref{lemma: hessian of f},
\begin{align*}
\nabla^2 f(\vartheta)
\left[
J_i(g_i),
J_i(g_i)
\right]
&=
\frac12
\left\|
Z(\vartheta)^{-1/2}
A\bigl(J_i(g_i)\bigr)
Z(\vartheta)^{-1/2}
\right\|_F^2.
\end{align*}
Using the inequality
\(
\|AB\|_F\leq \|A\|_2\|B\|_F
\)
twice, we upper-bound  the right--hand side according to
\begin{align*}
\left\|
Z(\vartheta)^{-1/2}
A\bigl(J_i(g_i)\bigr)
Z(\vartheta)^{-1/2}
\right\|_F
&\leq
\|Z(\vartheta)^{-1/2}\|_2
\left\|
A\bigl(J_i(g_i)\bigr)
\right\|_F
\|Z(\vartheta)^{-1/2}\|_2.
\numberthis
\label{eqn: ZAZ on gi is bounded by Z-12 A Z-12}
\end{align*}

We have
\[
\|Z(\vartheta)^{-1/2}\|_2 = \lambda_{\max}(Z(\vartheta)^{-1/2}) = (\lambda_{\min}(Z(\vartheta)))^{-\frac{1}{2}}
\leq
m_K^{-\frac12}.
\]
where we used the fact that \(Z(\vartheta)\succ0\). 
Thus,
\eqref{eqn: ZAZ on gi is bounded by Z-12 A Z-12} gives
\begin{align*}
\left\|
Z(\vartheta)^{-1/2}
A\bigl(J_i(g_i)\bigr)
Z(\vartheta)^{-1/2}
\right\|_F
&\leq
m_K^{-1}
\left\|
A\bigl(J_i(g_i)\bigr)
\right\|_F.
\end{align*}
We further have
\begin{align*}
\left\|
A\bigl(J_i(g_i)\bigr)
\right\|_F
&\leq
2
\|g_i\|_{\Theta_i},
\end{align*}
due to \eqref{eqn: A first block upper bound} when \(i=1\), and
\eqref{eqn: A second block upper bound} when \(i=2\).
Consequently,
\begin{equation}
\left\|
Z(\vartheta)^{-1/2}
A\bigl(J_i(g_i)\bigr)
Z(\vartheta)^{-1/2}
\right\|_F
\leq
2\,m_K^{-1}
\|g_i\|_{\Theta_i}.
\label{eqn: ZAZ norm for gi }
\end{equation}
concluding the following bound on the Hessian: 
\begin{align*}
\nabla^2 f(\vartheta)
\left[
J_i(g_i),
J_i(g_i)
\right]\numberthis
\label{eqn: upper bound on block Hessian quadratic form}
&=
\frac12
\left\|
Z(\vartheta)^{-1/2}
A\bigl(J_i(g_i)\bigr)
Z(\vartheta)^{-1/2}
\right\|_F^2\\
&\leq
2\,m_K^{-2}
\|g_i\|_{\Theta_i}^2\\
&=
L_i\|g_i\|_{\Theta_i}^2.
\end{align*}
Meanwhile,
\begin{equation*}
\nabla^2 f(\vartheta)
\left[
J_i(g_i),J_i(g_i)
\right]
=
\frac12
\left\|
Z(\vartheta)^{-1/2}
A\bigl(J_i(g_i)\bigr)
Z(\vartheta)^{-1/2}
\right\|_F^2
\geq0.
\end{equation*}
Thus, applying the standard Hessian criterion for Lipschitz continuity of the gradient to each coordinate restriction, as recorded in
Lemma~\ref{lemma: Hessian bounds imply block Lipschitz and strong convexity},
shows that the gradient of \(f\) is block--coordinate Lipschitz continuous
on \(K\) with constants \(L_1,L_2\).


\item    
We establish a uniform positive lower bound for the Hessian evaluated along perturbations in \(\Theta\). 
For every \(\vartheta\in K\), and for any direction
\(\delta\theta\in\Theta\), by Lemma \ref{lemma: hessian of f}, 
\begin{equation}
    \nabla^2 f(\vartheta)[\delta\theta,\delta\theta]
=
\frac12
\left\|
Z(\vartheta)^{-1/2}A(\delta\theta)Z(\vartheta)^{-1/2}
\right\|_F^2, \label{eqn: hessian formula at var phi used in strong convex}
\end{equation}
To obtain a lower bound on this quantity, we first relate it to
\(\|A(\delta\theta)\|_F\).
\begin{align*}
\|A(\delta \theta)\|_F\numberthis 
\label{eqn: A delta theta F is bounded by Z var theta 2  ZAZ F Z var theta 2}
&=
\|Z(\vartheta)^{1/2}
(Z(\vartheta)^{-1/2} A(\delta \theta) Z(\vartheta)^{-1/2})
Z(\vartheta)^{1/2}\|_F\\
&\leq
\|Z(\vartheta)^{1/2}\|_2
\|Z(\vartheta)^{-1/2} A(\delta \theta) Z(\vartheta)^{-1/2}\|_F
\|Z(\vartheta)^{1/2}\|_2 \\ &\leq M_K
\|Z(\vartheta)^{-1/2} A(\delta \theta) Z(\vartheta)^{-1/2}\|_F, 
\end{align*}
where the second inequality follows from
\[
\|Z(\vartheta)^{1/2}\|_2 = \lambda_{\max}(Z(\vartheta))^{\frac{1}{2}}
\leq
M_K^{\frac12}.
\]
Therefore \eqref{eqn: A delta theta F is bounded by Z var theta 2  ZAZ F Z var theta 2} concludes the lower-bound
\[
\left\|
Z(\vartheta)^{-1/2}
A(\delta\theta)
Z(\vartheta)^{-1/2}
\right\|_F \geq \frac{\|A(\delta\theta)\|_F}{M_K},
\]
which is used in 
\eqref{eqn: hessian formula at var phi used in strong convex} to conlcude the lower-bound on the Hessian:
\begin{align*}
\numberthis
\label{eqn: min eigvalues of the hessian in strong convex}\nabla^2 f(\vartheta)[\delta\theta,\delta\theta]
&\geq
\frac{1}{2M_K^2}\|A(\delta\theta)\|_F^2\\
&\geq
\frac{1}{M_K^2}\|\delta\theta\|_{\Theta}^2\\ &=
\sigma_K\|\delta\theta\|_\Theta^2,
\end{align*}
where the second inequality follows from
\eqref{eqn: A full space lower bound}. 
Thus, by the standard implication
from a uniform lower Hessian bound to strong convexity, as recorded in
Lemma~\ref{lemma: Hessian bounds imply block Lipschitz and strong convexity},
\(f\) is strongly convex on \(K\) with constant \(\sigma_K\).
\end{enumerate}

\end{proof}
\subsection{Supplementary Analysis for Section  \ref{section: sublevel set restriction}}\label{section: supp sublevel set restriction}
\begin{proof}[Proof of Lemma \ref{lemma: uniform spectral bound on level set}]
First, we establish that $-e^{-C_Q}$ belongs to the domain of $W_0$ and $W_{-1}$, i.e. $-e^{-C_Q} \in (-e^{-1},0)$, by showing that 
\(C_Q>1\). To see this, let
\(q_1,\ldots,q_{2Td}\) be the eigenvalues of \(Q\).
We have \(q_i>1\) for every \(i=1,\ldots,2Td\), because $Q=(\Sigma^\gamma)^{-1}+I_{2Td}$ by \eqref{eqn:primal loss} and \((\Sigma^\gamma)^{-1}\succ0\) by Assumption \ref{assumption:cost}. Since the trace of a square matrix is equal to the sum of its eigenvalues
and its determinant is equal to the product of its eigenvalues, we have
\begin{equation*}
\operatorname{Tr}(Q)-\log\det Q
=
\sum_{i=1}^{2Td}
\left(q_i-\log q_i\right)
>
2Td,
\end{equation*}
where we used \(x-\log x>1\) for \(x>1\). This proves \(C_Q>1\).

Next, we establish a condition for all $\theta \in L(\theta^{\{0\}}) $ when $\theta^{\{0\}}=0$.
Express the objective function as
\begin{align*}
   f(\theta)
&=
\frac12\operatorname{Tr}(Z(\theta))
-
\frac12\log\det Z(\theta)
-
\frac12\operatorname{Tr}(Q)
-Td 
\end{align*}
where we used $\trace(Z(\theta)) = \trace(Q) + 2 \trace(\Lambda_u) + 2 \trace(\Lambda_y)$ according to~\eqref{eqn: Z from theta}. Simplifying further by expressing the trace and determinant in terms of eigenvalues yields 
 \begin{align*}
f(\theta)
&=
\frac12
\sum_{i=1}^n
\left(
\lambda_i(\theta)-\log\lambda_i(\theta)
\right)
-
\frac12\operatorname{Tr}(Q)
-Td \\
&= \frac12
\sum_{i=1}^n
h(\lambda_i(\theta))
-
\frac12\operatorname{Tr}(Q)
-Td, 
\end{align*}
where 
\(\lambda_1(\theta),\ldots,\lambda_n(\theta)>0\)
are the eigenvalues of \(Z(\theta)\succ0\), $n=2Td$, and we introduced the function 
\begin{equation*}
h(\lambda):=\lambda-\log\lambda,
\qquad
\lambda>0.
\end{equation*}
According to this expression, the sub-level condition $f(\theta) \leq f(\theta^{\{0\}})=f(0)$ is equivalent to
\begin{equation}
\sum_{i=1}^n h(\lambda_i(\theta))
\leq
\operatorname{Tr}(Q)-\log\det Q.
\label{eqn: sum h eigenvalues bounded}
\end{equation}
because $Z(0)=Q$. Moreover, $h(\lambda) \geq 1$ for $\lambda>0$ due to the inequality $\log(\lambda) \leq \lambda -1$. Therefore, for every \(i=1,\ldots,n\),
\begin{align*}
h(\lambda_i(\theta))
&=
\sum_{j=1}^n h(\lambda_j(\theta))
-
\sum_{j\neq i}h(\lambda_j(\theta))\\
&\leq
\operatorname{Tr}(Q)-\log\det Q-(n-1)\\
&=
C_Q,
\end{align*}
concluding the sub-level condition
\begin{equation}
h(\lambda_i(\theta))
\leq C_Q,
\qquad
i=1,\ldots,n,
\qquad
\forall\theta\in L(\theta^{\{0\}}).
\label{eqn: h eigenvalue bounded by CQ}
\end{equation}
The graph of the function $h(\lambda)=\lambda - \log(\lambda)$ attains its minimum $1$ at $\lambda=1$, monotonically increases to $+\infty$ as $\lambda \downarrow 0$  and $\lambda \to \infty$, so it crosses $C_Q>1$ at exactly two points, given by $a_Q=-W_0(-e^{-C_Q})<1,
$ and $
b_Q=-W_{-1}(-e^{-C_Q})>1.$  Therefore, the sub-level condition implies 
\begin{align*}
    \lambda_i(\theta) \in [a_Q,b_Q],\qquad
i=1,\ldots,n,
\qquad
\forall\theta\in L(\theta^{\{0\}}),
\end{align*}
concluding the proof. 

\end{proof}

The next lemma shows that the block subproblems \eqref{eq:update-lambday} and 
\eqref{eq:update-upsi} in the BCD scheme decrease the objective and, consequently, that all BCD iterates remain in the initial sublevel set \(L(\theta^{\{0\}})\).
\begin{lemma} \label{lemma: monotone}
Let $\{\theta^{\{k\}}, \theta^{\{k+\frac{1}{2}\}}\}_{k \geq 0}$ be the sequence of iterates generated by Algorithm \ref{alg:bcd-dual}, initialized with $\theta^{\{0\}} = \mathbf{0}$.\\
For every \(k\geq0\),
 \begin{equation*}f(\theta^{\{k+1\}}) \le f(\theta^{\{k+\frac{1}{2}\}}) \le f(\theta^{\{k\}}).\end{equation*}
 Moreover, all iterates remain in the level set ${L}(\theta^{\{0\}})$:
\begin{equation} \label{eqn: all iterations are in the level set}
   \forall k, \quad \theta^{\{k\}},\theta^{\{k+\frac{1}{2}\}}  \in L(\theta^{\{0\}}) 
\end{equation}
In particular, since \(L(\theta^{\{0\}})\subset\Thetafeasible\), all iterates remain in \(\Thetafeasible\).
\end{lemma}
\begin{proof}[Proof of Lemma~\ref{lemma: monotone}]
The optimality of the block subproblems \eqref{eq:update-upsi} and \eqref{eq:update-lambday} imply 
\begin{equation*}
f(\theta^{\{k+1\}})
\leq
f(\theta^{\{k+\frac12\}})
\leq
f(\theta^{\{k\}}).
\end{equation*}
The proof follows by induction starting from $k=0$ and the definition of the sub-level set $L(\theta^{\{0\}})$. 




\end{proof}

\begin{proof}[Proof of Corollary~\ref{corollary: lip and strong convexity on sublevel set}]
The proof follows by applying Lemma~\ref{lemma: local lip and str convex}
with $K=L(\theta^{\{0\}})$ and using
Lemma~\ref{lemma: uniform spectral bound on level set}
to bound $m_K$ and $M_K$.

\end{proof}
\subsection{Supplementary Analysis for Section~\ref{subsection: linear convergence final rate}}\label{subsec: linear convergence final rate}
As noted in the main text, our proof adapts the linear convergence argument of  \cite[Theorem~5.2]{BeckAmirConvergenceofBlockCoordinateDescent}, which is not directly applicable to our setting. In particular, the strong
convexity of \(f\) and the block--coordinate Lipschitz continuity of its
gradient can only be established on the sublevel set, rather than globally. Therefore, to carry out the standard descent argument on this subset, we need to show that the points obtained by taking the block gradient steps used in the descent estimates remain in the sublevel set (see Lemma \ref{lemma: block descent estimate on the descent direction}). This is the main additional step needed to adapt the standard linear convergence argument
to our setting.
We first recall a standard consequence from the Lipschitz continuity of the gradient.
\begin{lemma}
\label{lemma: block descent lemma on sublevel set}
Let
\(
\theta=(\theta_1,\theta_2)\in L(\theta^{\{0\}})
\).
Then:
\begin{enumerate}[label=\alph*.]
    \item If
    \(
    (\bar\theta_1,\theta_2)\in L(\theta^{\{0\}})
    \),
    then
    \begin{equation}
    \begin{aligned}
    f(\bar\theta_1,\theta_2)
    &\leq
    f(\theta_1,\theta_2)
    +
    \left\langle
    \nabla_{\theta_1}f(\theta_1,\theta_2),
    \bar\theta_1-\theta_1
    \right\rangle_{\Theta_1}
    +
    \frac{L^Q}{2}
    \|\bar\theta_1-\theta_1\|_{\Theta_1}^2.
    \end{aligned}
    \label{eqn: block descent lemma theta1}
    \end{equation}

    \item If
    \(
    (\theta_1,\bar\theta_2)\in L(\theta^{\{0\}})
    \),
    then
    \begin{equation}
    \begin{aligned}
    f(\theta_1,\bar\theta_2)
    &\leq
    f(\theta_1,\theta_2)
    +
    \left\langle
    \nabla_{\theta_2}f(\theta_1,\theta_2),
    \bar\theta_2-\theta_2
    \right\rangle_{\Theta_2}
    +
    \frac{L^Q}{2}
    \|\bar\theta_2-\theta_2\|_{\Theta_2}^2.
    \end{aligned}
    \label{eqn: block descent lemma theta2}
    \end{equation}
\end{enumerate}
\end{lemma}
\begin{proof}[Proof of Lemma~\ref{lemma: block descent lemma on sublevel set}]
For each part, the line segment joining
\((\theta_1,\theta_2)\) to
\((\bar\theta_1,\theta_2)\), or joining
\((\theta_1,\theta_2)\) to
\((\theta_1,\bar\theta_2)\),
is contained in \(L(\theta^{\{0\}})\), since
\(L(\theta^{\{0\}})\) is convex.
Both inequalities follow from the standard descent lemma applied to the
corresponding one-block restriction of \(f\), using the block--coordinate
Lipschitz continuity from
Corollary~\ref{corollary: lip and strong convexity on sublevel set} and
the convexity of \(L(\theta^{\{0\}})\).
\end{proof}

\begin{lemma}
\label{lemma: block descent estimate on the descent direction}
Let
\(\theta=(\theta_1,\theta_2)\in L(\theta^{\{0\}})\).

\begin{enumerate}[label=\alph*.]
    \item Define
    \begin{equation*}
    \widetilde\theta_1
    :=
    \theta_1
    -
    \frac{1}{L^Q}
    \nabla_{\theta_1}f(\theta_1,\theta_2).
    \end{equation*}
    Then
    \begin{equation*}
    (\widetilde\theta_1,\theta_2)
    \in
    L(\theta^{\{0\}}),
    \end{equation*}
    and
    \begin{equation}
    f(\theta_1,\theta_2)
    -
    f(\widetilde\theta_1,\theta_2)
    \geq
    \frac{1}{2L^Q}
    \left\|
    \nabla_{\theta_1}f(\theta_1,\theta_2)
    \right\|_{\Theta_1}^2.
    \label{eqn: block descent estimate theta1}
    \end{equation}

    \item Define
    \begin{equation*}
    \widetilde\theta_2
    :=
    \theta_2
    -
    \frac{1}{L^Q}
    \nabla_{\theta_2}f(\theta_1,\theta_2).
    \end{equation*}
    Then
    \begin{equation*}
    (\theta_1,\widetilde\theta_2)
    \in
    L(\theta^{\{0\}}),
    \end{equation*}
    and
    \begin{equation}
    f(\theta_1,\theta_2)
    -
    f(\theta_1,\widetilde\theta_2)
    \geq
    \frac{1}{2L^Q}
    \left\|
    \nabla_{\theta_2}f(\theta_1,\theta_2)
    \right\|_{\Theta_2}^2.
    \label{eqn: block descent estimate theta2}
    \end{equation}
\end{enumerate}
\end{lemma}
\begin{proof}[Proof of Lemma~\ref{lemma: block descent estimate on the descent direction}]
We prove part (a) first.
Let
\begin{equation*}
h_1
:=
-\frac{1}{L^Q}
\nabla_{\theta_1}f(\theta_1,\theta_2),
\end{equation*}
so that
\(
\widetilde\theta_1=\theta_1+h_1
\).
Define the affine path
\begin{equation*}
\theta(t):=(\theta_1+th_1,\theta_2),
\qquad
t\in\mathbb R.
\end{equation*}

We first show that the entire segment remains in the sublevel set, i.e.,
\begin{equation*}
\theta(t)\in L(\theta^{\{0\}})
\qquad
\text{for every }t\in[0,1].
\end{equation*}

If
\(
\nabla_{\theta_1}f(\theta_1,\theta_2)=0
\),
then \(h_1=0\), so
\(
\theta(t)=\theta\in L(\theta^{\{0\}})
\)
for every \(t\in[0,1]\). It remains to consider the case
\begin{equation}
\nabla_{\theta_1}f(\theta_1,\theta_2)\neq0.
\label{eqn: assumtion that block 1 derivaitve of f is not zero in the proof of stay in level set}
\end{equation}

Consider the portion of the segment that remains in the sublevel set.
\begin{equation*}
\mathcal I_L
:=
\{t\in[0,1]:
\theta(t)\in L(\theta^{\{0\}})\}.
\end{equation*}
We have
\begin{equation*}
\mathcal I_L=[0,\tau]
\end{equation*}
for some \(\tau\in[0,1]\). To see this, first, $0\in\mathcal I_L,$
since
\(\theta(0)=\theta\in L(\theta^{\{0\}})\). Second, the set \(\mathcal I_L\) is closed in \([0,1]\), since
\(L(\theta^{\{0\}})\) is closed by
Lemma~\ref{lemma: compact level set} and
\(\theta(\cdot)\) is continuous.
Third, the set \(\mathcal I_L\) is an interval in \(\mathbb R\), since it
is convex by the convexity of \(L(\theta^{\{0\}})\) and the affineness of
\(\theta(\cdot)\).

We next show that
\begin{equation*}
\tau>0.
\end{equation*}
We will use the derivative of \(f\) along the feasible portion of the path. Define
\begin{equation*}
\mathcal I_{\mathrm{feasible}}
:=
\{r\in\mathbb R:
\theta(r)\in\Thetafeasible\}.
\end{equation*}
The set \(\mathcal I_{\mathrm{feasible}}\) is open, since
\(\theta(\cdot):\mathbb R\to\Theta\) is continuous and
\(\Thetafeasible\) is open. Moreover, \(\theta(\cdot)\) is
Fr\'echet differentiable on \(\mathbb R\), with
\(
\frac{d}{dr}\theta(r)=(h_1,0)
\),
since \(\theta(\cdot)\) is affine. Since \(f\) is Fr\'echet
differentiable on \(\Thetafeasible\),
\eqref{eqn: time derivative of composition is D  in directionn of inner function time derivative} applied to
the restriction of \(\theta(\cdot)\) to
\(\mathcal I_{\mathrm{feasible}}\), gives
$\forall r\in\mathcal I_{\text{feasible}} $,
\begin{align}
\frac{d}{dr}f(\theta(r))
&=
Df(\theta(r))[(h_1,0)] \notag\\
&=
\left\langle
\nabla f(\theta(r)),
(h_1,0)
\right\rangle_\Theta \notag\\
&=
\left\langle
(\nabla_{\theta_1}f(\theta_1+rh_1,\theta_2),\nabla_{\theta_2}f(\theta_1+rh_1,\theta_2) ),
(h_1,0)
\right\rangle_{\Theta}\notag\\
&=
\left\langle
\nabla_{\theta_1}f(\theta_1+rh_1,\theta_2),
h_1
\right\rangle_{\Theta_1}.
\label{eqn: derivative of f along block one path}
\end{align}
We have
\(0\in\mathcal I_{\mathrm{feasible}}\), since
\(\theta(0)=\theta\in L(\theta^{\{0\}})\subset\Thetafeasible\).
Therefore, evaluating
\eqref{eqn: derivative of f along block one path}
at \(r=0\) gives
\begin{align*}
\left.
\frac{d}{dr}f(\theta(r))
\right|_{r=0}
&=
\left\langle
\nabla_{\theta_1}f(\theta_1,\theta_2),
h_1
\right\rangle_{\Theta_1}\\
&=
-\frac{1}{L^Q}
\left\|
\nabla_{\theta_1}f(\theta_1,\theta_2)
\right\|_{\Theta_1}^2\\
&<0,
\end{align*}
where the last inequality follows from
\eqref{eqn: assumtion that block 1 derivaitve of f is not zero in the proof of stay in level set}.
Since \(\mathcal I_{\mathrm{feasible}}\) is open and contains \(0\),
there exists \(\epsilon>0\) such that
\begin{equation*}
(-\epsilon,\epsilon)\subset\mathcal I_{\mathrm{feasible}}.
\end{equation*}
Since the derivative at \(r=0\) is strictly negative, there exists
\(0<\delta\leq\min\{\epsilon,1\}\) such that
\(
f(\theta(t))<f(\theta(0))\leq f(\theta^{\{0\}})
\)
for every \(t\in(0,\delta]\).
Hence,
\begin{equation*}
\theta(t)\in L(\theta^{\{0\}})
\qquad
\text{for every }t\in[0,\delta].
\end{equation*}
Thus,
\begin{equation*}
\tau\geq\delta>0.
\end{equation*}

We now show that
\begin{equation*}
\tau=1.
\end{equation*}
Suppose, to the contrary, that \(\tau<1\). 
Since \(\theta(\tau )\in  L(\theta^{\{0\}})\), applying Lemma \ref{lemma: block descent lemma on sublevel set} to
\(\theta\) and
\(\theta(\tau )\), we have
\begin{align}
f(\theta(\tau))-f(\theta)
&\leq
\tau
\left\langle
\nabla_{\theta_1}f(\theta_1,\theta_2),
h_1
\right\rangle_{\Theta_1}
+
\frac{L^Q\tau^2}{2}
\|h_1\|_{\Theta_1}^2
\notag\\
&=
-\frac{\tau-\frac12\tau^2}{L^Q}
\left\|
\nabla_{\theta_1}f(\theta_1,\theta_2)
\right\|_{\Theta_1}^2.
\label{eqn: descent estimate from gradient block lip until t}
\end{align}
Hence,
\begin{equation}
f(\theta(\tau))
<
f(\theta)
\leq
f(\theta^{\{0\}}),
\label{eqn: f theta tau is strictly bounded by f theta initial}
\end{equation}
since \(\tau\in(0,1)\) and
\eqref{eqn: assumtion that block 1 derivaitve of f is not zero in the proof of stay in level set}
holds.
We next show that the path remains in the sublevel set for a small time
beyond \(\tau\). We have
\begin{equation*}
\theta(\tau)\in\Thetafeasible,
\end{equation*}
since
\(\tau\in\mathcal I_L\) and
\(L(\theta^{\{0\}})\subset\Thetafeasible\).
Therefore, there exists \(\eta_1>0\) such that
\begin{equation*}
\theta(r)\in\Thetafeasible
\qquad
\text{whenever }|r-\tau|<\eta_1,
\end{equation*}
since \(\Thetafeasible\) is open and
\(\theta(\cdot)\) is continuous.

Moreover, there exists \(0<\eta_2\leq\eta_1\) such that
\begin{equation*}
f(\theta(r))<f(\theta^{\{0\}})
\qquad
\text{whenever }|r-\tau|<\eta_2,
\end{equation*}
since \(f(\theta(\cdot))\) is continuous on
\((\tau-\eta_1,\tau+\eta_1)\) and
\eqref{eqn: f theta tau is strictly bounded by f theta initial}
holds.

Thus, we may choose \(\delta'>0\) such that
\begin{equation*}
\delta'<\eta_2
\qquad\text{and}\qquad
\tau+\delta'<1,
\end{equation*}
since \(\eta_2>0\) and \(\tau<1\).
It follows that
\begin{equation*}
\theta(\tau+\delta')
\in
L(\theta^{\{0\}}).
\end{equation*}
Therefore,
\begin{equation*}
\tau+\delta'\in\mathcal I_L,
\end{equation*}
which contradicts
\(\mathcal I_L=[0,\tau]\).
Hence,
\begin{equation*}
\tau=1.
\end{equation*}

Therefore, the entire segment remains in the sublevel set.
In particular,
\begin{equation*}
(\widetilde\theta_1,\theta_2)
=
(\theta_1+h_1,\theta_2)
\in
L(\theta^{\{0\}}).
\end{equation*}

Finally, applying Lemma \ref{lemma: block descent lemma on sublevel set}  to
\(\theta\) and \((\widetilde\theta_1,\theta_2)\) gives
\begin{align*}
f(\widetilde\theta_1,\theta_2)
-
f(\theta_1,\theta_2)
&\leq
-\frac{1}{2L^Q}
\left\|
\nabla_{\theta_1}f(\theta_1,\theta_2)
\right\|_{\Theta_1}^2.
\end{align*}
This proves part (a).

For part (b), let
\(
h_2
:=
-\frac{1}{L^Q}
\nabla_{\theta_2}f(\theta_1,\theta_2)
\)
and define
\(
\theta(t):=(\theta_1,\theta_2+th_2)
\).
Applying the same argument to the second block gives
\begin{equation*}
(\theta_1,\widetilde\theta_2)
\in
L(\theta^{\{0\}})
\end{equation*}
and
\begin{equation*}
f(\theta_1,\theta_2)
-
f(\theta_1,\widetilde\theta_2)
\geq
\frac{1}{2L^Q}
\left\|
\nabla_{\theta_2}f(\theta_1,\theta_2)
\right\|_{\Theta_2}^2.
\end{equation*}
This proves part (b).
\end{proof}
The descent inequalities \eqref{eqn: block descent estimate theta1}, \eqref{eqn: block descent estimate theta2} are standard. 
The main purpose of the Lemma above is to show that the line segment from
\((\theta_1,\theta_2)\) to either
\((\widetilde\theta_1,\theta_2)\) or
\((\theta_1,\widetilde\theta_2)\)
remains in $L(\theta^{\{0\}})$. We argue by contradiction.
Suppose that this line segment exits \(L(\theta^{\{0\}})\), and let
\(\tau<1\) denote the exit time. Up to time \(\tau\), the segment
remains in \(L(\theta^{\{0\}})\), so the block--coordinate Lipschitz
continuity of the gradient applies along this portion of the segment. The
standard descent lemma then gives
\begin{equation*}
f(\theta(\tau))
<
f(\theta)
\leq
f(\theta^{\{0\}}).
\end{equation*}
Since the point at the exit time is feasible and the first inequality is
strict, the openness of \(\Thetafeasible\) and the continuity of \(f\)
imply that points slightly beyond the exit time still belong to
\(L(\theta^{\{0\}})\). This contradicts the definition of \(\tau\).
Hence, the entire line segment remains in \(L(\theta^{\{0\}})\).

Before establishing the linear convergence result, we recall a standard
consequence of strong convexity.
\begin{lemma}
\label{lemma: objective gap bounded by gradient norm from strong convexity}
For every \(\theta\in L(\theta^{\{0\}})\),
\begin{equation}
f(\theta)-f(\theta^*)
\leq
\frac{1}{2\sigma^Q}
\|\nabla f(\theta)\|_\Theta^2.
\label{eqn: objective gap bounded by gradient norm}
\end{equation}
\end{lemma}

\begin{proof}[Proof of Lemma~\ref{lemma: objective gap bounded by gradient norm from strong convexity}]
We have
\(\theta^*\in L(\theta^{\{0\}})\), since
\(f(\theta^*)\leq f(\theta^{\{0\}})\) and
\(\theta^*\in\Thetafeasible\).
Then \eqref{eqn: objective gap bounded by gradient norm} follows from
the standard Polyak--\L{}ojasiewicz inequality implied by strong convexity, applied
with the strong convexity constant \(\sigma^Q\) from
Corollary~\ref{corollary: lip and strong convexity on sublevel set}.
\end{proof}

\begin{proof}[Proof of Theorem~\ref{theorem: linear convergence}]
The proof adapts the argument of Theorem~5.2 in
\cite{BeckAmirConvergenceofBlockCoordinateDescent} to the
sublevel set \(L(\theta^{\{0\}})\).
By Lemma~\ref{lemma: monotone},
\begin{equation*}
\theta^{\{k\}},
\theta^{\{k+\frac12\}}
\in
L(\theta^{\{0\}})
\qquad
\text{for every }k\geq0.
\end{equation*}
Hence, the block descent estimates and the Polyak--\L{}ojasiewicz inequality established
above apply to all BCD iterates.

\textbf{Block descent estimates.}
For \(k\geq1\), we have
\begin{equation}
f(\theta^{\{k\}})
-
f(\theta^{\{k+\frac12\}})
\geq
\frac{1}{2L^Q}
\left\|
\nabla f(\theta^{\{k\}})
\right\|_{\Theta}^2.
\label{eqn: Lip property1}
\end{equation}
For \(k\geq0\), we have
\begin{equation}
f(\theta^{\{k+\frac12\}})
-
f(\theta^{\{k+1\}})
\geq
\frac{1}{2L^Q}
\left\|
\nabla f(\theta^{\{k+\frac12\}})
\right\|_{\Theta}^2.
\label{eqn: Lip property2}
\end{equation}
To prove \eqref{eqn: Lip property1}, we have
\(
\theta^{\{k\}}\in L(\theta^{\{0\}})
\)
by Lemma~\ref{lemma: monotone}. Applying
Lemma~\ref{lemma: block descent estimate on the descent direction} to
\(\theta^{\{k\}}\), define
\begin{equation*}
\widetilde\theta_1^{\{k\}}
:=
\theta_1^{\{k\}}
-
\frac{1}{L^Q}
\nabla_{\theta_1}f(\theta^{\{k\}}).
\end{equation*}
Then
\begin{equation*}
(\widetilde\theta_1^{\{k\}},\theta_2^{\{k\}})
\in
L(\theta^{\{0\}})
\subset
\Thetafeasible,
\end{equation*}
and
\begin{equation*}
f(\theta^{\{k\}})
-
f(\widetilde\theta_1^{\{k\}},\theta_2^{\{k\}})
\geq
\frac{1}{2L^Q}
\left\|
\nabla_{\theta_1}f(\theta^{\{k\}})
\right\|_{\Theta_1}^2.
\end{equation*}
Since
\((\widetilde\theta_1^{\{k\}},\theta_2^{\{k\}})\)
is feasible and differs from \(\theta^{\{k\}}\) only in the first
block, it is feasible for the first-block minimization
\eqref{eq:update-upsi}. Since
\(\theta^{\{k+\frac12\}}\) is the exact minimizer of this block,
\begin{equation*}
f(\theta^{\{k+\frac12\}})
\leq
f(\widetilde\theta_1^{\{k\}},\theta_2^{\{k\}}).
\end{equation*}
Therefore,
\begin{equation*}
f(\theta^{\{k\}})
-
f(\theta^{\{k+\frac12\}})
\geq
\frac{1}{2L^Q}
\left\|
\nabla_{\theta_1}f(\theta^{\{k\}})
\right\|_{\Theta_1}^2.
\end{equation*}
Moreover,
\begin{equation*}
\nabla_{\theta_2}f(\theta^{\{k\}})=0,
\qquad
k\geq1,
\end{equation*}
since \(\theta^{\{k\}}\) is obtained from the exact second-block
minimization. Hence,
\begin{equation*}
\left\|
\nabla f(\theta^{\{k\}})
\right\|_{\Theta}^2
=
\left\|
\nabla_{\theta_1}f(\theta^{\{k\}})
\right\|_{\Theta_1}^2,
\end{equation*}
which proves \eqref{eqn: Lip property1}.

Similarly, comparing the exact second-block minimizer
\(\theta^{\{k+1\}}\) with the feasible second-block gradient step from
Lemma~\ref{lemma: block descent estimate on the descent direction}
gives
\begin{equation*}
f(\theta^{\{k+\frac12\}})
-
f(\theta^{\{k+1\}})
\geq
\frac{1}{2L^Q}
\left\|
\nabla_{\theta_2}f(\theta^{\{k+\frac12\}})
\right\|_{\Theta_2}^2.
\end{equation*}
Moreover,
\begin{equation*}
\nabla_{\theta_1}f(\theta^{\{k+\frac12\}})=0,
\qquad
k\geq0,
\end{equation*}
since \(\theta^{\{k+\frac12\}}\) is obtained from the exact first-block
minimization. Hence,
\begin{equation*}
\left\|
\nabla f(\theta^{\{k+\frac12\}})
\right\|_{\Theta}^2
=
\left\|
\nabla_{\theta_2}f(\theta^{\{k+\frac12\}})
\right\|_{\Theta_2}^2,
\end{equation*}
which proves \eqref{eqn: Lip property2}.

\textbf{Polyak--\L{}ojasiewicz inequality}
For every \(k\geq0\), we have
\begin{align}
f(\theta^{\{k\}})-f(\theta^*)
&\leq
\frac{1}{2\sigma^Q}
\left\|
\nabla f(\theta^{\{k\}})
\right\|_\Theta^2,
\label{eqn: strong convex property integer iterate}\\
f(\theta^{\{k+\frac12\}})-f(\theta^*)
&\leq
\frac{1}{2\sigma^Q}
\left\|
\nabla f(\theta^{\{k+\frac12\}})
\right\|_\Theta^2,
\label{eqn: strong convex property half iterate}
\end{align}
by Lemma~\ref{lemma: objective gap bounded by gradient norm from strong convexity},
since
\(\theta^{\{k\}},\theta^{\{k+\frac12\}}\in L(\theta^{\{0\}})\).
We next combine these with the block descent estimates above.

We first have the contraction from the integer iterate to the half
iterate:
\begin{equation}
f(\theta^{\{k+\frac12\}})-f(\theta^*)
\leq
\rho^Q
\left[
f(\theta^{\{k\}})-f(\theta^*)
\right],
\qquad
k\geq1.
\label{eqn: contraction from integer to half iterate}
\end{equation}
To see this, by
\eqref{eqn: strong convex property integer iterate} and
\eqref{eqn: Lip property1},
\begin{align*}
f(\theta^{\{k\}})-f(\theta^*)
&\leq
\frac{1}{2\sigma^Q}
\left\|
\nabla f(\theta^{\{k\}})
\right\|_{\Theta}^2\\
&\leq
\frac{L^Q}{\sigma^Q}
\left[
f(\theta^{\{k\}})
-
f(\theta^{\{k+\frac12\}})
\right]\\
&=
\frac{L^Q}{\sigma^Q}
\left\{
\left[
f(\theta^{\{k\}})-f(\theta^*)
\right]
-
\left[
f(\theta^{\{k+\frac12\}})-f(\theta^*)
\right]
\right\}.
\end{align*}
Rearranging gives
\eqref{eqn: contraction from integer to half iterate}, since
\(
\rho^Q=1-\frac{\sigma^Q}{L^Q}.
\)

We next have the contraction from the half iterate to the next integer
iterate:
\begin{equation}
f(\theta^{\{k+1\}})-f(\theta^*)
\leq
\rho^Q
\left[
f(\theta^{\{k+\frac12\}})-f(\theta^*)
\right],
\qquad
k\geq0.
\label{eqn: contraction from half to integer iterate}
\end{equation}
This follows by the same argument as
\eqref{eqn: contraction from integer to half iterate}, using
\eqref{eqn: strong convex property half iterate} and
\eqref{eqn: Lip property2}.

Combining
\eqref{eqn: contraction from integer to half iterate} and
\eqref{eqn: contraction from half to integer iterate}, we obtain
\begin{equation}
f(\theta^{\{k+1\}})-f(\theta^*)
\leq
(\rho^Q)^2
\left[
f(\theta^{\{k\}})-f(\theta^*)
\right],
\qquad
k\geq1.
\label{eqn: full iteration contraction}
\end{equation}

We first prove the bound for the integer iterates in
\eqref{eqn: linear conv bound integer iterate}.
For \(k=0\), we have
\begin{equation*}
f(\theta^{\{\frac12\}})-f(\theta^*)
\leq
f(\theta^{\{0\}})-f(\theta^*),
\end{equation*}
by Lemma~\ref{lemma: monotone}. Therefore,
\eqref{eqn: contraction from half to integer iterate} with \(k=0\)
gives
\begin{align*}
f(\theta^{\{1\}})-f(\theta^*)
&\leq
\rho^Q
\left[
f(\theta^{\{\frac12\}})-f(\theta^*)
\right]\\
&\leq
\rho^Q
\left[
f(\theta^{\{0\}})-f(\theta^*)
\right].
\end{align*}
Hence,
\eqref{eqn: linear conv bound integer iterate} holds for \(k=0\).

For \(k\geq1\), iterating
\eqref{eqn: full iteration contraction} gives
\begin{equation*}
f(\theta^{\{k+1\}})-f(\theta^*)
\leq
(\rho^Q)^{2k}
\left[
f(\theta^{\{1\}})-f(\theta^*)
\right].
\end{equation*}
Using the preceding bound for
\(f(\theta^{\{1\}})-f(\theta^*)\), we obtain
\begin{equation*}
f(\theta^{\{k+1\}})-f(\theta^*)
\leq
(\rho^Q)^{2k+1}
\left[
f(\theta^{\{0\}})-f(\theta^*)
\right].
\end{equation*}
Thus,
\eqref{eqn: linear conv bound integer iterate} holds for every
\(k\geq0\).

We next prove the bound for the half iterates in
\eqref{eqn: linear conv bound half iterate}.
For \(k=0\), we have
\begin{equation*}
f(\theta^{\{\frac12\}})-f(\theta^*)
\leq
f(\theta^{\{0\}})-f(\theta^*),
\end{equation*}
by Lemma~\ref{lemma: monotone}.
For \(k\geq1\),
\eqref{eqn: contraction from integer to half iterate} and
\eqref{eqn: linear conv bound integer iterate} with \(k-1\) give
\begin{align*}
f(\theta^{\{k+\frac12\}})-f(\theta^*)
&\leq
\rho^Q
\left[
f(\theta^{\{k\}})-f(\theta^*)
\right]\\
&\leq
\rho^Q(\rho^Q)^{2k-1}
\left[
f(\theta^{\{0\}})-f(\theta^*)
\right]\\
&=
(\rho^Q)^{2k}
\left[
f(\theta^{\{0\}})-f(\theta^*)
\right].
\end{align*}
Thus,
\eqref{eqn: linear conv bound half iterate} holds for every
\(k\geq0\).

Finally, we show that
\begin{equation*}
\rho^Q\in(0,1).
\end{equation*}
We have
\begin{equation*}
0<
\frac{\sigma^Q}{L^Q}
=
\frac{a_Q^2}{2b_Q^2}
<
\frac12,
\end{equation*}
by
\(
L^Q=2a_Q^{-2},
\)
\(
\sigma^Q=b_Q^{-2},
\)
and
Lemma \ref{lemma: uniform spectral bound on level set}.
Hence,
\begin{equation*}
\frac12<
\rho^Q
=
1-\frac{\sigma^Q}{L^Q}
<1.
\end{equation*}
\end{proof}

\begin{proof}[Proof of Corollary~\ref{corollary: linear convergence of the dual iterations}]
Let
\begin{equation*}
\Delta_0
:=
f(\theta^{\{0\}})-f(\theta^*).
\end{equation*}
By Theorem~\ref{theorem: linear convergence}, for every \(k\geq1\),
\begin{align*}
f(\theta^{\{k\}})-f(\theta^*)
&\leq
\Delta_0(\rho^Q)^{2k-1}\\
&=
\frac{\Delta_0}{\rho^Q}\rho^k.
\end{align*}
For the iteration \(k=0\), since \(\rho^Q\in(\frac12,1)\),
\begin{equation*}
f(\theta^{\{0\}})-f(\theta^*)
=
\Delta_0
\leq
\frac{\Delta_0}{\rho^Q}.
\end{equation*}
Therefore, for every \(k\geq0\),
\begin{equation}
f(\theta^{\{k\}})-f(\theta^*)
\leq
\frac{\Delta_0}{\rho^Q}\rho^k.
\label{eqn: objective bound for dual iterates}
\end{equation}

Moreover, Theorem~\ref{theorem: linear convergence} gives directly
\begin{equation}
f(\theta^{\{k+\frac12\}})-f(\theta^*)
\leq
\Delta_0\rho^k,
\qquad
k\geq0.
\label{eqn: objective bound for half dual iterates}
\end{equation}

We have
\(
\theta^*\in L(\theta^{\{0\}})
\),
since
\(
f(\theta^*)\leq f(\theta^{\{0\}})
\) and \(\theta^*\in\Thetafeasible\) .
By Corollary~\ref{corollary: lip and strong convexity on sublevel set},
for any
\(\theta\in L(\theta^{\{0\}})\),
\begin{align*}
f(\theta)
&\geq
f(\theta^*)
+
\left\langle
\nabla f(\theta^*),
\theta-\theta^*
\right\rangle_\Theta
+
\frac{\sigma^Q}{2}
\|\theta-\theta^*\|_\Theta^2.
\end{align*}
Moreover,
\begin{equation*}
\nabla f(\theta^*)=0,
\end{equation*}
since \(\theta^*\) is the minimizer of \(f\) on the open set
\(\Thetafeasible\).
Therefore,
\begin{equation*}
f(\theta)-f(\theta^*)
\geq
\frac{\sigma^Q}{2}
\|\theta-\theta^*\|_\Theta^2,
\qquad
\theta\in L(\theta^{\{0\}}).
\end{equation*}

Setting $\theta=\theta^{\{k\}}$ and combining this inequality with
\eqref{eqn: objective bound for dual iterates}, we obtain, for every
\(k\geq0\),
\begin{equation*}
\|\theta^{\{k\}}-\theta^*\|_\Theta^2
\leq
\frac{2\Delta_0}
{\rho^Q\sigma^Q}
\rho^k.
\end{equation*}
Similarly, by
\eqref{eqn: objective bound for half dual iterates},
\begin{equation*}
\|\theta^{\{k+\frac12\}}-\theta^*\|_\Theta^2
\leq
\frac{2\Delta_0}
{\sigma^Q}
\rho^k
\leq
\frac{2\Delta_0}
{\rho^Q\sigma^Q}
\rho^k,
\end{equation*}
where the last inequality follows from \(\rho^Q<1\).

Thus, choosing \(C_0>0\) such that
\begin{equation*}
C_0
\geq
\frac{2\Delta_0}
{\rho^Q\sigma^Q}
\end{equation*}
gives the claimed result.
\end{proof}
\section{Supplementary Analysis for Section~\ref{sec: stability}}\label{sec: supp of stabiility}
For any pair of marginal covariance matrices satisfying Assumption \ref{assumption:maringal}, we use \(L\) to denote the
corresponding block-diagonal Cholesky matrix, following
\eqref{eqn: marginal cholesky factors} and
\eqref{eqn: def of L for two marginal cholesky}.
Define
\begin{equation*}
\mathcal L^{n}
:=
\left\{
M\in\mathbb R^{n\times n}:
M \text{ is lower triangular}
\right\},
\end{equation*}
and 
\begin{equation*}
\mathcal L_{++}^{n}
:=
\left\{
M\in\mathcal L^{n}:
M_{ii}>0,\quad i=1,\ldots,n
\right\}.
\end{equation*}
Define
\begin{equation*}
\mathcal L_{\mathrm{blk-diag}}^{2Td}
:=
\left\{
\begin{bmatrix}
M^\mu&0\\
0&M^\nu
\end{bmatrix}
:
M^\mu,M^\nu\in\mathcal L^{Td}
\right\},
\end{equation*}
and
\begin{equation*}
\mathcal L_{\mathrm{blk-diag},++}^{2Td}
:=
\left\{
\begin{bmatrix}
M^\mu&0\\
0&M^\nu
\end{bmatrix}
:
M^\mu,M^\nu\in\mathcal L_{++}^{Td}
\right\}.
\end{equation*}
By \eqref{eqn: marginal cholesky factors} and
\eqref{eqn: def of L for two marginal cholesky},
we know 
\begin{equation}
L\in\mathcal L_{\mathrm{blk-diag},++}^{2Td}.
\label{eqn: population L belongs to positive cholesky space}
\end{equation}

We now make explicit the dependence of the matrices introduced in
Section~\ref{sec: optimization in cov space} on the marginal covariance
matrices through their Cholesky factors. Following
\eqref{eqn:primal loss}, define
\begin{equation}
Q(L)
:=
L^\top(\Sigma^\gamma)^{-1}L+I_{2Td},
\label{eqn: Q depending on cholesky matrix L}
\end{equation}
where \(\Sigma^\gamma\) denotes the covariance matrix of the quadratic
cost in the original coordinates, as in Assumption~\ref{assumption:cost}.
To see that \eqref{eqn: Q depending on cholesky matrix L} agrees with
\eqref{eqn:primal loss}, recall
from \eqref{eqn: quadratic c tilde in cholesky} that the transformed
cost covariance corresponding to \(L\) is $L^{-1}\Sigma^\gamma L^{-\top}.$
Therefore, its inverse is $L^\top(\Sigma^\gamma)^{-1}L,$
which gives \eqref{eqn: Q depending on cholesky matrix L}.

Recalling the linear map \(A\) of the dual variables \(\theta\), defined
in \eqref{eqn: A from theta}, define
\begin{equation}
Z(\theta;L)
:=
Q(L)+A(\theta),
\qquad
\theta\in\Theta,
\label{eqn: Z depending on theta and cholesky matrix L}
\end{equation}
and the corresponding feasible dual domain, 
\begin{equation}
\Thetafeasible(L)
:=
\left\{
\theta\in\Theta:
Z(\theta;L)\succ0
\right\}.
\label{eqn: feasible dual domain depending on L}
\end{equation}
For
\(\theta=(\Lambda_u,\Lambda_y,\Psi)\in\Thetafeasible(L)\), as in \eqref{eqn: dual objective}, define
\begin{equation}
f(\theta;L)
:=
-\frac12\log\det Z(\theta;L)
+
\operatorname{Tr}(\Lambda_u)
+
\operatorname{Tr}(\Lambda_y)
-
Td.
\label{eqn: dual objective depending on L}
\end{equation}
Finally, denote the unique dual optimizer corresponding to \(L\) by
\begin{equation}
\theta^*(L)
:=
\underset{\theta\in\Thetafeasible(L)}{\arg\min}\,
f(\theta;L).
\label{eqn: dual optimizer as function of L}
\end{equation}
The existence and uniqueness of \(\theta^*(L)\) follow from
Theorem~\ref{theorem: dual problem of ecot gaussian}.
Define
\begin{equation}
m_{Z^*}(L)
:=
\lambda_{\min}
\left(
Z(\theta^*(L);L)
\right)
>0.
\label{eqn: m Z star depending on L}
\end{equation}
The strict positivity follows from
\(
\theta^*(L)\in\Thetafeasible(L)
\) and \eqref{eqn: feasible dual domain depending on L}.

We recall the Cholesky perturbation bound from
\cite[Theorem~1.4]{sun1991perturbation}, specialized to the real-valued setting.
\begin{lemma}[Cholesky perturbation bound {\cite[Theorem~1.4]{sun1991perturbation}}]
\label{lemma: cholesky perturbation bound}
Let \(S\in\mathbb S_{++}^n\), and let $S=BB^\top$
be its Cholesky factorization, where \(B\in\mathcal L_{++}^n\).
Suppose that \(\Delta S\in\mathbb S^n\)
satisfies
\begin{equation}
\label{eqn: cholesky perturbation condition in general}
\|S^{-1}\|_2\|\Delta S\|_2<1.
\end{equation}
Then \(S+\Delta S\in\mathbb S_{++}^n\) admits the unique Cholesky
factorization $S+\Delta S
=
(B+\Delta B)(B+\Delta B)^\top,$
where \(B+\Delta B\in\mathcal L_{++}^n\).
Moreover,
\begin{equation}
\label{eqn: cholesky perturbation bound in general}
{\|\Delta B\|_F}
\leq
\frac{{\|B\|_p}
\kappa_2(S)
}{
\sqrt{2}
\left(
1-\|S^{-1}\|_2\|\Delta S\|_2
\right)
\|S\|_p
}\|\Delta S\|_F,
\qquad
p=2,F,
\end{equation}
where $\kappa_2(S):=\|S\|_2\|S^{-1}\|_2.$
\end{lemma}
\begin{lemma}
\label{lemma: local stability of dual optimizer under cholesky perturbation}
There exist constants
\begin{equation*}
r_\theta(L)>0,
\qquad
C_\theta(L)\geq0,
\end{equation*}
depending on \(L,\Sigma^\gamma,T,\) and \(d\), such that for every $\Delta L
\in
\mathcal L_{\mathrm{blk-diag}}^{2Td}$
satisfying
\begin{equation}
\|\Delta L\|_F
\leq
r_\theta(L),
\label{eqn: IFT cholesky perturbation radius}
\end{equation}
the following hold:
\begin{align}
L+\Delta L
&\in
\mathcal L_{\mathrm{blk-diag},++}^{2Td},
\label{eqn: perturbed L remains positive cholesky}\\
\left\|
\theta^*(L+\Delta L)-\theta^*(L)
\right\|_\Theta
&\leq
C_\theta(L)\|\Delta L\|_F,
\label{eqn: dual optimizer locally Lipschitz in cholesky factor}\\
\lambda_{\min}
\left(
Z\left(
\theta^*(L+\Delta L);
L+\Delta L
\right)
\right)
&\geq
\frac12m_{Z^*}(L).
\label{eqn: perturbed optimal Z uniform positive lower bound}
\end{align}
\end{lemma}
\begin{proof}[Proof of Lemma~\ref{lemma: local stability of dual optimizer under cholesky perturbation}]
We prove the three conclusions separately and choose a common
radius \(r_\theta(L)\).

\begin{enumerate}[label=\alph*.]

\item
We first show that \(L+\Delta L\) remains in
\(\mathcal L_{\mathrm{blk-diag},++}^{2Td}\) for all sufficiently small
\(\Delta L\in\mathcal L_{\mathrm{blk-diag}}^{2Td}\).

By \eqref{eqn: population L belongs to positive cholesky space},
\(L\in\mathcal L_{\mathrm{blk-diag},++}^{2Td}\). Hence all diagonal entries
of \(L\) are strictly positive. Define
\begin{equation}
r_{\mathrm{diag}}(L)
:=
\frac12\min_{1\leq i\leq 2Td}L_{ii}
>0,
\label{eqn: radius preserving positive diagonal of L}
\end{equation}
where $M_{ii}$ is the $i-$th diagonal component of a square matrix $M$.
For any
\(\Delta L\in\mathcal L_{\mathrm{blk-diag}}^{2Td}\) satisfying
\(\|\Delta L\|_F\leq r_{\mathrm{diag}}(L)\), we have, for \(i=1,\ldots,2Td\), the $i-$th diagonal component of $L+\Delta L$ satisfies
\begin{align*}
L_{ii}+(\Delta L)_{ii}
&\geq
L_{ii}-|(\Delta L)_{ii}|\\
&\geq
L_{ii}-r_{\mathrm{diag}}(L)\\
&\geq
\frac12 L_{ii}
>0.
\end{align*}
Moreover, \(L+\Delta L\in\mathcal L_{\mathrm{blk-diag}}^{2Td}\). Hence,
\begin{equation}
L+\Delta L
\in
\mathcal L_{\mathrm{blk-diag},++}^{2Td}
\qquad
\text{whenever }
\|\Delta L\|_F\leq r_{\mathrm{diag}}(L).
\label{eqn: local radius preserves positive cholesky space}
\end{equation}

\item
We next show that the optimizer map is locally Lipschitz with respect to
the Cholesky matrix. To that end, we apply the Implicit Function Theorem to show the map $\theta^*(\cdot)$ is continuously differentiable on an open domain of $L$.

Equip \(\mathcal L_{\mathrm{blk-diag}}^{2Td}\) with the Frobenius norm and recall that
\(\Theta\) is equipped with the norm \(\|\cdot\|_\Theta\). Since both spaces are
finite-dimensional, they are Banach spaces.

Define
\begin{equation}
\mathcal O
:=
\left\{
(L',\theta)
\in
\mathcal L_{\mathrm{blk-diag},++}^{2Td}\times\Theta:
Z(\theta;L')\succ0
\right\}.
\label{eqn: open domain for IFT stability}
\end{equation}
The set \(\mathcal O\) is open in
\(\mathcal L_{\mathrm{blk-diag}}^{2Td}\times\Theta\).
First,
\(\mathcal L_{\mathrm{blk-diag},++}^{2Td}\) is open in
\(\mathcal L_{\mathrm{blk-diag}}^{2Td}\).
To see this, take any
\(L'\in\mathcal L_{\mathrm{blk-diag},++}^{2Td}\) and define
$\varepsilon
:=
\frac12
\min_{1\leq i\leq 2Td}L'_{ii}
>0.$ We can show 
\begin{equation*}
    L'+H\in\mathcal L_{\mathrm{blk-diag},++}^{2Td},
\end{equation*}
for any
\(H\in\mathcal L_{\mathrm{blk-diag}}^{2Td}\) satisfying
\(\|H\|_F<\varepsilon\). This is because
\begin{equation*}
L'_{ii}+H_{ii}
>
L'_{ii}-\varepsilon
\geq
\frac12L'_{ii}
>0,
\qquad
i=1,\ldots,2Td,
\end{equation*}
where we use
\(
|H_{ii}|\leq\|H\|_F<\varepsilon
\)
for each $i$.
Since
\(L'+H\in\mathcal L_{\mathrm{blk-diag}}^{2Td}\),
it follows that
\(L'+H\in\mathcal L_{\mathrm{blk-diag},++}^{2Td}\).
Thus,
\(\mathcal L_{\mathrm{blk-diag},++}^{2Td}\) is open in
\(\mathcal L_{\mathrm{blk-diag}}^{2Td}\). Second,
\begin{equation*}
\left\{
(L',\theta):
Z(\theta;L')\succ0
\right\}
\end{equation*}
is open in
\(\mathcal L_{\mathrm{blk-diag}}^{2Td}\times\Theta\), because this is the preimage under the map
\(
(L',\theta)\mapsto Z(\theta;L')
\)
which is continuous by
\eqref{eqn: Q depending on cholesky matrix L},
\eqref{eqn: Z depending on theta and cholesky matrix L}, and
\eqref{eqn: A from theta}, and \(\mathbb S_{++}^{2Td}\) is open in
\(\mathbb S^{2Td}\).
Therefore,
\(\mathcal O\), being the intersection of this open set with
\(
\mathcal L_{\mathrm{blk-diag},++}^{2Td}\times\Theta
\),
is open in
\(\mathcal L_{\mathrm{blk-diag}}^{2Td}\times\Theta\).

Define
\begin{equation}
F:\mathcal O\to\Theta,
\qquad
F(L',\theta)
:=
\nabla_\theta f(\theta;L').
\label{eqn: stationary map for IFT stability}
\end{equation}
By \eqref{eqn: dual optimizer as function of L},
\(\theta^*(L)\) is the optimizer of \(f(\cdot;L)\). Moreover,
\(\Thetafeasible(L)\) is open and
\(\theta^*(L)\in\Thetafeasible(L)\). Therefore, the first-order
optimality condition gives
\begin{equation}
F(L,\theta^*(L))
=
\nabla_\theta f(\theta^*(L);L)
=
0.
\label{eqn: stationary map vanishes at optimizer L}
\end{equation}
We next verify that \(F\) is continuously differentiable on
\(\mathcal O\).
By \eqref{eqn: Q depending on cholesky matrix L}, the map $L'
\longmapsto
Q(L')$
from \(\mathcal L_{\mathrm{blk-diag}}^{2Td}\) to
\(\mathbb S^{2Td}\) is polynomial in the entries of \(L'\), and hence
is \(C^\infty\). By \eqref{eqn: A from theta}, the map
 $\theta
\longmapsto
A(\theta)$
from \(\Theta\) to \(\mathbb S^{2Td}\) is linear, and hence is
\(C^\infty\). Therefore, by
\eqref{eqn: Z depending on theta and cholesky matrix L}, the map $(L',\theta)
\longmapsto
Z(\theta;L')$
from
\(\mathcal L_{\mathrm{blk-diag}}^{2Td}\times\Theta\) to
\(\mathbb S^{2Td}\) is \(C^\infty\). Since restricting a
\(C^\infty\) map to an open subset of its domain preserves
\(C^\infty\)-smoothness, its restriction to \(\mathcal O\) is also
\(C^\infty\).

Since \(Z(\theta;L')\succ0\) for every
\((L',\theta)\in\mathcal O\), and
\(\log\det(\cdot)\) is \(C^\infty\) on
\(\mathbb S_{++}^{2Td}\), it follows from
\eqref{eqn: dual objective depending on L} that the map
\begin{equation*}
(L',\theta)
\longmapsto
f(\theta;L')
\end{equation*}
from \(\mathcal O\) to \(\mathbb R\) is \(C^\infty\).
Consequently,
\begin{equation*}
F:\mathcal O\to\Theta,
\qquad
F(L',\theta)
=
\nabla_\theta f(\theta;L'),
\end{equation*}
is \(C^\infty\), and in particular continuously differentiable, on
\(\mathcal O\).
 By \eqref{eqn: population L belongs to positive cholesky space} and
\eqref{eqn: m Z star depending on L},
\begin{equation*}
(L,\theta^*(L))\in\mathcal O.
\end{equation*}

We next verify the invertibility condition in the Implicit Function
Theorem. Consider the linear operator \(D_\theta F(L,\theta^*(L)):\Theta\to\Theta\), the Fr\'echet derivative with respect to \(\theta\). This operator is injective, because by
\eqref{eqn: Hessian bilinear form via gradient derivative},
for every \(\delta\theta\in\Theta\),
\begin{align*}
\left\langle
\delta\theta,
D_\theta F(L,\theta^*(L))[\delta\theta]
\right\rangle_{\Theta}
&=
\left\langle
\delta\theta,
D_\theta\nabla_\theta
f(\theta^*(L);L)[\delta\theta]
\right\rangle_{\Theta}\\
&=
\nabla^2_{\theta\theta}f(\theta^*(L);L)
[\delta\theta,\delta\theta],
\end{align*}
which is strictly positive for every
\(\delta\theta\neq0\), by Lemma~\ref{lemma: hessian of f}.

Therefore,
\(D_\theta F(L,\theta^*(L))\) is injective. Since it is a linear
operator from the finite-dimensional space \(\Theta\) to itself,
injectivity implies surjectivity, and hence it is bijective. Thus its
inverse exists and is linear. And this inverse operator is bounded since every linear operator between
finite-dimensional normed spaces is bounded.


Thus, we may apply the Implicit Function Theorem (see \cite[Theorem 5.7]{Sideris2013Ordinary}) to
\(F\) at
\((L,\theta^*(L))\).

There exist open neighborhoods
\(\mathcal V_L\subset\mathcal L_{\mathrm{blk-diag}}^{2Td}\) of \(L\) and
\(\mathcal V_\theta\subset\Theta\) of \(\theta^*(L)\),  such that
\begin{equation*}
\mathcal V_L\times\mathcal V_\theta
\subset\mathcal O,
\end{equation*}
and a unique continuously differentiable map 
\begin{equation}
\phi:\mathcal V_L\to\mathcal V_\theta
\label{eqn: IFT local optimizer map}
\end{equation}
such that
\begin{equation}
\phi(L)
=
\theta^*(L)
\label{eqn: IFT local map agrees at L}
\end{equation}
and
\begin{equation}
F(L',\phi(L'))
=
\nabla_\theta
f(\phi(L');L')
=
0,
\qquad
L'\in\mathcal V_L.
\label{eqn: IFT local map stationary}
\end{equation}

We now identify \(\phi\) with the optimizer map.
For every \(L'\in\mathcal V_L\),
Theorem~\ref{theorem: dual problem of ecot gaussian} gives the unique
optimizer \(\theta^*(L')\), and hence
\begin{equation*}
\nabla_\theta
f(\theta^*(L');L')
=
0.
\end{equation*}
On the other hand, \(\phi(L')\) satisfies
\eqref{eqn: IFT local map stationary}.
By Lemma~\ref{lemma: hessian of f},
\(f(\cdot;L')\) is strictly convex with respect to $\theta$ on
\(\Thetafeasible(L')\), and therefore has at most one stationary point.
Consequently,
\begin{equation}
\phi(L')
=
\theta^*(L'),
\qquad
L'\in\mathcal V_L.
\label{eqn: IFT map equals optimizer map}
\end{equation}
Hence,
\begin{equation*}
L'
\longmapsto
\theta^*(L')
\end{equation*}
is continuously differentiable on \(\mathcal V_L\), hence its derivative is continuous on \(\mathcal V_L\).

Since \(\mathcal V_L\) is open and contains \(L\), there exists
\(h_\theta(L)>0\), depending on
\(L,\Sigma^\gamma,T,\) and \(d\), such that
\begin{equation}
\left\{
L'\in\mathcal L_{\mathrm{blk-diag}}^{2Td}:
\|L'-L\|_F
\leq
h_\theta(L)
\right\}
\subset
\mathcal V_L.
\label{eqn: closed ball contained in IFT neighborhood}
\end{equation}
Since \(\mathcal L_{\mathrm{blk-diag}}^{2Td}\) is finite-dimensional, the
closed and bounded set in
\eqref{eqn: closed ball contained in IFT neighborhood} is compact.
And the derivative $D\theta^*(L')$ is continuous on the compact set in
\eqref{eqn: closed ball contained in IFT neighborhood}. Since the
operator norm is continuous, the map
\(
L'\mapsto \|D\theta^*(L')\|_{\mathrm{op}}
\)
is a continuous real--valued function on this compact set. Therefore, by
the Extreme Value Theorem, it attains a finite maximum. Define
\begin{equation}
C_\theta(L)
:=
\max_{\substack{
L'\in\mathcal L_{\mathrm{blk-diag}}^{2Td}\\
\|L'-L\|_F\leq h_\theta(L)
}}
\|D\theta^*(L')\|_{\mathrm{op}}
<\infty.
\label{eqn: definition of C theta L}
\end{equation}

Now let
\(\Delta L\in\mathcal L_{\mathrm{blk-diag}}^{2Td}\) satisfy
\(
\|\Delta L\|_F\leq h_\theta(L)
\).
The Mean value theorem (see \cite[Theorem~1.8]{ambrosetti1995primer}) gives
\begin{align*}
\left\|
\theta^*(L+\Delta L)-\theta^*(L)
\right\|_\Theta
&\leq
\sup_{L'\in[L,L+\Delta L]}
\left\|
D_G\theta^*(L')
\right\|_{\mathrm{op}}
\|\Delta L\|_F\\
&=
\sup_{L'\in[L,L+\Delta L]}
\left\|
D\theta^*(L')
\right\|_{\mathrm{op}}
\|\Delta L\|_F\\
&\leq
C_\theta(L)\|\Delta L\|_F,
\end{align*}
 Here \(D_G\theta^*(L')\) denotes the G\^ateaux derivative, and the
equality follows because \(\theta^*(\cdot)\) is continuously
Fr\'echet differentiable on \(\mathcal V_L\). And the second inequality follows from \eqref{eqn: definition of C theta L} and
\begin{equation*}
    [L,L+\Delta L]
\subset
\left\{
L'\in\mathcal L_{\mathrm{blk-diag}}^{2Td}:
\|L'-L\|_F\leq h_\theta(L)
\right\}
\end{equation*}

This proves
\begin{equation}
\left\|
\theta^*(L+\Delta L)-\theta^*(L)
\right\|_\Theta
\leq
C_\theta(L)\|\Delta L\|_F,
\qquad
\|\Delta L\|_F\leq h_\theta(L).
\label{eqn: local optimizer perturbation before final radius}
\end{equation}

\item
We finally establish a uniform positive lower bound for the smallest
eigenvalue of the perturbed optimal \(Z\).

By \eqref{eqn: IFT map equals optimizer map}, the map
\(L'\mapsto\theta^*(L')\) is continuous on \(\mathcal V_L\).
Together with the continuity of
\((\theta,L')\mapsto Z(\theta;L')\), established above, and the
continuity of the minimum eigenvalue, the map
\begin{equation*}
L'
\longmapsto
\lambda_{\min}
\left(
Z(\theta^*(L');L')
\right)
\end{equation*}
is continuous at \(L\).

At \(L\), its value is \(m_{Z^*}(L)>0\) by
\eqref{eqn: m Z star depending on L}. Hence there exists
\(h_Z(L)>0\), depending on \(L,\Sigma^\gamma,T,\) and \(d\), such that if
\(
\|L'-L\|_F\leq h_Z(L)
\),
\begin{equation*}
    |\lambda_{\min}
\left(
Z(\theta^*(L');L')
\right)-m_{Z^*}(L)|
<
\frac12m_{Z^*}(L),
\end{equation*}
hence
\begin{equation}
\lambda_{\min}
\left(
Z(\theta^*(L');L')
\right)\geq
\frac12m_{Z^*}(L).
\label{eqn: local spectral lower bound before final radius}
\end{equation}

\end{enumerate}

Define
\begin{equation}
r_\theta(L)
:=
\min
\left\{
r_{\mathrm{diag}}(L),
h_\theta(L),
h_Z(L)
\right\}
>0.
\label{eqn: definition of final IFT stability radius}
\end{equation}
Then every
\(\Delta L\in\mathcal L_{\mathrm{blk-diag}}^{2Td}\) satisfying
\eqref{eqn: IFT cholesky perturbation radius} also satisfies the
hypotheses of
\eqref{eqn: local radius preserves positive cholesky space},
\eqref{eqn: local optimizer perturbation before final radius}, and
\eqref{eqn: local spectral lower bound before final radius}.
These three inequalities give
\eqref{eqn: perturbed L remains positive cholesky},
\eqref{eqn: dual optimizer locally Lipschitz in cholesky factor}, and
\eqref{eqn: perturbed optimal Z uniform positive lower bound},
respectively.
\end{proof}

\begin{proof}[Proof of Theorem \ref{theorem: stability under marginal covariance perturbations}]
Let
\(
L=\operatorname{diag}(L^\mu,L^\nu)
\)
be the block-diagonal Cholesky matrix associated with \(\Sigma_0\), as in
\eqref{eqn: marginal cholesky factors} and
\eqref{eqn: def of L for two marginal cholesky}.
We prove the result by propagating the perturbation
\(\Delta\Sigma_0\) through the Cholesky factor, the dual optimizer, the
matrix \(Z\), and finally the optimal covariance.

\begin{enumerate}[label=\alph*.]

\item
We first control the perturbation of the Cholesky factor.
Assume
\begin{equation}
\|\Delta\Sigma_0\|_F
\leq
r_{\mathrm{chol}}(\Sigma_0)
:=
\frac{1}{
2\max\left\{
\|(\Sigma^\mu)^{-1}\|_2,
\|(\Sigma^\nu)^{-1}\|_2
\right\}
}
>0.
\label{eqn: radius for marginal cholesky perturbation}
\end{equation}
Then, for each \(\alpha\in\{\mu,\nu\}\),
we have
\begin{align*}
\numberthis\label{eqn: Sigma alpg inv 2 norm times delta sigma alpha 2 nrom bounded by half} \|(\Sigma^\alpha)^{-1}\|_2 
\|\Delta\Sigma^\alpha\|_2 &\leq \|(\Sigma^\alpha)^{-1}\|_2 
\|\Delta\Sigma^\alpha\|_F\\
&\leq  \|(\Sigma^\alpha)^{-1}\|_2 
\|\Delta\Sigma_0\|_F\\
&\leq \frac12
\end{align*}

Therefore, the condition
\eqref{eqn: cholesky perturbation condition in general} of
Lemma~\ref{lemma: cholesky perturbation bound} holds for both
\(\Sigma^\mu\) and \(\Sigma^\nu\).
It follows that
\(
\Sigma^\mu+\Delta\Sigma^\mu\succ0
\)
and
\(
\Sigma^\nu+\Delta\Sigma^\nu\succ0
\).
Let
\(
L^\mu+\Delta L^\mu
\)
and
\(
L^\nu+\Delta L^\nu
\)
be their Cholesky factors, and define
\(
\Delta L
:=
\operatorname{diag}(\Delta L^\mu,\Delta L^\nu)
\).
We next obtain
\begin{equation}
\|\Delta L\|_F
\leq
C_{\mathrm{chol}}(\Sigma_0)
\|\Delta\Sigma_0\|_F,
\label{eqn: marginal covariance perturbation implies cholesky perturbation}
\end{equation}
where
\begin{equation}
C_{\mathrm{chol}}(\Sigma_0)
:=
\sqrt{2}
\max_{\alpha\in\{\mu,\nu\}}
\left\{
\|L^\alpha\|_2
\|(\Sigma^\alpha)^{-1}\|_2
\right\}.
\label{eqn: constant for marginal cholesky perturbation}
\end{equation}
For each \(\alpha\in\{\mu,\nu\}\), applying
\eqref{eqn: cholesky perturbation bound in general} with \(p=2\) gives
\begin{align*}
\|\Delta L^\alpha\|_F
&\leq
\frac{
\|L^\alpha\|_2
\kappa_2(\Sigma^\alpha)
}{
\sqrt{2}
\left(
1-\|(\Sigma^\alpha)^{-1}\|_2
\|\Delta\Sigma^\alpha\|_2
\right)
\|\Sigma^\alpha\|_2
}
\|\Delta\Sigma^\alpha\|_F\\
&=
\frac{
\|L^\alpha\|_2
\|(\Sigma^\alpha)^{-1}\|_2
}{
\sqrt{2}
\left(
1-\|(\Sigma^\alpha)^{-1}\|_2
\|\Delta\Sigma^\alpha\|_2
\right)
}
\|\Delta\Sigma^\alpha\|_F\\
&\leq
\sqrt{2}
\|L^\alpha\|_2
\|(\Sigma^\alpha)^{-1}\|_2
\|\Delta\Sigma^\alpha\|_F\\
&\leq
C_{\mathrm{chol}}(\Sigma_0)
\|\Delta\Sigma^\alpha\|_F.
\end{align*}
Here, the equality follows from
\(
\kappa_2(\Sigma^\alpha)
=
\|\Sigma^\alpha\|_2
\|(\Sigma^\alpha)^{-1}\|_2
\),
and the third line follows from \eqref{eqn: Sigma alpg inv 2 norm times delta sigma alpha 2 nrom bounded by half}.

Therefore,
\begin{align*}
\|\Delta L\|_F^2
&=
\|\Delta L^\mu\|_F^2
+
\|\Delta L^\nu\|_F^2\\
&\leq
C_{\mathrm{chol}}(\Sigma_0)^2
\left(
\|\Delta\Sigma^\mu\|_F^2
+
\|\Delta\Sigma^\nu\|_F^2
\right)\\
&=
C_{\mathrm{chol}}(\Sigma_0)^2
\|\Delta\Sigma_0\|_F^2.
\end{align*}
Taking square roots gives
\eqref{eqn: marginal covariance perturbation implies cholesky perturbation}.
\item
We next control the perturbation of the dual optimizer.
Assume
\begin{equation}
\|\Delta\Sigma_0\|_F
\leq
r_{\mathrm{IFT}}(\Sigma_0)
:=
\frac{r_\theta(L)}
{C_{\mathrm{chol}}(\Sigma_0)}
>0,
\label{eqn: covariance radius for applying IFT stability lemma}
\end{equation}
where \(r_\theta(L)\) is the perturbation radius in
Lemma~\ref{lemma: local stability of dual optimizer under cholesky perturbation}.
Then, by
\eqref{eqn: marginal covariance perturbation implies cholesky perturbation},
\begin{align*}
\numberthis\label{eqn: delta L F norm is bouned by r theta (L) proved from assumption in staiblity thoerem}\|\Delta L\|_F
&\leq
C_{\mathrm{chol}}(\Sigma_0)
\|\Delta\Sigma_0\|_F\\
&\leq
C_{\mathrm{chol}}(\Sigma_0)
r_{\mathrm{IFT}}(\Sigma_0)\\
&=
r_\theta(L). 
\end{align*}
Therefore, the condition
\eqref{eqn: IFT cholesky perturbation radius} of
Lemma~\ref{lemma: local stability of dual optimizer under cholesky perturbation}
holds. Applying
\eqref{eqn: dual optimizer locally Lipschitz in cholesky factor} and then
\eqref{eqn: marginal covariance perturbation implies cholesky perturbation}
gives
\begin{align*}
\numberthis\label{eqn: marginal covariance perturbation implies dual optimizer perturbation}\left\|
\theta^*(L+\Delta L)-\theta^*(L)
\right\|_\Theta
&\leq
C_\theta(L)\|\Delta L\|_F\\
&\leq
C_\theta(L)C_{\mathrm{chol}}(\Sigma_0)
\|\Delta\Sigma_0\|_F. 
\end{align*}

Combining \eqref{eqn: radius for marginal cholesky perturbation} and
\eqref{eqn: covariance radius for applying IFT stability lemma}, we now choose
\begin{equation}
r_\Sigma
:=
\min
\left\{
r_{\mathrm{chol}}(\Sigma_0),
r_{\mathrm{IFT}}(\Sigma_0)
\right\}
>0.
\label{eqn: choice of marginal covariance perturbation radius}
\end{equation}
Under the condition $\|\Delta\Sigma_0\|_F\leq r_\Sigma$,
the perturbed marginal covariance matrices
$\Sigma^\mu+\Delta\Sigma^\mu$ and
$\Sigma^\nu+\Delta\Sigma^\nu$ remain positive definite.
Moreover,
\eqref{eqn: marginal covariance perturbation implies cholesky perturbation},
\eqref{eqn: marginal covariance perturbation implies dual optimizer perturbation},
and
\eqref{eqn: perturbed optimal Z uniform positive lower bound}
all hold.
\item
We next control the perturbation of \(Q\).
We obtain
\begin{equation}
\left\|
Q(L+\Delta L)-Q(L)
\right\|_F
\leq
C_{\Delta Q}
\|\Delta\Sigma_0\|_F,
\label{eqn: Q perturbation controlled by marginal covariance perturbation}
\end{equation}
where
\begin{equation}
C_{\Delta Q}
:=
\left(
2\|L\|_2+r_\theta(L)
\right)
\|(\Sigma^\gamma)^{-1}\|_2
C_{\mathrm{chol}}(\Sigma_0).
\label{eqn: constant for Q perturbation}
\end{equation}
To see why, by \eqref{eqn: Q depending on cholesky matrix L},
\begin{align*}
Q(L+\Delta L)-Q(L)
&=
(L+\Delta L)^\top
(\Sigma^\gamma)^{-1}
(L+\Delta L)
-
L^\top(\Sigma^\gamma)^{-1}L\\
&=
(\Delta L)^\top
(\Sigma^\gamma)^{-1}
(L+\Delta L)
+
L^\top(\Sigma^\gamma)^{-1}\Delta L.
\end{align*}
Therefore,
\begin{align*}
&\left\|
Q(L+\Delta L)-Q(L)
\right\|_F\\
&\leq
\left\|
(\Delta L)^\top
(\Sigma^\gamma)^{-1}
(L+\Delta L)
\right\|_F +
\left\|
L^\top
(\Sigma^\gamma)^{-1}
\Delta L
\right\|_F\\
&\leq
\|\Delta L\|_F
\|(\Sigma^\gamma)^{-1}\|_2
\|L+\Delta L\|_2 +
\|L\|_2
\|(\Sigma^\gamma)^{-1}\|_2
\|\Delta L\|_F\\
&\leq \|\Delta L\|_F
\|(\Sigma^\gamma)^{-1}\|_2 (\|L\|_2+\|\Delta L\|_F)
 +
\|L\|_2
\|(\Sigma^\gamma)^{-1}\|_2
\|\Delta L\|_F
\\
&\leq \|\Delta L\|_F
\|(\Sigma^\gamma)^{-1}\|_2 (\|L\|_2+r_\theta(L))
 +
\|L\|_2
\|(\Sigma^\gamma)^{-1}\|_2
\|\Delta L\|_F
\\
&\leq
\left(
2\|L\|_2+r_\theta(L)
\right)
\|(\Sigma^\gamma)^{-1}\|_2
\|\Delta L\|_F\\
&\leq
\left(
2\|L\|_2+r_\theta(L)
\right)
\|(\Sigma^\gamma)^{-1}\|_2
C_{\mathrm{chol}}(\Sigma_0)
\|\Delta\Sigma_0\|_F\\
&=
C_{\Delta Q}
\|\Delta\Sigma_0\|_F.
\end{align*}
Here, the second inequality follows from
\(
\|ABC\|_F
\leq
\|A\|_F\|B\|_2\|C\|_2
\), and the fourth inequality follows from
\eqref{eqn: delta L F norm is bouned by r theta (L) proved from assumption in staiblity thoerem}.
The sixth inequality follows from
\eqref{eqn: marginal covariance perturbation implies cholesky perturbation}.
This proves
\eqref{eqn: Q perturbation controlled by marginal covariance perturbation}.
\item
We next control the perturbation of the optimal \(Z\).
We first obtain, for every \(h=(h_1,h_2)\in\Theta_1\times\Theta_2=\Theta\),
\begin{equation}
\|A(h)\|_F
\leq
2\|h\|_\Theta.
\label{eqn: A Frobenius upper bound for stability}
\end{equation}
This is because
\(
A(h)=A(h_1,0)+A(0,h_2)
\) by the linearity of \(A\).
Moreover, by \eqref{eqn: A from theta}, the nonzero entries of
\(A(h_1,0)\) and \(A(0,h_2)\) lie in disjoint matrix blocks. Hence,
\begin{align*}
\|A(h)\|_F^2
&=
\|A(h_1,0)\|_F^2
+
\|A(0,h_2)\|_F^2\\
&\leq
4\|h_1\|_{\Theta_1}^2
+
4\|h_2\|_{\Theta_2}^2\\
&=
4\|h\|_\Theta^2,
\end{align*}
where the inequality follows from
\eqref{eqn: A first block upper bound} and
\eqref{eqn: A second block upper bound}.
Taking square roots gives
\eqref{eqn: A Frobenius upper bound for stability}.

We now obtain
\begin{equation}
\left\|
Z\left(
\theta^*(L+\Delta L);L+\Delta L
\right)
-
Z\left(
\theta^*(L);L
\right)
\right\|_F
\leq
C_{\Delta Z}
\|\Delta\Sigma_0\|_F,
\label{eqn: optimal Z perturbation controlled by marginal covariance perturbation}
\end{equation}
where
\begin{equation}
C_{\Delta Z}
:=
C_{\Delta Q}
+
2C_\theta(L)C_{\mathrm{chol}}(\Sigma_0).
\label{eqn: constant for optimal Z perturbation}
\end{equation}
By \eqref{eqn: Z depending on theta and cholesky matrix L} and the
linearity of \(A\),
\begin{align*}
&
Z\left(
\theta^*(L+\Delta L);L+\Delta L
\right)
-
Z\left(
\theta^*(L);L
\right)\\
&=
Q(L+\Delta L)-Q(L)
+
A\left(
\theta^*(L+\Delta L)-\theta^*(L)
\right).
\end{align*}
Therefore,
\begin{align*}
&
\left\|
Z\left(
\theta^*(L+\Delta L);L+\Delta L
\right)
-
Z\left(
\theta^*(L);L
\right)
\right\|_F\\
&\leq
\left\|
Q(L+\Delta L)-Q(L)
\right\|_F
+
\left\|
A\left(
\theta^*(L+\Delta L)-\theta^*(L)
\right)
\right\|_F\\
&\leq
C_{\Delta Q}\|\Delta\Sigma_0\|_F
+
2
\left\|
\theta^*(L+\Delta L)-\theta^*(L)
\right\|_\Theta\\
&\leq
C_{\Delta Q}\|\Delta\Sigma_0\|_F
+
2C_\theta(L)C_{\mathrm{chol}}(\Sigma_0)
\|\Delta\Sigma_0\|_F\\
&=
C_{\Delta Z}
\|\Delta\Sigma_0\|_F.
\end{align*}
Here, the second inequality follows from
\eqref{eqn: Q perturbation controlled by marginal covariance perturbation}
and
\eqref{eqn: A Frobenius upper bound for stability}, and the third
inequality follows from
\eqref{eqn: marginal covariance perturbation implies dual optimizer perturbation}.
This proves
\eqref{eqn: optimal Z perturbation controlled by marginal covariance perturbation}.

\item
We next control the perturbation of the optimal covariance in the
Cholesky coordinates.
Denote the optimal covariance matrices of Problem \eqref{eqn:problem-over-cov} in the Cholesky coordinates corresponding to \(L\) and \(L+\Delta L\) by \(\widetilde\Sigma^*\) and \(\widetilde\Sigma_\Delta^*\), respectively.
We obtain
\begin{equation}
\left\|
\widetilde\Sigma_\Delta^*
-
\widetilde\Sigma^*
\right\|_F
\leq
C_{\Delta\widetilde\Sigma}
\|\Delta\Sigma_0\|_F,
\label{eqn: optimal covariance perturbation in cholesky coordinates}
\end{equation}
where
\begin{equation}
C_{\Delta\widetilde\Sigma}
:=
\frac{
2C_{\Delta Z}
}{
m_{Z^*}(L)^2
}.
\label{eqn: constant for optimal covariance perturbation in cholesky coordinates}
\end{equation}

By Theorem~\ref{theorem: dual problem of ecot gaussian}, the primal and
dual optimizers satisfy
\begin{equation*}
\widetilde\Sigma^*
=
Z(\theta^*(L);L)^{-1},
\qquad
\widetilde\Sigma_\Delta^*
=
Z(\theta^*(L+\Delta L);L+\Delta L)^{-1}.
\end{equation*}

By \eqref{eqn: m Z star depending on L},
\begin{align*}
\numberthis\label{eqn: operator norm of optimal covariance in cholesky coordinates}\|\widetilde\Sigma^*\|_2
&=
\left\|
Z(\theta^*(L);L)^{-1}
\right\|_2\\
&=
\frac{1}{
\lambda_{\min}
\left(
Z(\theta^*(L);L)
\right)
}\\
&=
\frac{1}{m_{Z^*}(L)}.
\end{align*}

Moreover, by
\eqref{eqn: perturbed optimal Z uniform positive lower bound},
\begin{align*}
\numberthis\label{eqn: operator norm of perturbed optimal covariance in cholesky coordinates} \|\widetilde\Sigma_\Delta^*\|_2
&=
\left\|
Z(\theta^*(L+\Delta L);L+\Delta L)^{-1}
\right\|_2\\
&=
\frac{1}{
\lambda_{\min}
\left(
Z(\theta^*(L+\Delta L);L+\Delta L)
\right)
}\\
&\leq
\frac{2}{m_{Z^*}(L)}. 
\end{align*}

Using the identity
\(
X^{-1}-Y^{-1}
=
X^{-1}(Y-X)Y^{-1}
\),
with
\(
X=
Z(\theta^*(L+\Delta L);L+\Delta L)
\)
and
\(
Y=
Z(\theta^*(L);L)
\),
we have
\begin{align*}
&
\left\|
\widetilde\Sigma_\Delta^*
-
\widetilde\Sigma^*
\right\|_F\\
&=
\left\|
Z(\theta^*(L+\Delta L);L+\Delta L)^{-1}
-
Z(\theta^*(L);L)^{-1}
\right\|_F\\
&=
\|
Z(\theta^*(L+\Delta L);L+\Delta L)^{-1}
\left(
Z(\theta^*(L);L)
-
Z(\theta^*(L+\Delta L);L+\Delta L)
\right) \\&
Z(\theta^*(L);L)^{-1}
\|_F\\
&\leq
\|
Z(\theta^*(L+\Delta L);L+\Delta L)^{-1}
\|_2
\|
Z(\theta^*(L+\Delta L);L+\Delta L)
-
Z(\theta^*(L);L)
\|_F \\&
\|
Z(\theta^*(L);L)^{-1}
\|_2\\
&\leq
\frac{2}{m_{Z^*}(L)}
C_{\Delta Z}
\|\Delta\Sigma_0\|_F
\frac{1}{m_{Z^*}(L)}\\
&=
C_{\Delta\widetilde\Sigma}
\|\Delta\Sigma_0\|_F.
\end{align*}
Here, the first inequality follows from
\(
\|ABC\|_F
\leq
\|A\|_2\|B\|_F\|C\|_2
\),
and the second inequality follows from
\eqref{eqn: operator norm of optimal covariance in cholesky coordinates},
\eqref{eqn: operator norm of perturbed optimal covariance in cholesky coordinates},
and
\eqref{eqn: optimal Z perturbation controlled by marginal covariance perturbation}.
This proves
\eqref{eqn: optimal covariance perturbation in cholesky coordinates}.

\item
We finally return to the original coordinates.
We obtain 
\begin{equation*}
\left\|
\Sigma^{\pi_\Delta^*}
-
\Sigma^{\pi^*}
\right\|_F
\leq
C_\Sigma
\|\Delta\Sigma_0\|_F,
\end{equation*}
where
\begin{equation}
C_\Sigma
:=
\frac{
\left(
3\|L\|_2+2r_\theta(L)
\right)
C_{\mathrm{chol}}(\Sigma_0)
}{
m_{Z^*}(L)
}
+
\|L\|_2
\left(
\|L\|_2+r_\theta(L)
\right)
C_{\Delta\widetilde\Sigma}.
\label{eqn: deterministic marginal covariance stability constant}
\end{equation}

By Lemma \ref{lemma: pi start is minimizer of original problem if and ony if its pushfowrad by T chol is optimizer of the problem in chol coordinate},
\begin{equation}
\Sigma^{\pi^*}
=
L\widetilde\Sigma^*L^\top,
\qquad
\Sigma^{\pi_\Delta^*}
=
(L+\Delta L)
\widetilde\Sigma_\Delta^*
(L+\Delta L)^\top.
\label{eqn: optimal covariance recovered from cholesky coordinates}
\end{equation}
Therefore,
\begin{align*}
\Sigma^{\pi_\Delta^*}
-
\Sigma^{\pi^*}
&=
(L+\Delta L)
\widetilde\Sigma_\Delta^*
(L+\Delta L)^\top
-
L\widetilde\Sigma^*L^\top\\
&=
\Delta L\,
\widetilde\Sigma_\Delta^*
(L+\Delta L)^\top +
L\left(
\widetilde\Sigma_\Delta^*
-
\widetilde\Sigma^*
\right)
(L+\Delta L)^\top +
L\widetilde\Sigma^*
(\Delta L)^\top.
\end{align*}
Hence,
\begin{align*}
&\left\|
\Sigma^{\pi_\Delta^*}
-
\Sigma^{\pi^*}
\right\|_F\\
&\leq
\|\Delta L\|_F
\|\widetilde\Sigma_\Delta^*\|_2
\|L+\Delta L\|_2 +
\|L\|_2
\left\|
\widetilde\Sigma_\Delta^*
-
\widetilde\Sigma^*
\right\|_F
\|L+\Delta L\|_2 \\&\quad +
\|L\|_2
\|\widetilde\Sigma^*\|_2
\|\Delta L\|_F\\
&\leq
\|\Delta L\|_F
\|\widetilde\Sigma_\Delta^*\|_2
(\|L\|_2+\|\Delta L\|_F) +
\|L\|_2
\left\|
\widetilde\Sigma_\Delta^*
-
\widetilde\Sigma^*
\right\|_F
(\|L\|_2+\|\Delta L\|_F)  \\&\quad+
\|L\|_2
\|\widetilde\Sigma^*\|_2
\|\Delta L\|_F\\
&\leq
\|\Delta L\|_F
\|\widetilde\Sigma_\Delta^*\|_2
(\|L\|_2+r_\theta(L)) +
\|L\|_2
\left\|
\widetilde\Sigma_\Delta^*
-
\widetilde\Sigma^*
\right\|_F
(\|L\|_2+r_\theta(L))  \\&\quad+
\|L\|_2
\|\widetilde\Sigma^*\|_2
\|\Delta L\|_F\\
&\leq
C_{\mathrm{chol}}(\Sigma_0)
\|\Delta\Sigma_0\|_F
\frac{2}{m_{Z^*}(L)}
\left(
\|L\|_2+r_\theta(L)
\right)\\
&\quad+
\|L\|_2
C_{\Delta\widetilde\Sigma}
\|\Delta\Sigma_0\|_F
\left(
\|L\|_2+r_\theta(L)
\right)\\
&\quad+
\|L\|_2
\frac{1}{m_{Z^*}(L)}
C_{\mathrm{chol}}(\Sigma_0)
\|\Delta\Sigma_0\|_F\\
&=
C_\Sigma
\|\Delta\Sigma_0\|_F,
\end{align*}
where  the third inequality follows from
\eqref{eqn: delta L F norm is bouned by r theta (L) proved from assumption in staiblity thoerem}, and the fourth inequality follows from
\eqref{eqn: marginal covariance perturbation implies cholesky perturbation},
\eqref{eqn: operator norm of optimal covariance in cholesky coordinates},
\eqref{eqn: operator norm of perturbed optimal covariance in cholesky coordinates},
and
\eqref{eqn: optimal covariance perturbation in cholesky coordinates}.
This proves \eqref{eqn: stability of optimal covariance under marginal perturbation}.
\end{enumerate}

It remains only to verify the dependence of the constants
\(r_\Sigma\) and \(C_\Sigma\) stated in the theorem.
The constants
\(r_{\mathrm{chol}}(\Sigma_0)\) and
\(C_{\mathrm{chol}}(\Sigma_0)\) depend only on \(\Sigma_0\), while
\(r_\theta(L)\) and \(C_\theta(L)\) have the dependence stated in
Lemma~\ref{lemma: local stability of dual optimizer under cholesky perturbation}.
Moreover, $m_{Z^*}(L)$
depends only on \(L,\Sigma^\gamma,T,\) and \(d\).
Since \(L\) is determined by \(\Sigma_0\), it follows from
\eqref{eqn: covariance radius for applying IFT stability lemma} and
\eqref{eqn: choice of marginal covariance perturbation radius} that
\(r_\Sigma\) depends only on
\(\Sigma_0,\Sigma^\gamma,T,\) and \(d\).
Similarly, by
\eqref{eqn: constant for Q perturbation},
\eqref{eqn: constant for optimal Z perturbation},
\eqref{eqn: constant for optimal covariance perturbation in cholesky coordinates},
and
\eqref{eqn: deterministic marginal covariance stability constant},
\(C_\Sigma\) depends only on
\(\Sigma_0,\Sigma^\gamma,T,\) and \(d\).
This proves the result.

\end{proof}

\begin{lemma}
\label{lemma:concentration inequality}
Let
\(\{\mathbf u^i\}_{i=1}^N\subset\mathbb R^{Td}\) and
\(\{\mathbf y^i\}_{i=1}^N\subset\mathbb R^{Td}\)
be i.i.d. samples drawn from
\begin{equation*}
\mu=\mathcal N(0,\Sigma^\mu),
\qquad
\nu=\mathcal N(0,\Sigma^\nu),
\end{equation*}
respectively, and let
\begin{equation*}
\hat\Sigma^\mu
=
\frac1N\sum_{i=1}^N
\mathbf u^i(\mathbf u^i)^\top,
\qquad
\hat\Sigma^\nu
=
\frac1N\sum_{i=1}^N
\mathbf y^i(\mathbf y^i)^\top.
\end{equation*}
There exist  a constant \(C>0\) such that, for every
\(\epsilon\geq0\), with probability at least
\(1-2e^{-\epsilon}\),
\begin{equation}\label{eqn: for mu Sigma marignal hat minus Sigma marignal isbouned by Sigma maringal and fucntion of sample size }
\|\hat\Sigma^\mu-\Sigma^\mu\|_2
\leq
C \|\Sigma^\mu\|_2
\left(
\sqrt{\frac{Td+\epsilon}{N}}
+
\frac{Td+\epsilon}{N}
\right),
\end{equation}
and, with probability at least
\(1-2e^{-\epsilon}\),
\begin{equation}\label{eqn: for nu Sigma marignal hat minus Sigma marignal isbouned by Sigma maringal and fucntion of sample size }
\|\hat\Sigma^\nu-\Sigma^\nu\|_2
\leq
C \|\Sigma^\nu\|_2
\left(
\sqrt{\frac{Td+\epsilon}{N}}
+
\frac{Td+\epsilon}{N}
\right).
\end{equation}
\end{lemma}
\begin{proof}[Proof of Lemma~\ref{lemma:concentration inequality}]
By \cite[Chapter~2]{Vershynin_2026}, for every
\(x\in\mathbb R^{Td}\),
\begin{equation*}
    \|\langle \mathbf u^i,x\rangle\|_{\psi_2}
=
\sqrt{\frac{8}{3}}\,
\|\langle \mathbf u^i,x\rangle\|_{L^2},
\end{equation*}
where \(\|\cdot\|_{\psi_2}\) denotes the subgaussian norm.
Hence, condition~(4.29) of \cite[Theorem~4.7.1]{Vershynin_2026} holds with
\(K=\sqrt{8/3}\). Applying \cite[Theorem~4.7.1]{Vershynin_2026} and \cite[Remark~4.7.3]{Vershynin_2026}, we obtain
\eqref{eqn: for mu Sigma marignal hat minus Sigma marignal isbouned by Sigma maringal and fucntion of sample size }
with probability at least \(1-2e^{-\epsilon}\).
The same argument applied to \(\mathbf y^i\) gives
\eqref{eqn: for nu Sigma marignal hat minus Sigma marignal isbouned by Sigma maringal and fucntion of sample size }.
\end{proof}

\begin{proof}[Proof of Theorem \ref{theorem: error bound}]
Let \(C>0\) satisfy the bounds in
Lemma~\ref{lemma:concentration inequality}.
We first show that
\begin{equation}
\|\Delta\Sigma_0\|_F
\leq
\psi(N,\epsilon,C) \label{eqn: delta Sigma 0 F norm bouned by psi functiono}
\end{equation}
with probability at least \(1-4e^{-\epsilon}\).
By
\eqref{eqn: for mu Sigma marignal hat minus Sigma marignal isbouned by Sigma maringal and fucntion of sample size },
the probability that the bound for
\(\|\hat\Sigma^\mu-\Sigma^\mu\|_2\) fails is at most
\(2e^{-\epsilon}\).
Similarly, by
\eqref{eqn: for nu Sigma marignal hat minus Sigma marignal isbouned by Sigma maringal and fucntion of sample size },
the probability that the bound for
\(\|\hat\Sigma^\nu-\Sigma^\nu\|_2\) fails is at most
\(2e^{-\epsilon}\).
Therefore, the probability that at least one of these two bounds fails is at most
\(4e^{-\epsilon}\).
Hence, with probability at least \(1-4e^{-\epsilon}\),
both
\eqref{eqn: for mu Sigma marignal hat minus Sigma marignal isbouned by Sigma maringal and fucntion of sample size }
and
\eqref{eqn: for nu Sigma marignal hat minus Sigma marignal isbouned by Sigma maringal and fucntion of sample size }
hold simultaneously.

On this event, we have
\begin{align*}
\|\Delta\Sigma_0\|_F^2
&=
\|\hat\Sigma^\mu-\Sigma^\mu\|_F^2
+
\|\hat\Sigma^\nu-\Sigma^\nu\|_F^2\\
&\leq
Td\|\hat\Sigma^\mu-\Sigma^\mu\|_2^2
+
Td\|\hat\Sigma^\nu-\Sigma^\nu\|_2^2\\
&\leq
Td\,C^2
\left(
\|\Sigma^\mu\|_2^2
+
\|\Sigma^\nu\|_2^2
\right)
\left(
\sqrt{\frac{Td+\epsilon}{N}}
+
\frac{Td+\epsilon}{N}
\right)^2\\
&=
\psi(N,\epsilon,C)^2.
\end{align*}
Here, the first inequality follows from
\(\|A\|_F\leq\sqrt{Td}\|A\|_2\) for every
\(A\in\mathbb R^{Td\times Td}\).
The second inequality follows from
\eqref{eqn: for mu Sigma marignal hat minus Sigma marignal isbouned by Sigma maringal and fucntion of sample size }
and
\eqref{eqn: for nu Sigma marignal hat minus Sigma marignal isbouned by Sigma maringal and fucntion of sample size }. We also have \(\psi(N,\epsilon,C)\geq0\), since every factor in its definition is
nonnegative.
Taking square roots gives \eqref{eqn: delta Sigma 0 F norm bouned by psi functiono}.

Under the sample-size condition
\eqref{eqn: sample size condition for statistical stability},
on the event where
\eqref{eqn: delta Sigma 0 F norm bouned by psi functiono} holds,
\begin{equation*}
\|\Delta\Sigma_0\|_F
\leq
\psi(N,\epsilon,C)
\leq
r_\Sigma.
\end{equation*}
This event has probability at least \(1-4e^{-\epsilon}\).
Hence, on this event,
Theorem~\ref{theorem: stability under marginal covariance perturbations}
applies.
In this case, the perturbed marginal covariance matrices are
\(\hat\Sigma^\mu\) and \(\hat\Sigma^\nu\), so the perturbed optimizer
\(\pi_\Delta^*\) in
Theorem~\ref{theorem: stability under marginal covariance perturbations}
is precisely \(\hat\pi^*\).
Therefore, by
\eqref{eqn: stability of optimal covariance under marginal perturbation},
\begin{align*}
\|\Sigma^{\hat\pi^*}-\Sigma^{\pi^*}\|_F
&\leq
C_\Sigma\|\Delta\Sigma_0\|_F\\
&\leq
C_\Sigma\psi(N,\epsilon,C).
\end{align*}
Thus,
\eqref{eqn: statistical stability of optimal covariance}
holds with probability at least \(1-4e^{-\epsilon}\).
\end{proof}
\section{Supplementary Analysis for Section~\ref{sec: sinkhorn equivalence}}\label{sec: supp sinkhorn equivalence}

\begin{proof}[Proof of Lemma \ref{lemma: gaussianity preserved under sinkhorn projections}]
We first consider the causal projection
\begin{equation*}
\min_{\pi\in\Pi_c(\mu,*)}D(\pi,\rho).
\end{equation*}
By Lemma~6.2 of
\cite{eckstein2023computationalmethodsadaptedoptimal}, its optimizer is
constructed backward in time, with
\begin{equation*}
\pi(y_t\mid u_{1:t},y_{1:t-1})
\propto
\rho(y_t\mid u_{1:t},y_{1:t-1})e^{-f_t},
\end{equation*}
where \(f_T=0\), and, for \(t<T\), \(f_t\) is the relative entropy
between the conditional distributions of the future variables under
\(\pi\) and \(\rho\), conditional on
\((u_{1:t},y_{1:t})\) (see \cite[(23)]{eckstein2023computationalmethodsadaptedoptimal}.

Since \(\mu\) and \(\rho\) are Gaussian, their conditional
distributions are Gaussian. At \(t=T\), \(f_T=0\), so the conditional
distribution constructed at the terminal step is Gaussian.
Proceeding backward, if the conditional distribution constructed at
time \(t+1\) is Gaussian, then \(f_t\), being the relative entropy
between conditional Gaussian distributions, is quadratic in the
conditioning variables. Hence,
\(\rho(y_t\mid u_{1:t},y_{1:t-1})e^{-f_t}\) is the exponential of a
quadratic function and, after normalization, defines a Gaussian
conditional distribution. Applying this argument recursively shows
that the optimizer of the causal projection is Gaussian.

We next consider the projection
\begin{equation*}
\min_{\pi\in\Pi(*,\nu)}D(\pi,\rho).
\end{equation*}
 Section~6 of
\cite{nutz2021introductionentropic} shows that the optimizer is given by
\begin{equation*}
\pi(du,dy)=\rho(du\mid y)\nu(dy).
\end{equation*}
Since \(\rho\) and \(\nu\) are Gaussian,
\(\rho(du\mid y)\) is Gaussian with mean affine in \(y\) and covariance
independent of \(y\). Hence, the resulting joint distribution \(\pi\) is
Gaussian.

It remains to show that the two Gaussian optimizers are non--degenerate. We first show that both projection problems have finite
value. This is true because $P$
belongs to both \(\Pi_c(\mu,*)\) and \(\Pi(*,\nu)\).
$P$ is non--degenerate Gaussian , and hence have finite
relative entropy with respect to the non--degenerate Gaussian
distribution \(\rho\). Therefore, both projection problems have finite
value, and their optimizers \(\pi\) satisfy
\(D(\pi,\rho)<\infty\), which implies \(\pi\ll\rho\).

If such an optimizer \(\pi\) were degenerate, it would assign probability
one to an affine hyperplane. Since the non--degenerate Gaussian
distribution \(\rho\) assigns probability zero to every affine
hyperplane, this would contradict \(\pi\ll\rho\).
Hence, both optimizers are non--degenerate.
\end{proof}

\begin{proof}[Proof of Proposition~\ref{proposition: sinkhorn is bcd}]
We prove the two identities simultaneously by induction on \(k\).
At \(k=0\), we have
\begin{equation}
\left(\Sigma^{\pi^{\{0\}}}\right)^{-1}
=
Z(\theta^{\{0\}}). \label{eqn: equivalcne between intilization of bcd and sinkhorn}
\end{equation}
To see this, \(\theta^{\{0\}}=\mathbf 0\) and
\eqref{eqn: Z from theta} give
\begin{equation*}
    Z(\theta^{\{0\}})
=
(\Sigma^\gamma)^{-1}+I_{2Td}.
\end{equation*}
Meanwhile, 
the density of \(\pi^{\{0\}}\) with respect to the Lebesgue
measure is proportional to
\begin{equation*}
    \exp\left(
-\frac12
\mathbf z^\top
\left(
(\Sigma^\gamma)^{-1}+I_{2Td}
\right)
\mathbf z
\right),\qquad
\mathbf z=
\begin{bmatrix}
\mathbf u\\
\mathbf y
\end{bmatrix}
\end{equation*}
because 
\eqref{eqn: init of sinkhorn} gives
\begin{equation*}
    \frac{d\pi^{\{0\}}}{dP}(\mathbf z)
\propto
\exp\left(
-\frac12
\mathbf z^\top(\Sigma^\gamma)^{-1}\mathbf z
\right),
\end{equation*}
and \(P=\mu\otimes\nu=\mathcal N(0,I_{2Td})\) in the Cholesky coordinates.
Hence,
\(\pi^{\{0\}}\) is centered and has precision matrix
\((\Sigma^\gamma)^{-1}+I_{2Td}=Z(\theta^{\{0\}})\), which proves
\eqref{eqn: equivalcne between intilization of bcd and sinkhorn}.

Now suppose that, for some \(k\geq0\), \(\pi^{\{k\}}\) is centered and
\begin{equation*}
\left(\Sigma^{\pi^{\{k\}}}\right)^{-1}
=
Z(\theta^{\{k\}}).
\end{equation*}
The iterate \(\pi^{\{k\}}\) is a non--degenerate Gaussian distribution,
because \(\pi^{\{0\}}\) is non--degenerate Gaussian and
Lemma~\ref{lemma: gaussianity preserved under sinkhorn projections}
shows recursively that both Sinkhorn projection steps preserve this
property.

We first show that the odd Sinkhorn projection
\eqref{eqn: odd sinkhorn} corresponds to the first BCD update
\eqref{eq:update-upsi}.
Define a quadratic function \(g^{\{k\}}\), up to an additive constant, by
\begin{equation*}
    \frac{d\pi^{\{k\}}}{dP}
\propto
\exp\left(g^{\{k\}}-c\right).
\end{equation*}
The function \(g^{\{k\}}\) has no linear term, because both
\(\pi^{\{k\}}\) and \(P\) are centered Gaussian distributions and
\(c\) has no linear term. Therefore, we can write
\begin{equation*}
    g^{\{k\}}(\mathbf z)
=
-\frac12
\mathbf z^\top Q_g^{\{k\}}\mathbf z
+
\mathrm{Const}.
\end{equation*}
The odd Sinkhorn projection \eqref{eqn: odd sinkhorn}  can be written as
\begin{equation}
\min_{\pi\in\Pi_c(\mu,*)}
\left\{
\int(c-g^{\{k\}})\,d\pi
+
D(\pi,P)
\right\},
\label{eqn: odd sinkhorn rewritten as entropy regularized problem}
\end{equation}
up to an additive constant independent of \(\pi\).
This follows from
\eqref{eqn: odd sinkhorn} and the definition of \(g^{\{k\}}\), because
for every feasible Gaussian coupling,
\begin{align*}
    D(\pi,\pi^{\{k\}})
&=
\int
\log\frac{d\pi}{d\pi^{\{k\}}}\,d\pi\\
&=
D(\pi,P)
-
\int
\log\frac{d\pi^{\{k\}}}{dP}\,d\pi\\
&=
D(\pi,P)
+
\int(c-g^{\{k\}})\,d\pi
+
\mathrm{const.}
\end{align*}
Restricting to Gaussian couplings does not change the optimizer,
since Lemma~\ref{lemma: gaussianity preserved under sinkhorn projections}
shows that the optimizer of \eqref{eqn: odd sinkhorn} is Gaussian.
The optimizer of \eqref{eqn: odd sinkhorn rewritten as entropy regularized problem}
is centered. To see this, let
\(\pi=\mathcal N(m,\Sigma)\) be a feasible Gaussian distribution.
Since \(g^{\{k\}}\) has no linear term, the standard quadratic--expectation
calculation and the formula for the relative entropy between
\(\mathcal N(m,\Sigma)\) and \(P=\mathcal N(0,I_{2Td})\) show that the
terms in \eqref{eqn: odd sinkhorn rewritten as entropy regularized problem}
depending on the mean \(m\) are
\begin{equation*}
\frac12
m^\top
\left(
Q_g^{\{k\}}+(\Sigma^\gamma)^{-1}+I_{2Td}
\right)m
=
\frac12
m^\top
\left(
Q_g^{\{k\}}+Q
\right)m.
\end{equation*}
Moreover, \(Q_g^{\{k\}}+Q
=
\left(\Sigma^{\pi^{\{k\}}}\right)^{-1}
\succ0\), by the definition of $g^{\{k\}}.$ 
Therefore, the mean-dependent term is uniquely minimized at
\(m=\mathbf 0\). Since \(m=\mathbf 0\) is feasible, the optimizer is centered. Following that, for a Gaussian distribution \(\pi=\mathcal N(0,\Sigma)\), the
covariance-dependent part of
\(\int(c-g^{\{k\}})\,d\pi\) is
\begin{equation*}
\frac12
\operatorname{Tr}\left(
\left(
(\Sigma^\gamma)^{-1}+Q_g^{\{k\}}
\right)\Sigma
\right),
\end{equation*}
using the standard identity for the expectation of a quadratic form.
Moreover, the covariance  dependent part of \(D(\pi,P)\) is
\begin{equation*}
\frac12\operatorname{Tr}(\Sigma)
-\frac12\log\det\Sigma,
\end{equation*}
by the formula for the KL divergence between centered Gaussian distributions.
Therefore, the covariance-dependent part of
\eqref{eqn: odd sinkhorn rewritten as entropy regularized problem} is
\begin{equation*}
\frac12
\operatorname{Tr}
\left(
(Q_g^{\{k\}}+Q)\Sigma
\right)
-
\frac12\log\det\Sigma.
\end{equation*}
The corresponding covariance constraints in \eqref{eqn: odd sinkhorn rewritten as entropy regularized problem} are
\begin{equation*}
    P_u\Sigma P_u^\top=I_{Td},
\qquad
S_tP_y\Sigma P_u^\top E_{t+1}=0,
\quad
t=1,\ldots,T-1,
\end{equation*}
Therefore, introducing the dual variables
\(\widetilde\Lambda_u\in \mathbb S^{Td}
\) and
\(\widetilde\Psi=\{\widetilde\Psi_t\}_{t=1}^{T-1}\in \prod_{t=1}^{T-1}\mathbb R^{td\times d}\) and following a similar calculation as in \eqref{eqn: lagrangian of cov optimizatio problem}, the Lagrangian of \eqref{eqn: odd sinkhorn rewritten as entropy regularized problem} has the form
\begin{equation*}
    \widetilde{\mathcal L}
(\Sigma;\widetilde\Lambda_u,\widetilde\Psi)
=
\frac12\operatorname{Tr}(\widetilde Z\Sigma)
-
\frac12\log\det\Sigma
-
\operatorname{Tr}(\widetilde\Lambda_u),
\end{equation*}
where \(\widetilde Z
=
Q_g^{\{k\}}
+
Q
+
2P_u^\top\widetilde\Lambda_uP_u
+
\begin{bmatrix}
0&\Gamma_{\widetilde\Psi}\\
\Gamma_{\widetilde\Psi}^\top&0
\end{bmatrix}.\)
Minimizing the Lagrangian over \(\Sigma\succ0\), following the same calculation as in the derivation of \eqref{eqn: dual objective}, gives
the unique minimizer
\begin{equation}
\Sigma=\widetilde Z^{-1},
\label{eqn: primal dual relation for odd sinkhorn projection}
\end{equation}
and, up to an additive constant independent of
\((\widetilde\Lambda_u,\widetilde\Psi)\), the dual problem of \eqref{eqn: odd sinkhorn rewritten as entropy regularized problem}
 is
\begin{equation}
\max_{\widetilde\Lambda_u,\widetilde\Psi}
\widetilde{\mathcal D}
(\widetilde\Lambda_u,\widetilde\Psi),
\label{eqn: dual problem of EOT version of the sinkhorn step}
\end{equation}
where
\begin{equation}
\widetilde{\mathcal D}
(\widetilde\Lambda_u,\widetilde\Psi)
=
\begin{cases}
\displaystyle
\frac12\log\det\widetilde Z
-\operatorname{Tr}(\widetilde\Lambda_u),
& \widetilde Z\succ0,\\[1ex]
-\infty,
& \text{o.w.}.
\end{cases} \label{eqn: def of tilde D on tilde lambda_u and tilde psi, the dual of the one maringal causal projection}
\end{equation}

We next show that
\eqref{eqn: dual problem of EOT version of the sinkhorn step}
is equivalent to the first BCD update
\eqref{eq:update-upsi}.
The quadratic coefficient \(Q_g^{\{k\}}\) satisfies
\begin{equation*}
Q_g^{\{k\}}
=
Z(\theta^{\{k\}})-Q.
\end{equation*}
To see this, the induction hypothesis gives
\begin{equation*}
\log\frac{d\pi^{\{k\}}}{dP}(\mathbf z)
=
-\frac12
\mathbf z^\top
\left(
Z(\theta^{\{k\}})-I_{2Td}
\right)
\mathbf z
+\mathrm{Const}.
\end{equation*}

Meanwhile, the definition of \(g^{\{k\}}\)  gives
\begin{equation*}
\log\frac{d\pi^{\{k\}}}{dP}(\mathbf z)
=
-\frac12
\mathbf z^\top
\left(
Q_g^{\{k\}}+(\Sigma^\gamma)^{-1}
\right)
\mathbf z
+\mathrm{Const}.
\end{equation*}
Therefore, comparing the quadratic coefficients and using
\eqref{eqn:primal loss} gives
\begin{equation*}
Q_g^{\{k\}}
=
Z(\theta^{\{k\}})-Q.
\end{equation*}
Hence, substituting this expression into \(\widetilde Z\) gives
\begin{equation*}
\begin{aligned}
\widetilde Z
&=
Q
+
2P_u^\top
\left(
\widetilde\Lambda_u+\Lambda_u^{\{k\}}
\right)
P_u
+
2P_y^\top\Lambda_y^{\{k\}}P_y +
\begin{bmatrix}
\mathbf 0
&
\Gamma_{\widetilde\Psi}
+
\Gamma_{\Psi^{\{k\}}}\\
\left(
\Gamma_{\widetilde\Psi}
+
\Gamma_{\Psi^{\{k\}}}
\right)^\top
&
\mathbf 0
\end{bmatrix}.
\end{aligned}
\end{equation*}

The dual problem
\eqref{eqn: dual problem of EOT version of the sinkhorn step}
is therefore equivalent to the first BCD update
\eqref{eq:update-upsi} under the change of variables
\begin{equation*}
\Lambda_u
=
\widetilde\Lambda_u+\Lambda_u^{\{k\}},
\qquad
\Psi
=
\widetilde\Psi+\Psi^{\{k\}}.
\end{equation*}
This follows because
\(\Psi\mapsto\Gamma_\Psi\) is linear, so the expressions in \eqref{eqn: def of tilde D on tilde lambda_u and tilde psi, the dual of the one maringal causal projection} above becomes
\begin{equation*}
\widetilde Z
=
Z(\Lambda_u,\Lambda_y^{\{k\}},\Psi),
\end{equation*}
and
\begin{equation*}
-\operatorname{Tr}(\widetilde\Lambda_u)
=
-\operatorname{Tr}(\Lambda_u)
+
\operatorname{Tr}(\Lambda_u^{\{k\}}).
\end{equation*}
Thus, after the change of variables,
\(\widetilde{\mathcal D}\) differs from
\(\mathcal D(\Lambda_u,\Lambda_y^{\{k\}},\Psi)\) in \eqref{eqn: dual objective} only by a constant
independent of \((\Lambda_u,\Psi)\). 
Consequently, the optimizer of
\eqref{eqn: dual problem of EOT version of the sinkhorn step}
satisfies
\begin{equation*}
\widetilde\Lambda_u+\Lambda_u^{\{k\}}
=
\Lambda_u^{\{k+1\}},
\qquad
\widetilde\Psi+\Psi^{\{k\}}
=
\Psi^{\{k+1\}}.
\end{equation*}
Therefore, using the definition of
\(\theta^{\{k+\frac12\}}\) in Algorithm~\ref{alg:bcd-dual} and \eqref{eqn: primal dual relation for odd sinkhorn projection}, we obtain
\begin{equation*}
\left(
\Sigma^{\pi^{\{k+\frac12\}}}
\right)^{-1}
=
Z(\theta^{\{k+\frac12\}}).
\end{equation*}

We next show that the even Sinkhorn projection
\eqref{eqn: even sinkhorn} corresponds to the second BCD update
\eqref{eq:update-lambday}.
The iterate \(\pi^{\{k+\frac12\}}\) is a centered non--degenerate
Gaussian distribution by
Lemma~\ref{lemma: gaussianity preserved under sinkhorn projections}
and the argument for the odd projection above.
Define a quadratic function \(g^{\{k+\frac12\}}\), up to an additive
constant, by
\begin{equation*}
\frac{d\pi^{\{k+\frac12\}}}{dP}
\propto
\exp\left(
g^{\{k+\frac12\}}-c
\right).
\end{equation*}
The function \(g^{\{k+\frac12\}}\) has no linear term, because
\(\pi^{\{k+\frac12\}}\) and \(P\) are centered and \(c\) has no linear
term. Therefore, we can write
\begin{equation*}
g^{\{k+\frac12\}}(\mathbf z)
=
-\frac12
\mathbf z^\top Q_g^{\{k+\frac12\}}\mathbf z
+
\mathrm{Const}.
\end{equation*}

The even Sinkhorn projection \eqref{eqn: even sinkhorn} can be rewritten,
up to an additive constant independent of \(\pi\), as
\begin{equation*}
\min_{\pi\in\Pi(*,\nu)}
\left\{
\int(c-g^{\{k+\frac12\}})\,d\pi
+
D(\pi,P)
\right\}.
\end{equation*}
This follows from the same calculation used to obtain
\eqref{eqn: odd sinkhorn rewritten as entropy regularized problem}.
The optimizer is centered by the same mean argument used for the odd
projection, because
\begin{equation*}
Q_g^{\{k+\frac12\}}+Q
=
\left(\Sigma^{\pi^{\{k+\frac12\}}}\right)^{-1}
\succ0.
\end{equation*}

Therefore, following the same covariance calculation used for the odd
projection, the dual problem of the even Sinkhorn projection is, up to
an additive constant independent of \(\widetilde\Lambda_y\),
\begin{equation}
\max_{\widetilde\Lambda_y}
\widetilde{\mathcal D}_y(\widetilde\Lambda_y),
\label{eqn: dual problem of even sinkhorn projection}
\end{equation}
where
\begin{equation*}
\widetilde{\mathcal D}_y(\widetilde\Lambda_y)
=
\begin{cases}
\displaystyle
\frac12\log\det\widetilde Z_y
-
\operatorname{Tr}(\widetilde\Lambda_y),
& \widetilde Z_y\succ0,\\[1ex]
-\infty,
& \text{o.w.},
\end{cases}
\end{equation*}
and
\begin{equation*}
\widetilde Z_y
=
Q_g^{\{k+\frac12\}}
+
Q
+
2P_y^\top\widetilde\Lambda_yP_y.
\end{equation*}
Moreover, the unique minimizer of the corresponding Lagrangian over
\(\Sigma\succ0\) satisfies
\begin{equation}
\Sigma=\widetilde Z_y^{-1}.
\label{eqn: primal dual relation for even sinkhorn projection}
\end{equation}

We next show that
\eqref{eqn: dual problem of even sinkhorn projection}
is equivalent to the second BCD update
\eqref{eq:update-lambday}.
The quadratic coefficient \(Q_g^{\{k+\frac12\}}\) satisfies
\begin{equation*}
Q_g^{\{k+\frac12\}}
=
Z(\theta^{\{k+\frac12\}})-Q.
\end{equation*}
This follows from the same comparison of the quadratic coefficients
used for \(Q_g^{\{k\}}\), together with
\begin{equation*}
\left(
\Sigma^{\pi^{\{k+\frac12\}}}
\right)^{-1}
=
Z(\theta^{\{k+\frac12\}}).
\end{equation*}
Hence, substituting this expression into \(\widetilde Z_y\) gives
\begin{equation*}
\begin{aligned}
\widetilde Z_y
&=
Q
+
2P_u^\top\Lambda_u^{\{k+1\}}P_u
+
2P_y^\top
\left(
\widetilde\Lambda_y+\Lambda_y^{\{k\}}
\right)
P_y +
\begin{bmatrix}
\mathbf 0 & \Gamma_{\Psi^{\{k+1\}}}\\
\Gamma_{\Psi^{\{k+1\}}}^\top & \mathbf 0
\end{bmatrix}.
\end{aligned}
\end{equation*}

The dual problem
\eqref{eqn: dual problem of even sinkhorn projection}
is therefore equivalent to the second BCD update
\eqref{eq:update-lambday} under the change of variables
\begin{equation*}
\Lambda_y
=
\widetilde\Lambda_y+\Lambda_y^{\{k\}}.
\end{equation*}
This follows because the change of variables gives
\begin{equation*}
\widetilde Z_y
=
Z(\Lambda_u^{\{k+1\}},\Lambda_y,\Psi^{\{k+1\}}),
\end{equation*}
and
\begin{equation*}
-\operatorname{Tr}(\widetilde\Lambda_y)
=
-\operatorname{Tr}(\Lambda_y)
+
\operatorname{Tr}(\Lambda_y^{\{k\}}).
\end{equation*}
Thus, after the change of variables,
\(\widetilde{\mathcal D}_y\) differs from
\(\mathcal D(\Lambda_u^{\{k+1\}},\Lambda_y,\Psi^{\{k+1\}})\) in
\eqref{eqn: dual objective} only by a constant independent of
\(\Lambda_y\).
Consequently, the optimizer of
\eqref{eqn: dual problem of even sinkhorn projection} satisfies
\begin{equation*}
\widetilde\Lambda_y+\Lambda_y^{\{k\}}
=
\Lambda_y^{\{k+1\}}.
\end{equation*}
Therefore, using
\eqref{eqn: primal dual relation for even sinkhorn projection} and the
definition of \(\theta^{\{k+1\}}\) in
Algorithm~\ref{alg:bcd-dual}, we obtain
\begin{equation*}
\left(
\Sigma^{\pi^{\{k+1\}}}
\right)^{-1}
=
Z(\theta^{\{k+1\}}).
\end{equation*}
The iterate \(\pi^{\{k+1\}}\) is centered and non--degenerate Gaussian
by the same mean argument above and
Lemma~\ref{lemma: gaussianity preserved under sinkhorn projections}.
Therefore, the induction hypothesis holds with \(k\) replaced by
\(k+1\), which completes the proof.
\end{proof}

\section{Supplementary Details for Numerical Experiments}
\subsection{Input Gaussian Process} \label{section: detail of the input gaussina process}

The input distribution $\mu$ is modeled as a $d$-dimensional Gaussian Process over the time grid:
The path-level mean vector ${m}^\mu \in \mathbb{R}^{Td}$ is time-varying. For the $j$-th dimension component at time $t$, the corresponding term in the mean is defined as:
\begin{equation*}
\mathbb{E}^\mu [{U_t}[j]] = \sin\left(t\Delta t + (j - 1)\frac{\pi}{2}\right), \quad j \in \{1, \dots, d\}.
\end{equation*}
Although the resulting marginal $\mu$ is not centered, Lemma \ref{lemma: centering gaussian marginals} ensures that the Gaussian ECOT formulation applies to non-centered marginals as well.
The global covariance matrix ${\Sigma}^\mu \in \mathbb{R}^{Td \times Td}$ is constructed via a Kronecker product:
\begin{equation*}
{\Sigma}^\mu = \mathbf{\Gamma} \otimes {I}_d,
\end{equation*}
where the identity matrix ${I}_d$ is the spatial covariance matrix, assuming zero correlation between different spatial channels, and $\mathbf{\Gamma} \in \mathbb{R}^{T \times T}$ is the temporal correlation matrix governed by the exponential kernel:
\begin{equation*}
\mathbf{\Gamma}_{n, m} = \exp\left( -\frac{1}{2} \frac{|t_n - t_m|}{\ell^2} \right).
\end{equation*}
$\ell = 0.1$ denotes the temporal length-scale. Or equivalently,
\begin{equation*}
\text{Cov}^{\mu}(U(t)_{i}, U(s)_{j})
=
\mathbf{\Gamma}_{t, s} \delta_{i,j},
\quad
\forall i ,j \in \{1,\cdots,d\},
\quad
\forall t,s \in \{1,\cdots,T\}.
\end{equation*}

\subsection{Construction of the cost function}
\label{sec: From a prior dynamic to a prior coupling}
We consider two prior models.
The first is structured as another linear Gaussian State-Space Model (SSM), with model parameters different from those of the true SSM \eqref{eqn: spring system}.
\begin{align*}
X_{t}&=F X_{t-1}+B U_{t}+(\Delta tQ)^{\frac{1}{2}} V_{t} \numberthis \label{eqn: ssm prior}\\
Y_t &=H X_t+(\Delta tR)^{\frac{1}{2}} W_t ,\\
X_0&=\mathbf 0.
\end{align*}
This model \eqref{eqn: ssm prior} represents a standard causal model where information propagates strictly forward in time.

To test our framework against a prior that violates the causal constraint, we introduce a Boundary Value Problem (BVP) prior, mathematically formulated as:
\begin{align*}
X_t  &= F\frac{1}{2}(X_{t-1} +X_{t+1})+B U_t +(\Delta tQ)^{\frac{1}{2}}  V_t \numberthis \label{eqn: bvp prior}\\
Y_t &= H X_t + (\Delta tR)^{\frac{1}{2}} {W}_t,\\
X_0&=X_{T+1}=\mathbf0,
\end{align*}
This BVP prior \eqref{eqn: bvp prior} conditions the latent state $X_t$ simultaneously on its immediate past $X_{t-1}$ and its immediate future $X_{t+1}$. This allows us to evaluate whether our proposed framework can successfully obtain causal predictive structures from a prior distribution that inherently encodes non-causal, global dependencies.

Both priors use i.i.d.\ standard Gaussian inputs
and independent Gaussian process and observation noises,
with \(Q\succeq0\) and \(R\succ0\).
The noise variables $\{V_t,W_t\}$ in the prior models are understood
to be defined under the respective prior distributions,
independently of the true SSM.
For the BVP prior, we additionally assume that all eigenvalues of $F$ have absolute value strictly less than $1$. 

We next show how these prior dynamics define Gaussian couplings $\gamma$ with the covariance $\Sigma^\gamma$, from which the corresponding quadratic costs in the form of \eqref{eqn: quad c with pd hessian} are obtained, and how the assumptions here ensure that $\Sigma^\gamma\succ0$, as required by Assumption \ref{assumption:cost}.

The dynamics in \eqref{eqn: ssm prior} and
\eqref{eqn: bvp prior} can both be written in matrix form as
\begin{align*}
\mathbf M_\gamma\mathbf X
&=\mathbf B_\gamma\mathbf U
+\sqrt{\Delta t}\mathbf Q_\gamma^{1/2}\mathbf V,\\
\mathbf Y
&=\mathbf H_\gamma\mathbf X
+\sqrt{\Delta t}\mathbf R_\gamma^{1/2}\mathbf W, \numberthis\label{eqn: matrix form of the prior system for both}
\end{align*}
where $\mathbf X$, $\mathbf V$, and $\mathbf W$ denote
the corresponding stacked vectors, and
\begin{equation*}
\mathbf B_\gamma
=\operatorname{diag}(B,\ldots,B)
\in\mathbb R^{Td_x\times Td},
\end{equation*}
\begin{equation*}
\mathbf H_\gamma
=\operatorname{diag}(H,\ldots,H)
\in\mathbb R^{Td\times Td_x},
\end{equation*}
\begin{equation*}
\mathbf Q_\gamma
=\operatorname{diag}(Q,\ldots,Q)
\in\mathbb R^{Td_x\times Td_x},
\qquad
\mathbf R_\gamma
=\operatorname{diag}(R,\ldots,R)
\in\mathbb R^{Td\times Td}.
\end{equation*}
The two models differ in the definition of
$\mathbf M_\gamma$, as specified below.
\begin{enumerate}
    \item \textbf{SSM prior.}
For the dynamics in \eqref{eqn: ssm prior},
\begin{equation*}
\mathbf M_\gamma=
\begin{bmatrix}
I_{d_x} & 0 & 0 & \cdots & 0\\
-F & I_{d_x} & 0 & \cdots & 0\\
0 & -F & I_{d_x} & \ddots & \vdots\\
\vdots & \ddots & \ddots & \ddots & 0\\
0 & \cdots & 0 & -F & I_{d_x}
\end{bmatrix}
\in\mathbb R^{Td_x\times Td_x}.
\end{equation*}
Since $\mathbf M_\gamma$ is lower triangular with identity
diagonal blocks, it is invertible.
    \item \textbf{BVP prior.}
For the dynamics in \eqref{eqn: bvp prior},
\begin{equation*}
\mathbf M_\gamma=
\begin{bmatrix}
I_{d_x} & -F/2 & 0 & \cdots\\
-F/2 & I_{d_x} & -F/2 & \cdots\\
0 & -F/2 & I_{d_x} & \cdots\\
\vdots & \vdots & \vdots & \ddots
\end{bmatrix}
\in\mathbb R^{Td_x\times Td_x}.
\end{equation*}
We verify that $\mathbf M_\gamma$ is invertible.
Write
$\mathbf M_\gamma=I_{Td_x}-J_T\otimes F,$
where $J_T\in\mathbb R^{T\times T}$ is the symmetric
tridiagonal matrix with zero diagonal entries and
$1/2$ on the two adjacent diagonals:
\begin{equation*}
J_T=
\frac12
\begin{bmatrix}
0&1&0&\cdots\\
1&0&1&\cdots\\
0&1&0&\cdots\\
\vdots&\vdots&\vdots&\ddots
\end{bmatrix}
\in\mathbb R^{T\times T}.
\end{equation*}
Every eigenvalue $\lambda$ of $J_T$
satisfies $|\lambda|\leq1$. To see why, suppose
$J_Tv=\lambda v$ with $v\in\mathbb R^T$, $v\neq0$,
and choose $j$ such that
$|v_j|=\max_{1\leq t\leq T}|v_t|$.
With the convention $v_0=v_{T+1}=0$, we obtain
\begin{equation*}
|\lambda||v_j|
=
\frac12|v_{j-1}+v_{j+1}|
\leq
\frac12(|v_{j-1}|+|v_{j+1}|)
\leq |v_j|.
\end{equation*}
Thus, $|\lambda|\leq1$.
Since the eigenvalues of $J_T\otimes F$ are products
of eigenvalues of $J_T$ and $F$, and all eigenvalues
of $F$ have absolute value strictly less than $1$,
all eigenvalues of $J_T\otimes F$ have absolute value
strictly less than $1$.
Therefore, $\mathbf M_\gamma=I_{Td_x}-J_T\otimes F$
is invertible.
\end{enumerate}

For either model, the dynamics in \eqref{eqn: matrix form of the prior system for both}, together with a Gaussian
distribution on $U_{1:T}$ with covariance $I_{Td}$,
define a Gaussian coupling $\gamma$ with covariance
\begin{align*}
\operatorname{Cov}^\gamma(\mathbf U,\mathbf U)
&=I_{Td},\\
\operatorname{Cov}^\gamma(\mathbf U,\mathbf Y)
&=(\mathbf M_\gamma^{-1}\mathbf B_\gamma)^\top
\mathbf H_\gamma^\top,\\
\operatorname{Cov}^\gamma(\mathbf Y,\mathbf Y)
&=\mathbf H_\gamma\mathbf M_\gamma^{-1}
\left(
\mathbf B_\gamma\mathbf B_\gamma^\top
+\Delta t\mathbf Q_\gamma
\right)
(\mathbf M_\gamma^{-1})^\top\mathbf H_\gamma^\top
+\Delta t\mathbf R_\gamma.
\end{align*}
To verify these covariance expressions, substituting
the first matrix equation into the second gives
\begin{align*}
\mathbf Y
={}&\mathbf H_\gamma\mathbf M_\gamma^{-1}
\mathbf B_\gamma\mathbf U
+\sqrt{\Delta t}\mathbf H_\gamma\mathbf M_\gamma^{-1}
\mathbf Q_\gamma^{1/2}\mathbf V
+\sqrt{\Delta t}\mathbf R_\gamma^{1/2}\mathbf W.
\end{align*}
Since $\mathbf U$, $\mathbf V$, and $\mathbf W$ are
mutually independent, with
$\operatorname{Cov}(\mathbf U)=I_{Td}$,
$\operatorname{Cov}(\mathbf V)=I_{Td_x}$, and
$\operatorname{Cov}(\mathbf W)=I_{Td}$,
the covariance expressions above follow directly.

To see that $\operatorname{Cov}^\gamma\succ0$, let
\begin{equation*}
\mathcal C=\mathbf H_\gamma\mathbf M_\gamma^{-1}\mathbf B_\gamma,
\qquad
\mathcal S=\Delta t\left(
\mathbf H_\gamma\mathbf M_\gamma^{-1}
\mathbf Q_\gamma(\mathbf M_\gamma^{-1})^\top
\mathbf H_\gamma^\top+\mathbf R_\gamma
\right).
\end{equation*}
Since $Q\succeq0$, we have $\mathbf Q_\gamma\succeq0$.
For any $y\in\mathbb R^{Td}$,
\begin{align*}
&y^\top\mathbf H_\gamma\mathbf M_\gamma^{-1}
\mathbf Q_\gamma(\mathbf M_\gamma^{-1})^\top
\mathbf H_\gamma^\top y\\
&\quad=
\left((\mathbf M_\gamma^{-1})^\top\mathbf H_\gamma^\top y\right)^\top
\mathbf Q_\gamma
\left((\mathbf M_\gamma^{-1})^\top\mathbf H_\gamma^\top y\right)
\geq0.
\end{align*}
Therefore,
\begin{equation*}
\mathbf H_\gamma\mathbf M_\gamma^{-1}
\mathbf Q_\gamma(\mathbf M_\gamma^{-1})^\top
\mathbf H_\gamma^\top\succeq0.
\end{equation*}
Since $R\succ0$ and $\Delta t>0$, it follows that
\begin{equation*}
\mathcal S\succeq\Delta t\,\mathbf R_\gamma\succ0.
\end{equation*}
Moreover,
\begin{equation*}
\operatorname{Cov}^\gamma
=
\begin{bmatrix}
I_{Td}&\mathcal C^\top\\
\mathcal C&\mathcal C\mathcal C^\top+\mathcal S
\end{bmatrix}.
\end{equation*}
For any $u,y\in\mathbb R^{Td}$ with $(u,y)\neq(0,0)$,
\begin{equation*}
\begin{bmatrix}u\\y\end{bmatrix}^\top
\operatorname{Cov}^\gamma
\begin{bmatrix}u\\y\end{bmatrix}
=
\|u+\mathcal C^\top y\|^2+y^\top\mathcal S y>0.
\end{equation*}
Thus, $\operatorname{Cov}^\gamma\succ0$.
The cost in \eqref{eqn: quad c with pd hessian} can be obtained
from $(\operatorname{Cov}^\gamma)^{-1}$.

\subsection{Model Parameters}\label{section: numercial experiment parameters}
{
For the true SSM in \eqref{eqn: spring system}, the model parameters in our experiments are
\begin{align*}
F_* &= \begin{bmatrix}0.65&-0.12\\0.12&0.65\end{bmatrix}\otimes I_d,
& B_* &= \begin{bmatrix}0.8\\0.8\end{bmatrix}\otimes I_d, \\
H_* &= \begin{bmatrix}1&1\end{bmatrix}\otimes I_d,
& Q_* &= 0.04 I_{2d}, \qquad R_*=0.03 I_d.
\end{align*}
For both priors (\eqref{eqn: ssm prior}, \eqref{eqn: bvp prior}), we use
\begin{align*}
F &= \begin{bmatrix}0.45&-0.04\\0.04&0.45\end{bmatrix}\otimes I_d,
& B &= \begin{bmatrix}0.60\\0.60\end{bmatrix}\otimes I_d, \\
H &= \begin{bmatrix}0.80&0.80\end{bmatrix}\otimes I_d,
& Q &= 0.12 I_{2d}, \qquad R=0.11 I_d.
\end{align*}
}

\subsection{Moments of the True State Space Model}
\label{section: for the simpler true ssm experiment-the ref coupling for stability numeical plot}
We derive the output marginal moments induced by the true SSM
in \eqref{eqn: spring system}, which are used to construct
the reference coupling in the covariance stability experiment. 
Recall that
\begin{equation*}
U_{1:T}
=
\begin{bmatrix}
U_1^\top & \cdots & U_T^\top
\end{bmatrix}^{\top}
\sim
\mathcal N(m^\mu,\Sigma^\mu),
\end{equation*}
where
\begin{equation*}
\Sigma^\mu=\Gamma\otimes I_d.
\end{equation*}

The dynamics in \eqref{eqn: spring system} can be written
in matrix form as
\begin{align*}
\mathbf M_*\mathbf X
&=\mathbf B_*\mathbf U
+\sqrt{\Delta t}\mathbf Q_*^{1/2}\mathbf V,\\
\mathbf Y
&=\mathbf H_*\mathbf X
+\sqrt{\Delta t}\mathbf R_*^{1/2}\mathbf W,
\end{align*}
where $\mathbf X$, $\mathbf U$, $\mathbf Y$, $\mathbf V$,
and $\mathbf W$ denote the corresponding stacked vectors, and
\begin{equation*}
\mathbf M_*=
\begin{bmatrix}
I_{d_x} & 0 & 0 & \cdots & 0\\
-F_* & I_{d_x} & 0 & \cdots & 0\\
0 & -F_* & I_{d_x} & \ddots & \vdots\\
\vdots & \ddots & \ddots & \ddots & 0\\
0 & \cdots & 0 & -F_* & I_{d_x}
\end{bmatrix}
\in\mathbb R^{Td_x\times Td_x}.
\end{equation*}
\begin{equation*}
\mathbf B_*
=\operatorname{diag}(B_*,\ldots,B_*)
\in\mathbb R^{Td_x\times Td},
\end{equation*}
\begin{equation*}
\mathbf H_*
=\operatorname{diag}(H_*,\ldots,H_*)
\in\mathbb R^{Td\times Td_x},
\end{equation*}
\begin{equation*}
\mathbf Q_*
=\operatorname{diag}(Q_*,\ldots,Q_*)
\in\mathbb R^{Td_x\times Td_x},
\qquad
\mathbf R_*
=\operatorname{diag}(R_*,\ldots,R_*)
\in\mathbb R^{Td\times Td}.
\end{equation*}
Since $\mathbf M_*$ is lower triangular with identity
diagonal blocks, it is invertible.

Substituting the first matrix equation into the second gives
\begin{equation}
\mathbf Y
=
\mathbf H_*\mathbf M_*^{-1}\mathbf B_*\mathbf U
+\eta,
\label{eqn: stacked true ssm output}
\end{equation}
where
\begin{equation*}
\eta
=
\sqrt{\Delta t}\mathbf H_*\mathbf M_*^{-1}
\mathbf Q_*^{1/2}\mathbf V
+\sqrt{\Delta t}\mathbf R_*^{1/2}\mathbf W.
\end{equation*}
Since $\mathbf V$ and $\mathbf W$ are independent standard
Gaussian vectors and are independent of $\mathbf U$,
the vector $\eta$ is centered Gaussian and independent of
$\mathbf U$, with covariance
\begin{equation}
\Sigma^\eta
=
\Delta t\left(
\mathbf H_*\mathbf M_*^{-1}
\mathbf Q_*(\mathbf M_*^{-1})^\top
\mathbf H_*^\top
+\mathbf R_*
\right).
\label{eqn: covariance of true ssm noise}
\end{equation}

Since $\mathbf U$ and $\eta$ are independent Gaussian vectors,
\eqref{eqn: stacked true ssm output} shows that $\mathbf Y$
is Gaussian with
\begin{equation*}
\nu=\mathcal N(m^\nu,\Sigma^\nu),
\end{equation*}
where
\begin{equation}
m^\nu
=
\mathbf H_*\mathbf M_*^{-1}\mathbf B_*m^\mu,
\label{eqn: ref mean}
\end{equation}
and
\begin{equation}
\Sigma^\nu
=
\mathbf H_*\mathbf M_*^{-1}\mathbf B_*
\Sigma^\mu
\mathbf B_*^\top(\mathbf M_*^{-1})^\top
\mathbf H_*^\top
+\Sigma^\eta.
\label{eqn: nu cov from true paramter and Sigma mu and Sigma eta}
\end{equation}

Finally, we verify that $\Sigma^\nu\succ0$.
First, since $Q_*\succeq0$, we have $\mathbf Q_*\succeq0$.
For any $y\in\mathbb R^{Td}$,
\begin{align*}
&y^\top\mathbf H_*\mathbf M_*^{-1}
\mathbf Q_*(\mathbf M_*^{-1})^\top
\mathbf H_*^\top y\\
&\quad=
\left((\mathbf M_*^{-1})^\top\mathbf H_*^\top y\right)^\top
\mathbf Q_*
\left((\mathbf M_*^{-1})^\top\mathbf H_*^\top y\right)
\geq0.
\end{align*}
Therefore,
\begin{equation*}
\mathbf H_*\mathbf M_*^{-1}
\mathbf Q_*(\mathbf M_*^{-1})^\top
\mathbf H_*^\top\succeq0.
\end{equation*}
Since $R_*\succ0$ and $\Delta t>0$, it follows from \eqref{eqn: covariance of true ssm noise} that
\begin{equation*}
\Sigma^\eta\succeq\Delta t\,\mathbf R_*\succ0.
\end{equation*}
Since $\Sigma^\mu\succeq0$, for any $y\in\mathbb R^{Td}$,
\begin{align*}
&y^\top\mathbf H_*\mathbf M_*^{-1}\mathbf B_*
\Sigma^\mu
\mathbf B_*^\top(\mathbf M_*^{-1})^\top
\mathbf H_*^\top y\\
&\quad=
\left(\mathbf B_*^\top(\mathbf M_*^{-1})^\top
\mathbf H_*^\top y\right)^\top
\Sigma^\mu
\left(\mathbf B_*^\top(\mathbf M_*^{-1})^\top
\mathbf H_*^\top y\right)
\geq0.
\end{align*}
Therefore,
\begin{equation*}
\mathbf H_*\mathbf M_*^{-1}\mathbf B_*
\Sigma^\mu
\mathbf B_*^\top(\mathbf M_*^{-1})^\top
\mathbf H_*^\top\succeq0.
\end{equation*}
It follows from \eqref{eqn: nu cov from true paramter and Sigma mu and Sigma eta} that
\begin{equation*}
\Sigma^\nu\succeq\Sigma^\eta\succ0.
\end{equation*}
Thus, the output marginal $\nu$ is a non-degenerate Gaussian distribution.

\subsection{Prefix Non--causal Coupling}
\label{section: definition of the prefix non--causal}
Consider the centered Gaussian marginals and reference
\begin{equation*}
\mu=\mathcal N(0,\Sigma^\mu),\qquad
\nu=\mathcal N(0,\Sigma^\nu),\qquad
\gamma=\mathcal N(0,\Sigma^\gamma).
\end{equation*}
For each \(t=1,\ldots,T\), let
\(P_t\in\mathbb R^{2td\times 2Td}\) be the coordinate-selection
matrix satisfying
\begin{equation*}
P_t
\begin{bmatrix}
\mathbf U_{1:T}\\
\mathbf Y_{1:T}
\end{bmatrix}
=
\begin{bmatrix}
\mathbf U_{1:t}\\
\mathbf Y_{1:t}
\end{bmatrix}.
\end{equation*}
Define the prefix covariance matrices
\begin{equation*}
\Sigma_t^{\gamma,\mathrm{prefix}}
=
P_t\Sigma^\gamma P_t^\top,
\qquad
\Sigma_t^{\mu,\mathrm{prefix}}
=
[\Sigma^\mu]_{1:td,1:td},
\qquad
\Sigma_t^{\nu,\mathrm{prefix}}
=
[\Sigma^\nu]_{1:td,1:td}.
\end{equation*}
Thus, the prefix marginals and nominal reference are
\begin{equation*}
\mu_t^{\mathrm{prefix}}
=
\mathcal N(0,\Sigma_t^{\mu,\mathrm{prefix}}),
\qquad
\nu_t^{\mathrm{prefix}}
=
\mathcal N(0,\Sigma_t^{\nu,\mathrm{prefix}}),
\qquad
\gamma_t^{\mathrm{prefix}}
=
\mathcal N(0,\Sigma_t^{\gamma,\mathrm{prefix}}).
\end{equation*}

Let
\begin{equation*}
\Sigma_t^{\mu,\mathrm{prefix}}
=
L_t^\mu(L_t^\mu)^\top,
\qquad
\Sigma_t^{\nu,\mathrm{prefix}}
=
L_t^\nu(L_t^\nu)^\top
\end{equation*}
be their Cholesky decompositions, and define
\begin{equation*}
L_t=\operatorname{diag}(L_t^\mu,L_t^\nu).
\end{equation*}
Under the transformation
\begin{equation*}
\begin{bmatrix}
\widetilde{\mathbf U}_{1:t}\\
\widetilde{\mathbf Y}_{1:t}
\end{bmatrix}
=
L_t^{-1}
\begin{bmatrix}
\mathbf U_{1:t}\\
\mathbf Y_{1:t}
\end{bmatrix},
\end{equation*}
the marginal distributions become
\begin{equation*}
\widetilde\mu_t^{\mathrm{prefix}}
=
\widetilde\nu_t^{\mathrm{prefix}}
=
\mathcal N(0,I_{td}),
\end{equation*}
and the nominal reference becomes
\begin{equation*}
\widetilde\gamma_t^{\mathrm{prefix}}
=
\mathcal N(0,\widetilde\Sigma_t^{\gamma,\mathrm{prefix}}),
\qquad
\widetilde\Sigma_t^{\gamma,\mathrm{prefix}}
=
L_t^{-1}\Sigma_t^{\gamma,\mathrm{prefix}}L_t^{-\top}.
\end{equation*}

For each \(t=1,\ldots,T\), we independently solve
the non--causal Gaussian ECOT problem
\begin{equation*}
\widetilde\pi_t^{*,\mathrm{prefix}}
=
\underset{
\widetilde\pi\in
\Pi(\widetilde\mu_t^{\mathrm{prefix}},
\widetilde\nu_t^{\mathrm{prefix}})
}{
\operatorname{argmin}
}
D_{\mathrm{KL}}
(\widetilde\pi\|\widetilde\gamma_t^{\mathrm{prefix}}).
\end{equation*}
No causality constraint is imposed.

By the Gaussian characterization, the optimizer is Gaussian,
\begin{equation*}
\widetilde\pi_t^{*,\mathrm{prefix}}
=
\mathcal N(0,\widetilde\Sigma_t^{*,\mathrm{prefix}}).
\end{equation*}
Transforming back to the original coordinates gives
\begin{equation*}
\pi_t^{*,\mathrm{prefix}}
=
\mathcal N(0,\Sigma_t^{*,\mathrm{prefix}}),
\qquad
\Sigma_t^{*,\mathrm{prefix}}
=
L_t\widetilde\Sigma_t^{*,\mathrm{prefix}}L_t^\top.
\end{equation*}
At \(t=T\), the optimization coincides with the full-horizon non--causal ECOT problem.

\subsection{Prediction Tasks}\label{eqn: section: prediction task setup}
As described in Section \ref{sec: numerical experiment}, the model observes only $U_{1:t}$ to
predict $Y_t$, namely  input-only prediction. We consider the Gaussian conditional model
\begin{equation*}
Y_t
=
\widetilde{\mathbf h}_t^\top U_{1:t}
+\widetilde b_t+\widetilde w_t,
\qquad
\widetilde w_t
\sim
\mathcal N(0,\widetilde\Sigma_t^{\mathrm{cond}}).
\end{equation*}
By conditioning $Y_t$ on $U_{1:t}$ under $\pi$, we obtain
\begin{align*}
\widetilde{\mathbf h}_t^\top
&=
\mathrm{Cov}^{\pi}(Y_t,U_{1:t})
\bigl[\mathrm{Cov}^{\pi}(U_{1:t},U_{1:t})\bigr]^{-1},\\
\widetilde b_t
&=
\mathbb E^{\pi}[Y_t]
-
\widetilde{\mathbf h}_t^\top
\mathbb E^{\pi}[U_{1:t}],\\
\widetilde\Sigma_t^{\mathrm{cond}}
&=
\mathrm{Cov}^{\pi}(Y_t,Y_t)
-
\mathrm{Cov}^{\pi}(Y_t,U_{1:t})
\bigl[\mathrm{Cov}^{\pi}(U_{1:t},U_{1:t})\bigr]^{-1}
\mathrm{Cov}^{\pi}(U_{1:t},Y_t).
\end{align*}

In addition to this input-only prediction task, we also consider one-step-ahead prediction, in which the model
observes $U_{1:t}$ and $Y_{1:t-1}$ to predict $Y_t$.
Specifically, for each coupling, we recover a linear Gaussian
input-output model
\begin{equation*}
Y_t = {\mathbf h_t}^{\top} U_{1:t}
+ {\mathbf f_t}^{\top} Y_{1:t-1}+b_t+w_t,
\qquad
w_t\sim \mathcal{N}(0,\Sigma_t^{\mathrm{cond}}),
\end{equation*}
by conditioning $Y_t$ on
$\text{past}_t=(U_{1:t},Y_{1:t-1})$:
\begin{align*}
&[\mathbf{h}_t^\top, \mathbf{f}_t^\top]
= \text{Cov}^{\pi}({Y_t,\text{past}_t})
  {\text{Cov}^{\pi}({\text{past}_t,\text{past}_t})}^{-1},\\
&b_t = \mathbb{E}^{\pi}[Y_t]
- \left(
\mathbf{h}_t^\top \mathbb{E}^{\pi}[U_{1:t}]
+ \mathbf{f}_t^\top \mathbb{E}^{\pi}[Y_{1:t-1}]
\right),\\
&\Sigma_t^{\mathrm{cond}}
= \text{Cov}^{\pi}(Y_t,Y_t)
-\text{Cov}^{\pi}(Y_t,\text{past}_t)
{\text{Cov}^{\pi}(\text{past}_t,\text{past}_t)}^{-1}
\text{Cov}^{\pi}(\text{past}_t,Y_t).
\end{align*}
Here we use the notation $\text{Cov}^\pi(u,v)$ for the covariance
matrix between $u$ and $v$ under the probability measure $\pi$.
We evaluate the same causal, non-causal, and prefix non-causal
couplings under the SSM and BVP priors, using the same experimental
setup and RMSE definition as in Section~\ref{sec: numerical experiment}.
The results are in Section \ref{extra section: numerical results} below.

\subsection{Additional Numerical Results}\label{extra section: numerical results}

The results of the one-step-ahead prediction exhibit trends similar to those observed in the input-only prediction setting. Under the BVP prior, the causal
coupling consistently achieves lower RMSE than both non-causal alternatives. Under the SSM prior, the causal coupling becomes more competitive as the marginal estimation sample size increases, eventually outperforming both non-causal methods.

\begin{figure}[!tbp]
    \centering
    \includegraphics[width=\textwidth]{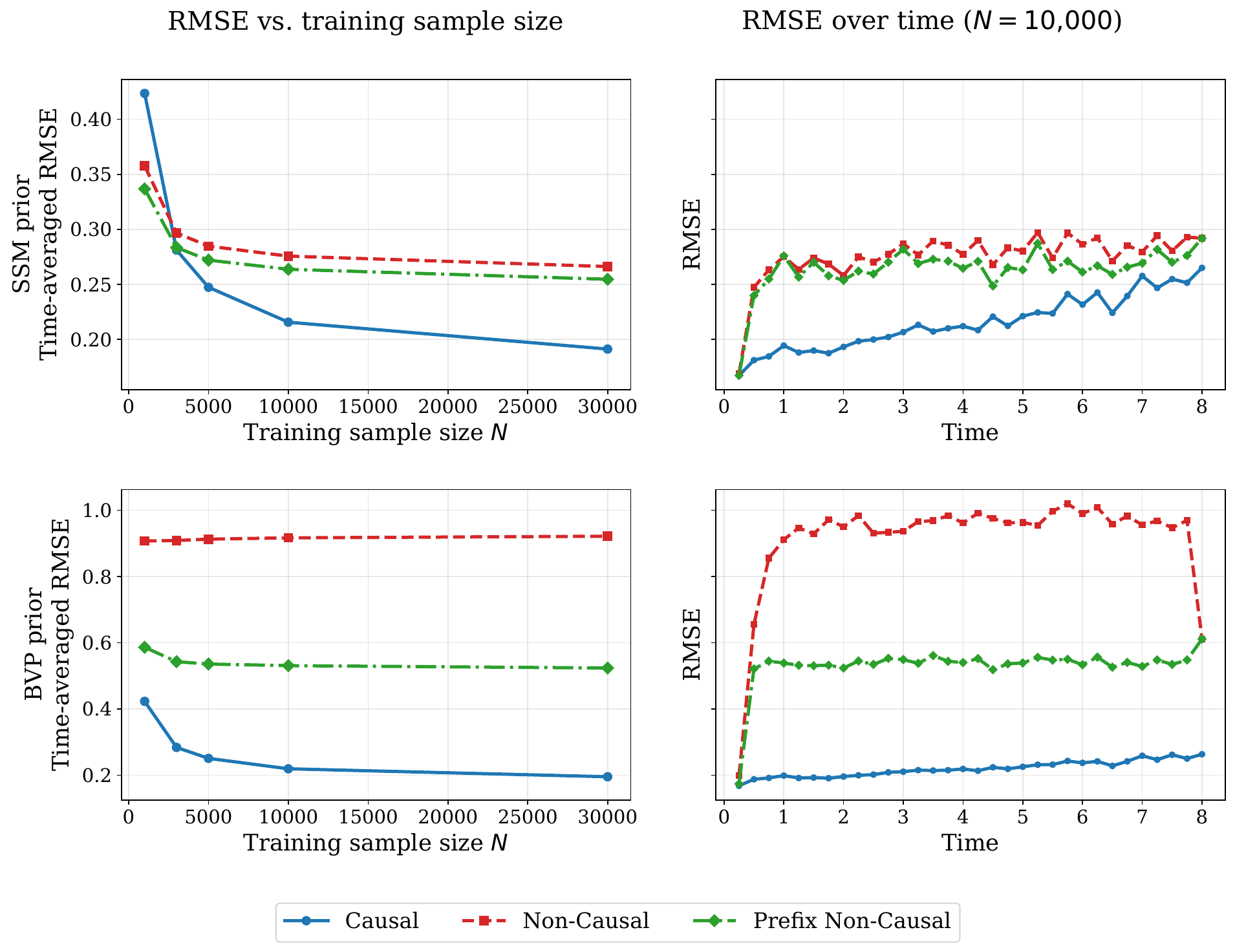}
    \caption{Prediction RMSE for the causal, non-causal, and prefix
    non-causal couplings under one-step-ahead prediction.
    The top row uses the SSM prior and the bottom row uses the BVP prior.
    Left: Time-averaged RMSE against the marginal estimation sample size $N$.
    Right: RMSE over time with $N=10000$.
    The total time horizon is 8, with $T=32$ time steps and
    data dimension $d=2$.}
    \label{fig:simulation_one_step_ahead}
\end{figure}

The corresponding numerical values of
Figure~\ref{fig:simulation_input_only} and Figure \ref{fig:simulation_one_step_ahead} are 
reported below.
\begin{table}[t]
\centering
\small
\caption{Time-averaged RMSE for one-step-ahead prediction
under different marginal estimation sample sizes $N$.}
\label{tab:one_step_ahead_rmse_results}
\begin{tabular}{lccccc}
\toprule
Prior & $N$ & Causal & Non-Causal & Prefix Non-Causal & $\gamma$ \\
\midrule
SSM & 1,000  & 0.4234 & 0.3589 & 0.3354 & 1.3490 \\
SSM & 3,000  & 0.2841 & 0.2972 & 0.2844 & 1.3490 \\
SSM & 5,000  & 0.2481 & 0.2844 & 0.2727 & 1.3490 \\
SSM & 10,000 & 0.2152 & 0.2759 & 0.2635 & 1.3490 \\
SSM & 30,000 & 0.1911 & 0.2664 & 0.2546 & 1.3490 \\
\midrule
BVP & 1,000  & 0.4225 & 0.9075 & 0.5831 & 1.3797 \\
BVP & 3,000  & 0.2855 & 0.9068 & 0.5431 & 1.3797 \\
BVP & 5,000  & 0.2511 & 0.9114 & 0.5357 & 1.3797 \\
BVP & 10,000 & 0.2190 & 0.9149 & 0.5306 & 1.3797 \\
BVP & 30,000 & 0.1952 & 0.9208 & 0.5238 & 1.3797 \\
\bottomrule
\end{tabular}
\end{table}
\begin{table}[t]
\centering
\small
\caption{Time-averaged RMSE for input-only prediction
under different marginal estimation sample sizes $N$.}
\label{tab:input_only_rmse_results}
\begin{tabular}{lccccc}
\toprule
Prior & $N$ & Causal & Non-Causal & Prefix Non-Causal & $\gamma$ \\
\midrule
SSM & 1,000  & 0.5138 & 0.4030 & 0.4014 & 1.8303 \\
SSM & 3,000  & 0.3431 & 0.3549 & 0.3510 & 1.8303 \\
SSM & 5,000  & 0.2959 & 0.3419 & 0.3363 & 1.8303 \\
SSM & 10,000 & 0.2553 & 0.3323 & 0.3253 & 1.8303 \\
SSM & 30,000 & 0.2260 & 0.3236 & 0.3167 & 1.8303 \\
\midrule
BVP & 1,000  & 0.5018 & 0.7243 & 0.6941 & 2.0768 \\
BVP & 3,000  & 0.3423 & 0.7000 & 0.6528 & 2.0768 \\
BVP & 5,000  & 0.2967 & 0.6937 & 0.6399 & 2.0768 \\
BVP & 10,000 & 0.2578 & 0.6883 & 0.6297 & 2.0768 \\
BVP & 30,000 & 0.2284 & 0.6819 & 0.6210 & 2.0768 \\
\bottomrule
\end{tabular}
\end{table}

\section{A Gaussian ECOT problem with a known optimizer}\label{section: a problem with known optimizer}
In the Cholesky coordinates introduced in
Section~\ref{section: cholesky Coordinates}, let
$\pi^{\mathrm{truth}}=\mathcal N(0,\Sigma^{\mathrm{truth}})
\in\Pi_c(\mu,\nu)$
be a non--degenerate causal Gaussian coupling, with the covariance
$\Sigma^{\mathrm{truth}}\succ0$.
For $t=2,\ldots,T$, choose matrices
$A_t\in\mathbb R^{d\times d}$ and define
\begin{equation}
S(\mathbf u,\mathbf y)
:=
\sum_{t=2}^T u_t^\top A_t y_{t-1}.
\label{eq: manufactured causal test function}
\end{equation}
Since $\mu=\mathcal N(0,I_{Td})$, the random variables
$U_t$ and $U_{1:t-1}$ are independent and the marginal means are zero. Consequently,
for every $\pi\in\Pi_c(\mu,\nu)$,
\begin{equation}
\begin{aligned}
\mathbb E_\pi[S(U,Y)]
&=
\sum_{t=2}^T
\mathbb E_\pi\left[
Y_{t-1}^\top A_t^\top
\mathbb E_\pi[U_t\mid U_{1:t-1},Y_{1:t-1}]
\right]\\
&=0.
\end{aligned}
\label{eq: manufactured zero expectation}
\end{equation}

Define the prior $\gamma$ by
\begin{equation}
d\gamma
:=
\frac{e^{-S}}{Z}\,d\pi^{\mathrm{truth}},
\qquad
Z:=\mathbb E_{\pi^{\mathrm{truth}}}[e^{-S(U,Y)}].
\label{eq: manufactured gaussian prior}
\end{equation}
\begin{lemma}
For the prior $\gamma$ defined in
\eqref{eq: manufactured gaussian prior}, if it has a positive definite covariance matrix, then  
$\pi^{\mathrm{truth}}$ is the unique optimizer of
\begin{equation}
\min_{\pi\in\Pi_c(\mu,\nu)}
D(\pi,\gamma).
\label{eq: manufactured exact KL problem}
\end{equation}
Moreover, the optimal value is $\log Z$.
\end{lemma}

\begin{proof}
For any $\pi\in\Pi_c(\mu,\nu)$ with
$D(\pi,\gamma)<\infty $, we have
\begin{equation}
\begin{aligned}
D(\pi,\gamma)
&=
D(\pi,\pi^{\mathrm{truth}})
+\mathbb E_\pi[S(U,Y)]
+\log Z\\
&=
D(\pi,\pi^{\mathrm{truth}})+\log Z\\
&\geq \log Z,
\end{aligned}
\end{equation}
where the second equality follows from
\eqref{eq: manufactured zero expectation}.
Equality holds if and only if $\pi=\pi^{\mathrm{truth}}$.
Since $\pi^{\mathrm{truth}}$ is feasible, it is the unique optimizer,
and the optimal value is $\log Z$.
\end{proof}

We next express the construction in terms of Gaussian
covariance matrices. 
Let $H\in\mathbb S^{2Td}$ have zero diagonal blocks, with nonzero off-diagonal blocks 
\begin{equation*}
H_{U_t,Y_{t-1}}=A_t,
\qquad
H_{Y_{t-1},U_t}=A_t^\top,
\qquad t=2,\ldots,T,
\end{equation*}
so that we have $S(\mathbf u,\mathbf y)
=
\frac12
\begin{bmatrix}
\mathbf u\\
\mathbf y
\end{bmatrix}^{\top}
H
\begin{bmatrix}
\mathbf u\\
\mathbf y
\end{bmatrix}$.
Assume that $(\Sigma^{\mathrm{truth}})^{-1}+H\succ0.$

Then $\gamma$ is a centered non--degenerate Gaussian
distribution with covariance $\Sigma^\gamma
=
\left((\Sigma^{\mathrm{truth}})^{-1}+H\right)^{-1}.$
The cost associated with $\gamma$
takes the form of Assumption~\ref{assumption:cost},
and the objective in \eqref{eq:problem-formulation-2}
differs from $D(\pi,\gamma)$ only by a constant
independent of $\pi$.
Therefore, $\pi^{\mathrm{truth}}$ is also the unique optimizer
of \eqref{eq:problem-formulation-2}.
In particular, the exact optimizer of the covariance
problem \eqref{eqn:problem-over-cov} is
$\Sigma^*=\Sigma^{\mathrm{truth}}.$

The corresponding dual optimizer is also known explicitly:
\begin{equation*}
\Lambda_u^*=\Lambda_y^*=-\frac12 I_{Td},
\qquad
\Psi_t^*
=
\begin{bmatrix}
0_{(t-1)d\times d}\\
-A_{t+1}^\top
\end{bmatrix},
\quad t=1,\ldots,T-1.
\end{equation*}

\end{document}